\documentclass[a4paper,reqno,12pt]{amsart}
\usepackage{mathrsfs}
\usepackage{stmaryrd}
\usepackage{cases}
\usepackage{epic}
\usepackage{amsfonts}
\usepackage{graphicx}
\usepackage{amsmath}
\usepackage{amssymb, upgreek}
\usepackage{mathtools}
\usepackage{bm}
\usepackage{array}
\usepackage{chngcntr} 
\counterwithin{figure}{section}
\counterwithin{table}{section}
\usepackage{latexsym,todonotes}
\usepackage{pdflscape}
\usepackage[all]{xypic}
\usepackage[all]{xy}
\usepackage{color}
\usepackage{colordvi}
\usepackage{multicol}
\usepackage{tikz}
\usepackage{tikz-cd}
\usepackage[linktocpage=true]{hyperref}
\hypersetup{colorlinks,linkcolor=blue,urlcolor=cyan,citecolor=blue}
\newcommand{\bvs}{\mathbf{\varsigma}}
\newcommand{\vs}{\varsigma}
\allowdisplaybreaks

\DeclareMathAlphabet{\mathbbb}{U}{bbold}{m}{n}

\begin{document}
	\input xy
	\xyoption{all}

	\newtheorem{innercustomthm}{{\bf Main Theorem}}
	\newenvironment{customthm}[1]
	{\renewcommand\theinnercustomthm{#1}\innercustomthm}
	{\endinnercustomthm}
	
	\newtheorem{innercustomcor}{{\bf Corollary}}
	\newenvironment{customcor}[1]
	{\renewcommand\theinnercustomcor{#1}\innercustomcor}
	{\endinnercustomthm}
	
	\newtheorem{innercustomprop}{{\bf Proposition}}
	\newenvironment{customprop}[1]
	{\renewcommand\theinnercustomprop{#1}\innercustomprop}
	{\endinnercustomthm}
	\def \wt{\mathrm{wt}}
	\def \wItau{\I_{\circ,\tau}}
	\def \brW{\mathrm{Br}(W_\btau)}
	\def \br{\mathbf{r}}
	\newcommand{\iadd}{\operatorname{iadd}\nolimits}
	\newcommand{\Gr}{\operatorname{Gr}\nolimits}
	\newcommand{\FGS}{\operatorname{FGS}\nolimits}
	\def \tbU{\tU_{\bullet}}
	\def \bU{\U_{\bullet}}
	\renewcommand{\mod}{\operatorname{mod^{\rm nil}}\nolimits}
	\newcommand{\proj}{\operatorname{proj}\nolimits}
	\newcommand{\inj}{\operatorname{inj.}\nolimits}
	\newcommand{\rad}{\operatorname{rad}\nolimits}
	\newcommand{\Span}{\operatorname{Span}\nolimits}
	\newcommand{\soc}{\operatorname{soc}\nolimits}
	\newcommand{\ind}{\operatorname{inj.dim}\nolimits}
	\newcommand{\Ginj}{\operatorname{Ginj}\nolimits}
	\newcommand{\res}{\operatorname{res}\nolimits}
	\newcommand{\np}{\operatorname{np}\nolimits}
	\newcommand{\Fac}{\operatorname{Fac}\nolimits}
	\newcommand{\Aut}{\operatorname{Aut}\nolimits}
	\newcommand{\DTr}{\operatorname{DTr}\nolimits}
	\newcommand{\TrD}{\operatorname{TrD}\nolimits}
	\newcommand{\Mod}{\operatorname{Mod}\nolimits}
	\newcommand{\R}{\operatorname{R}\nolimits}
	\newcommand{\End}{\operatorname{End}\nolimits}
	\newcommand{\lf}{\operatorname{l.f.}\nolimits}
	\newcommand{\Iso}{\operatorname{Iso}\nolimits}
	\newcommand{\aut}{\operatorname{Aut}\nolimits}
	\newcommand{\Ui}{{\mathbf U}^\imath}	\newcommand{\UU}{{\mathbf U}\otimes {\mathbf U}}
	\newcommand{\UUi}{(\UU)^\imath}
	\newcommand{\tUU}{{\tU}\otimes {\tU}}
	\newcommand{\tUUi}{(\tUU)^\imath}
	\newcommand{\tUi}{\widetilde{{\mathbf U}}^\imath}
	\newcommand{\sqq}{{\bf v}}
	\newcommand{\sqvs}{\sqrt{\vs}}
	\newcommand{\dbl}{\operatorname{dbl}\nolimits}
	\newcommand{\swa}{\operatorname{swap}\nolimits}
	\newcommand{\Gp}{\operatorname{Gp}\nolimits}
	\def \tTD{\widetilde{T}}
	\def \fX{\Upsilon}
	\newcommand{\U}{{\mathbf U}}
	\newcommand{\tU}{\widetilde{\mathbf U}}
	\newcommand{\fgm}{{\rm mod}^{{\rm fg}}}
	\newcommand{\fgmz}{\mathrm{mod}^{{\rm fg},\Z}}
	\newcommand{\fdmz}{\mathrm{mod}^{{\rm nil},\Z}}
	
	\newcommand{\ov}{\overline}
	\newcommand{\und}{\underline}
	\newcommand{\tk}{\widetilde{k}}
	\newcommand{\tK}{K}
	\newcommand{\tTT}{\operatorname{\widetilde{\mathbf{T}}}\nolimits}
	
	\def \hU{\widehat{\U}}
	\def \hUi{\widehat{\U}^\imath}
	\def \bfk{\mathbf{k}}
	\def \wI{\I_\circ}
	\def \cR{\mathcal{R}}
	\newcommand{\tH}{\operatorname{{\ch}_{\rm{tw}}}\nolimits}
	
	\newcommand{\utM}{\operatorname{\cm\ch}\nolimits}
	\newcommand{\tM}{\operatorname{\cs\cd\widetilde{\ch}}\nolimits}
	\newcommand{\rM}{\operatorname{\cm\ch_{\rm{red}}}\nolimits}
	\newcommand{\utMH}{\cs\cd\ch(\Lambda^\imath)}
	\newcommand{\tMH}{\cs\cd\widetilde{\ch}(\Lambda^\imath)}
	\newcommand{\tCMH}{{\cc\widetilde{\ch}(\K Q,\btau)}}
	
	\newcommand{\rMH}{\operatorname{\cs\cd\ch_{\rm{red}}(\Lambda^\imath)}\nolimits}
	\newcommand{\utMHg}{\operatorname{\ch(Q,\btau)}\nolimits}
	\newcommand{\tMHg}{\operatorname{\widetilde{\ch}(Q,\btau)}\nolimits}
	\newcommand{\tMHk}{{\widetilde{\ch}(\K Q,\btau)}}
	\newcommand{\rMHg}{\operatorname{\ch_{\rm{red}}(Q,\btau)}\nolimits}
	
	\newcommand{\rMHd}{\operatorname{\cm\ch_{\rm{red}}(\Lambda^\imath)_{\bvsd}}\nolimits}
	\newcommand{\tMHd}{\operatorname{\cs\cd\widetilde{\ch}(\Lambda^\imath)_{\bvsd}}\nolimits}
	
	\def \h{h} 
	\newcommand{\tMHl}{\cs\cd\widetilde{\ch}({\bs}_\ell\Lambda^\imath)}
	\newcommand{\rMHl}{\cm\ch_{\rm{red}}({\bs}_\ell\Lambda^\imath)_{\bvsd}}
	\newcommand{\tMHi}{\cs\cd\widetilde{\ch}({\bs}_i\Lambda^\imath)}
	\newcommand{\rMHi}{\cm\ch_{\rm{red}}({\bs}_i\Lambda^\imath)_{\bvsd}}
	\newcommand{\tMHgi}{\widetilde{\ch}({\bs}_i Q,\btau)}
	\def \bI{\I_\bullet}
	\newcommand{\utGpg}{\operatorname{\ch^{\rm Gp}(Q,\btau)}\nolimits}
	\newcommand{\tGpg}{\operatorname{\widetilde{\ch}^{\rm Gp}(Q,\btau)}\nolimits}
	\newcommand{\rGpg}{\operatorname{\ch_{red}^{\rm Gp}(Q,\btau)}\nolimits}

	\newcommand{\colim}{\operatorname{colim}\nolimits}
	\newcommand{\gldim}{\operatorname{gl.dim}\nolimits}
	\newcommand{\cone}{\operatorname{cone}\nolimits}
	\newcommand{\rep}{\operatorname{rep}\nolimits}
	\newcommand{\Ext}{\operatorname{Ext}\nolimits}
	\newcommand{\Tor}{\operatorname{Tor}\nolimits}
	\newcommand{\Hom}{\operatorname{Hom}\nolimits}
	\newcommand{\Top}{\operatorname{top}\nolimits}
	\newcommand{\Coker}{\operatorname{Coker}\nolimits}
	\newcommand{\thick}{\operatorname{thick}\nolimits}
	\newcommand{\rank}{\operatorname{rank}\nolimits}
	\newcommand{\Gproj}{\operatorname{Gproj}\nolimits}
	\newcommand{\Len}{\operatorname{Length}\nolimits}
	\newcommand{\RHom}{\operatorname{RHom}\nolimits}
	\renewcommand{\deg}{\operatorname{deg}\nolimits}
	\renewcommand{\Im}{\operatorname{Im}\nolimits}
	\newcommand{\Ker}{\operatorname{Ker}\nolimits}
	\newcommand{\Coh}{\operatorname{Coh}\nolimits}
	\newcommand{\Id}{\operatorname{Id}\nolimits}
	\newcommand{\Qcoh}{\operatorname{Qch}\nolimits}
	\newcommand{\CM}{\operatorname{CM}\nolimits}
	\newcommand{\sgn}{\operatorname{sgn}\nolimits}
	\newcommand{\Gdim}{\operatorname{G.dim}\nolimits}
	\newcommand{\fpr}{\operatorname{\mathcal{P}^{\leq1}}\nolimits}
	
	\newcommand{\For}{\operatorname{{\bf F}or}\nolimits}
	\newcommand{\coker}{\operatorname{Coker}\nolimits}
	\renewcommand{\dim}{\operatorname{dim}\nolimits}
	\newcommand{\rankv}{\operatorname{\underline{rank}}\nolimits}
	\newcommand{\dimv}{{\operatorname{\underline{dim}}\nolimits}}
	\newcommand{\diag}{{\operatorname{diag}\nolimits}}
	\newcommand{\qbinom}[2]{\begin{bmatrix} #1\\#2 \end{bmatrix} }
	
	\renewcommand{\Vec}{{\operatorname{Vec}\nolimits}}
	\newcommand{\pd}{\operatorname{proj.dim}\nolimits}
	\newcommand{\gr}{\operatorname{gr}\nolimits}
	\newcommand{\id}{\operatorname{Id}\nolimits}
	\newcommand{\Res}{\operatorname{Res}\nolimits}
	\def \tT{\widetilde{\mathscr{T} }}
	
	\def \bwi{w_{\bullet,i}}
	\def \tTL{\tT(\Lambda^\imath)}
	
	\newcommand{\mbf}{\mathbf}
	\newcommand{\mbb}{\mathbb}
	\newcommand{\mrm}{\mathrm}
	\newcommand{\cbinom}[2]{\left\{ \begin{matrix} #1\\#2 \end{matrix} \right\}}
	\newcommand{\dvev}[1]{{B_1|}_{\ev}^{{(#1)}}}
	\newcommand{\dv}[1]{{B_1|}_{\odd}^{{(#1)}}}
	\newcommand{\dvd}[1]{t_{\odd}^{{(#1)}}}
	\newcommand{\dvp}[1]{t_{\ev}^{{(#1)}}}
	\newcommand{\ev}{\bar{0}}
	\newcommand{\odd}{\bar{1}}
	\newcommand{\Iblack}{\I_{\bullet}}
	\newcommand{\wb}{w_\bullet}
	\newcommand{\Uidot}{\dot{\bold{U}}^{\imath}}
	
	\newcommand{\kk}{h}
	\newcommand{\la}{\lambda}
	\newcommand{\LR}[2]{\left\llbracket \begin{matrix} #1\\#2 \end{matrix} \right\rrbracket}
	\newcommand{\ff}{\mathbf{f}}
	\newcommand{\pdim}{\operatorname{proj.dim}\nolimits}
	\newcommand{\idim}{\operatorname{inj.dim}\nolimits}
	\newcommand{\Gd}{\operatorname{G.dim}\nolimits}
	\newcommand{\Ind}{\operatorname{Ind}\nolimits}
	\newcommand{\add}{\operatorname{add}\nolimits}
	\newcommand{\ad}{\operatorname{ad}\nolimits}
	\newcommand{\pr}{\operatorname{pr}\nolimits}
	\newcommand{\oR}{\operatorname{R}\nolimits}
	\newcommand{\oL}{\operatorname{L}\nolimits}
	\newcommand{\ext}{{ \mathfrak{Ext}}}
	\newcommand{\Perf}{{\mathfrak Perf}}
	\def\scrP{\mathscr{P}}
	\newcommand{\bk}{{\mathbb K}}
	\newcommand{\cc}{{\mathcal C}}
	\newcommand{\gc}{{\mathcal GC}}
	\newcommand{\dg}{{\rm dg}}
	\newcommand{\ce}{{\mathcal E}}
	\newcommand{\cs}{{\mathcal S}}
	\newcommand{\cl}{{\mathcal L}}
	\newcommand{\cf}{{\mathcal F}}
	\newcommand{\cx}{{\mathcal X}}
	\newcommand{\cy}{{\mathcal Y}}
	\newcommand{\ct}{{\mathcal T}}
	\newcommand{\cu}{{\mathcal U}}
	\newcommand{\cv}{{\mathcal V}}
	\newcommand{\cn}{{\mathcal N}}
	\newcommand{\mcr}{{\mathcal R}}
	\newcommand{\ch}{{\mathcal H}}
	\newcommand{\ca}{{\mathcal A}}
	\newcommand{\cb}{{\mathcal B}}
	\newcommand{\ci}{{\I}_{\btau}}
	\newcommand{\cj}{{\mathcal J}}
	\newcommand{\cm}{{\mathcal M}}
	\newcommand{\cp}{{\mathcal P}}
	\newcommand{\cg}{{\mathcal G}}
	\newcommand{\cw}{{\mathcal W}}
	\newcommand{\co}{{\mathcal O}}
	\newcommand{\cq}{{Q^{\rm dbl}}}
	\newcommand{\cd}{{\mathcal D}}
	\newcommand{\ck}{\widetilde{\mathcal K}}
	\newcommand{\calr}{{\mathcal R}}
	\newcommand{\iLa}{\Lambda^{\imath}}
	\newcommand{\La}{\Lambda}
	\newcommand{\ol}{\overline}
	\newcommand{\ul}{\underline}
	\newcommand{\st}{[1]}
	\newcommand{\ow}{\widetilde}
	\renewcommand{\P}{\mathbf{P}}
	\newcommand{\pic}{\operatorname{Pic}\nolimits}
	\newcommand{\Spec}{\operatorname{Spec}\nolimits}
	
	\newtheorem{theorem}{Theorem}[section]
	\newtheorem{acknowledgement}[theorem]{Acknowledgement}
	\newtheorem{conjecture}[theorem]{Conjecture}
	\newtheorem{corollary}[theorem]{Corollary}
	\newtheorem{definition}[theorem]{Definition}
	\newtheorem{example}[theorem]{Example}
	\newtheorem{lemma}[theorem]{Lemma}
	\newtheorem{notation}[theorem]{Notation}
	\newtheorem{problem}[theorem]{Problem}
	\newtheorem{proposition}[theorem]{Proposition}
	\newtheorem{summary}[theorem]{Summary}
	\numberwithin{equation}{section}

	\newtheorem{alphatheorem}{Theorem}
	\newtheorem{alphacorollary}[alphatheorem]{Corollary}
	\newtheorem{alphaproposition}[alphatheorem]{Proposition}
	\renewcommand*{\thealphatheorem}{\Alph{alphatheorem}}
	
	\theoremstyle{remark}
	\newtheorem{remark}[theorem]{Remark}

	\newcommand{\Pd}{\pi_*}
	\def \bvs{{\boldsymbol{\varsigma}}}
	\def \bvsd{{\boldsymbol{\varsigma}_{\diamond}}}
	\def \btau{{{\tau}}}
	
	\def \Br{\mathrm{Br}}
	\def \bp{{\mathbf p}}
	\def \bq{{\bm q}}
	\def \bv{{v}}
	\def \bs{{r}}
	
	\def \bfK{{\mathbf K}}
	
	\def \bw{w_\bullet}
	
	\newcommand{\tCMHg}{\cc\widetilde{\ch}(Q,\btau)}
	\newcommand{\bfv}{\mathbf{v}}
	\def \bA{{\mathbf A}}
	\def \ba{{\mathbf a}}
	\def \bi{{\mathbf i}}
	\def \bL{{\mathbf L}}
	\def \bF{{\mathbf F}}
	\def \bS{{\mathbf S}}
	\def \bC{{\mathbf C}}
	\def \bU{{\mathbf U}}
	\def \bc{{\mathbf c}}
	\def \fpi{\mathfrak{P}^\imath}
	\def \Ni{N^\imath}
	\def \fp{\mathfrak{P}}
	\def \fg{\mathfrak{g}}
	\def \fk{\fg^\theta}  
	
	\def\Inv{\mathrm{Inv}}
	\def \bbW{W^{\circ}}
	\def \bbw{{\boldsymbol{w}}}
	\def \BB{\mathbf{B}}
	\def \tB{\widetilde{\mathbf{B}}}
	\def \hB{\widehat{\mathbf{B}}}
	\def \reW{W^\circ} 
	\def \fn{\mathfrak{n}}
	\def \fh{\mathfrak{h}}
	\def \fu{\mathfrak{u}}
	\def \fv{\mathfrak{v}}
	\def \fa{\mathfrak{a}}
	\def \fq{\mathfrak{q}}
	\def \Z{{\Bbb Z}}
	\def \F{{\Bbb F}}
	\def \D{{\Bbb D}}
	\def \C{{\Bbb C}}
	\def \N{{\Bbb N}}
	\def \Q{{\Bbb Q}}
	\def \G{{\Bbb G}}
	\def \P{{\Bbb P}}
	\def \K{{\mathbb K}}
	\def \bK{{\Bbb K}}
	\def \J{{\Bbb J}}
	\def \E{{\Bbb E}}
	\def \A{{\Bbb A}}
	\def \L{{\Bbb L}}
	\def \I{{\Bbb I}}
	\def \BH{{\Bbb H}}
	\def \T{{\Bbb T}}
	\def \upLa {\raisebox{8pt}{\rotatebox{180}{$\Lambda$}}}
	\newcommand{\TT}
	{\operatorname{\mathbf{T}}\nolimits}
	\newcommand {\lu}[1]{\textcolor{red}{$\clubsuit$: #1}}
	
	\newcommand{\nc}{\newcommand}
	\newcommand{\browntext}[1]{\textcolor{brown}{#1}}
	\newcommand{\greentext}[1]{\textcolor{green}{#1}}
	\newcommand{\redtext}[1]{\textcolor{red}{#1}}
	\newcommand{\bluetext}[1]{\textcolor{blue}{#1}}
	\newcommand{\brown}[1]{\browntext{ #1}}
	\newcommand{\green}[1]{\greentext{ #1}}
	\newcommand{\red}[1]{\redtext{ #1}}
	\newcommand{\blue}[1]{\bluetext{ #1}}
	
	\newcommand{\wtodo}{\todo[inline,color=orange!20, caption={}]}
	\newcommand{\lutodo}{\todo[inline,color=green!20, caption={}]}

	\title[Dual canonical bases for $\mathrm i$Quantum groups]{$\mathrm i$Quantum groups and $\mathrm i$Hopf algebras II: dual canonical bases}
	
	\author[Jiayi Chen]{Jiayi Chen}
	\address{Department of Mathematics, Shantou University, Shantou 515063, P.R.China}
	\email{chenjiayi@stu.edu.cn}
	
	\author[Ming Lu]{Ming Lu}
	\address{Department of Mathematics, Sichuan University, Chengdu 610064, P.R.China}
	\email{luming@scu.edu.cn}

	\author[Xiaolong Pan]{Xiaolong Pan}
	\address{Department of Mathematics, Sichuan University, Chengdu 610064, P.R.China}
	\email{xiaolong\_pan@stu.scu.edu.cn}

	\author[Shiquan Ruan]{Shiquan Ruan}
	\address{ School of Mathematical Sciences,
		Xiamen University, Xiamen 361005, P.R.China}
	\email{sqruan@xmu.edu.cn}

	\author[Weiqiang Wang]{Weiqiang Wang}
	\address{Department of Mathematics, University of Virginia, Charlottesville, VA 22904, USA}
	\email{ww9c@virginia.edu}
	\subjclass[2020]{Primary 17B37, 20G42, 81R50}
	\keywords{iquantum groups, quantum groups, iHopf algebras, dual canonical bases, braid group symmetries}
	
	\begin{abstract}
		Building on the iHopf algebra realization of quasi-split universal iquantum groups developed in a prequel, we construct the dual canonical basis for a universal iquantum group of arbitrary finite type, which are further shown to be preserved by the ibraid group action; this recovers the results of Lu-Pan in ADE type obtained earlier in a geometric approach. Moreover, we identify the dual canonical basis for the Drinfeld double quantum group of arbitrary finite type, which is realized via iHopf algebra on the double Borel, with Berenstein-Greenstein's double canonical basis, settling several of their conjectures.
	\end{abstract}
	
	\maketitle
	\setcounter{tocdepth}{1}
	\tableofcontents

	\section{Introduction}
	
	
	Let $\U$ be the Drinfeld-Jimbo quantum group, $\tB$ be a Borel subalgebra of $\U$, and $\tU$ be the Drinfeld double on $\tB$. Qin \cite{Qin16} constructed a dual canonical basis on $\tU$ of ADE type with positive property, which contains the (rescaled) dual canonical bases of $\U^+$ and $\U^-$ of Lusztig \cite{Lus90a} over $\Z[v^{\frac12},v^{-\frac12}]$. Generalizing Hernandez-Leclerc \cite{HL15}, Qin's construction uses the quantum Grothendieck ring of perverse sheaves on Nakajima quiver varieties (cf. \cite{Na04, VV03}), which can also be viewed as a geometric counterpart for Bridgeland's Hall algebra construction of $\tU$ \cite{Br13}. Around the same time, Berenstein and Greenstein \cite{BG17a} constructed double canonical bases for $\tU$ of finite type algebraically; the relation between these two bases for $\tU$ remained unclear until recently. 
	
	Associated to any quasi-split Satake diagram $(\I,\tau)$ (with no $\tau$-fixed edge), a universal quasi-split iquantum group $\tUi$ is formulated and realized via an iHall algebra in \cite{LW22a, LW23}, which is a generalization of Bridgeland's Hall algebra. We note that $\tUi$ is a coideal subalgebra of $\tU$ and admits relative braid (=ibraid) group symmetries \cite{LW22b, WZ23, Z23}; also see \cite{KP11}. The Drinfeld double $\tU$ can be identified as a universal iquantum group associated to the diagonal Satake diagram $(\I \sqcup \I, \swa)$.   
	A dual canonical basis for $\tUi$ has been constructed in \cite{LW21b} for all quasi-split Satake diagrams of ADE type (other than type AIII$_{2r}$ which has a $\tau$-fixed edge) via quantum Grothendieck ring of perverse sheaves on Nakajima-Keller-Scherotzke quiver varieties. This generalizes \cite{Qin16}. 
	
	Two of the authors \cite{LP25} have recently made substantial progress further along the geometric directions. They gave a new construction of the dual canonical basis of $\tUi$ of ADE type (other than type AIII$_{2r}$) via rescaled iHall basis and Lusztig's lemma, connecting earlier constructions in \cite{LW22a} and \cite{LW21b}; moreover, they showed that the dual canonical basis is preserved by (rescaled) ibraid group symmetries $\tTT_i$ of $\tUi$. Specializing to the iquantum group of diagonal type, they show that the dual canonical basis on $\tU$ constructed by Qin \cite{Qin16} coincides with (two variants of) double canonical bases due to Berenstein and Greenstein. This allows them to settle several conjectures in \cite{BG17a} for ADE type; in particular, the dual canonical basis of $\tU$ is preserved by (rescaled) Lusztig's braid group symmetries $\widetilde{T}_i$. 
	
	The algebraic constructions in \cite{BG17a} are valid for $\tU$ of all finite type, and their construction of double canonical basis involves another ingenious yet quite complicated construction of Heisenberg doubles; again it contains (rescaled) dual canonical bases of $\U^+$ and $\U^-$ from \cite{Lus90a, Ka91}. Berenstein and Greenstein made several conjectures including that the double canonical basis in any finite type is preserved by the braid group action. A compatibility \cite[Theorem 1.2]{Lus96} between canonical bases in subalgebras of $\U^+$ under the braid group action was reformulated in \cite[Proposition 5.14]{BG17a} as that dual canonical bases of subalgebras of $\U^+$ are matched by the rescaled braid group action. 
	
	In a prequel \cite{CLPRW} we formulated a notion of iHopf algebras, a new associative algebra structure defined on Hopf algebras with Hopf pairings, and showed that the iHopf algebra $\tB^\imath_\tau$ on the Borel $\tB$ provides a realization of the $\imath$quantum groups $\tUi$; we shall identify $\tB^\imath_\tau\equiv \tUi$ hereafter. 
	
	The goal of this paper is to construct the dual canonical basis for $\tUi$ of arbitrary finite type in the framework of iHopf algebras, generalizing the main results of \cite{LP25}. Along the way, we develop direct connections between Lusztig's braid group action and ibraid group action. In particular, as the iHopf algebra defined on the double Borel $\tB\otimes \tB$ provides a realization of the Drinfeld double $\tU$, we construct the dual canonical basis of $\tU$ of arbitrary finite type, settling the main conjectures in \cite{BG17a}.
	\vspace{2mm}
	
	The presentations for $\tU$ and $\tUi$ used in this paper look a bit unusual, as they use dual Chevalley generators following \cite{BG17a, LP25, CLPRW}. We first strengthen the connection between braid and ibraid group symmetries initiated in \cite{CLPRW}. There is a natural embedding of Lusztig's algebra $\ff$ into $\tB$, $\iota:\ff \rightarrow \tB$, and also an embedding $\ff \rightarrow \tB^\imath_\tau\equiv \tUi$. Recall $\vartheta_i$ in \eqref{rescale:theta} is a rescaling of $\theta_i\in \ff$, $\tau_i$ in \eqref{taui} is a rank one analogue of the involution $\tau$, and $\ff_{i,\tau i}$ and ${}^\sigma\ff[i,\tau i]$ in \eqref{eq:firank1}--\eqref{eq:fitaui} are subalgebras of $\ff$. There exist (rescaled) braid group symmetries $\widetilde{T}_{i}$ and $\widetilde{T}_{w}$ in $\tU$ \cite{Lus90a, Lus90b, Lus93} as well as (rescaled) ibraid group symmetries $\tTT_i$ on $\tUi$ \cite{KP11, LW21a, WZ23, Z23}. 
	It was established in \cite[Theorem C]{CLPRW} that, for any $i,j\in \I$ such that $i\neq j,\tau j$,  
	\begin{align} \label{eq:TiPart1}
		\tTT_i(\vartheta_j)=\iota\big(\widetilde{T}_{r_i}(\vartheta_j)\big).
	\end{align}
	See \eqref{def:ri} for notation $r_i$.
	
	\begin{alphatheorem} [Theorems \ref{thm:Tifi} and \ref{thm:Tifif_i}]
		\label{thm:ibraid}
		\begin{enumerate}
			\item[(i)] 
			We have $\tTT_i^{-1}(x)=\widetilde{T}_{r_i}^{-1}(x)$, for any $x\in{^\sigma\ff}[i,\tau i]$.
			\item[(ii)] 
			For any $x\in{^\sigma\ff}[i,\tau i]$ and $u\in\ff_{i,\tau i}$, we have 
			\begin{align*}
				\widetilde{\TT}_i^{-1}(xu)=v^{-\frac{1}{2}(\tau\wt(u)+\wt(u),\wt(x))}\K_{\tau\tau_i\wt(u)}^{-1}\diamond \big(\tau\tau_i(u)\widetilde{T}_{r_i}^{-1}(x)\big).
			\end{align*}
			\item[(iii)] 
			For any $x\in\ff[i,\tau i]$ and $u\in\ff_{i,\tau i}$, we have 
			\begin{align*}
				\widetilde{\TT}_i(ux)=v^{-\frac{1}{2}(\tau\wt(u)+\wt(u),\wt(x))}\K_{\tau_i\wt(u)}^{-1}\diamond \big(\widetilde{T}_{r_i}(x)\tau\tau_i(u)\big).
			\end{align*}
		\end{enumerate}
	\end{alphatheorem}
	Theorem \ref{thm:ibraid} indicates that Lusztig's braid group action $\widetilde{T}_{r_i}$  on a subalgebra of $\ff$ is matched with the ibraid group action of $\tTT_i$ on part of $\tUi$, substantially improving \eqref{eq:TiPart1}; it is new even in the quantum group setting. This result is most naturally understood in the context of Hall and iHall algebras; compare \cite{LP25}. Recall \cite{LW22a,LW23} that iquiver algebras are defined by adding some arrows to the quivers and then used to realize iquantum groups. By restricting to the modules of quivers, the reflection functors of iquiver algebras \cite{LW21a,LW22b} coincide with the ones of quivers \cite{Rin96}. Therefore, both sides of the equalities in (i) above correspond to the same modules (though in different Hall algebras with different multiplications). For the statements (ii)--(iii), we recall the reflection functors reverse some arrows, and then change the structure of module categories. Take (ii) for example: $xu$ (not $ux$) corresponds naturally to a module $M$, and $\tau\tau_i(u)\widetilde{T}^{-1}_{r_i}(x)$ corresponds to the module $M$ acted by the reflection functor. 
	
	Theorem \ref{thm:ibraid} is valid in the Kac-Moody setting. 
	In the remainder of the Introduction, let us restrict ourselves to Drinfeld doubles and iquantum groups of arbitrary {\em finite} type. 
	
	The algebras $\tU$ (and resp. $\tB$, or $\tUi$) admits variants $\hU$ (and resp. $\hB$, or $\hUi$) where the generators $K_i, K_i'$ (and resp. $K_i$, or $\K_i$) of the Cartan subalgebras are not required to be invertible. When discussing about braid group symmetries, we need these Cartan generators to be invertible and so work with the tilde versions. The hat versions are natural from the viewpoints of Hall algebras and dual canonical bases. We shall identify the iHopf algebra $\hB^\imath_\tau$ on $\hB$ with $\hUi$: $\hB^\imath_\tau\equiv \hUi$. 
	
	For dual canonical basis, we use a version of bar involution (which is an anti-involution) on $\hUi$; cf. Lemma \ref{lem:bar}. Via the iHopf algebra construction, we import the dual canonical basis or a dual PBW basis of $\ff$ to $\hB^\imath_\tau\equiv \hUi$ via a linear embedding. We view such a basis of $\ff$ (after adjoining by Cartan) as a standard basis for $\hB^\imath_\tau$, and apply the bar involution to them. In this way, we are able to apply Lusztig's lemma to construct the dual canonical basis of $\hB^\imath_\tau$. Let $\K_\alpha$ denote an element in the Cartan subalgebra of $\hB^\imath_\tau$, and let $\mathbf{C}$ denote the dual canonical basis of $\ff$ with respect to the bilinear form $\varphi$ in \eqref{eq:hopf-pairing}. We refer to \eqref{eq:diamond} for the $\diamond$-action. 
	
	\begin{alphatheorem} [Theorem \ref{thm:dCB}, Proposition \ref{prop:dCB:PBW}]
		\label{thm:iDCB}
		There exists a unique bar-invariant element $C_{\alpha,b}\in \hB^\imath_\tau$ such that 
		$C_{\alpha,b}\in\K_\alpha\diamond \iota(b) +\sum_{(\alpha,b) \prec(\beta,b')}v^{-1}\Z[v^{-1}]\cdot \K_\beta\diamond \iota(b'),$  
		for $\alpha\in\N^\I$ and $b\in \mathbf{C}$.
		In addition, $C_{\alpha,b}=\K_\alpha\diamond C_{0,b}$. Then $\{C_{\alpha,b} \mid \alpha\in\N^\I,b\in \mathbf{C}\}$ forms the dual canonical basis for $\hB^\imath_\tau\equiv \hUi$.
	\end{alphatheorem}
	
	For ADE type (excluding AIII$_{2r}$), it follows from the new geometric construction in \cite{LP25} that the basis constructed in the theorem above matches with the dual canonical basis constructed \cite{LW21b, LP25}. Our algebraic approach is not sufficient to re-establish the positivity property of dual canonical bases for ADE type {\em loc. cit.} though. By construction here (and also in \cite{LP25}), there is an algorithm to compute the dual canonical basis for $\hUi$ in Theorem \ref{thm:iDCB}, which was missing in earlier works \cite{Qin16, LW21b}. 
	
	The iquantum groups $\Ui_\bvs$ with parameter $\bvs$ introduced earlier by G.~Letzter \cite{Let99} (see Kolb \cite{Ko14}) are recovered by central reductions from $\tUi$. By Proposition~ \ref{prop:dCB:parameter}, the dual canonical basis on $\tUi$ descends to a dual canonical basis on the iquantum group with a distinguished parameter ${\bvs_\diamond}$ defined in \eqref{eq:disting-para}. 
	
	Recall an anti-involution $\sigma^\imath$ on $\tUi$ from Lemma \ref{lem:involution-iQG}. The dual canonical basis of $\tUi$ constructed in Theorem \ref{thm:dCB} admits several symmetries.
	
	\begin{alphatheorem} [Theorem \ref{thm:dCB-braid}, Proposition \ref{prop:dualCB-anti-inv}, Corollary \ref{cor:dCB-tau}]
		\label{thm:DCB:ibraid}
		The dual canonical basis of $\tB_\tau^\imath\equiv \tUi$ is preserved under the ibraid group action. Moreover, it is also preserved by the anti-involution $\sigma^\imath$ and by the involution $\tau$. 
	\end{alphatheorem}
	As a consequence of Theorem \ref{thm:DCB:ibraid} (see Corollary \ref{cor:braid-simple}), we easily recover a difficult result (see \cite[Theorem 7.13]{WZ23}) in case of quasi-split iquantum groups:  $\tTT_{w}(B_i) = B_{wi}$ if $wi\in\I$, for $w \in W_\tau$ and $i \in \I$. (This includes the well-known quantum group counterpart \cite{Lus93, Jan96}.)  
	
	Theorem \ref{thm:ibraid} plays an essential role in the proof of Theorem \ref{thm:DCB:ibraid}. In order to prove Theorem~ \ref{thm:DCB:ibraid}, we first use a dual PBW basis of $\ff$ to construct the dual canonical basis of $\widehat{\mathbf{B}}^\imath_\tau \equiv\tUi$; see Proposition \ref{prop:dCB:PBW}. Then we apply Theorem \ref{thm:ibraid} to prove that the braid group action $\tTT_i$ sends a dual PBW basis to another dual PBW basis. Theorem \ref{thm:DCB:ibraid} follows then from Proposition \ref{prop:dCB:PBW} and the uniqueness of the dual canonical basis.  
	
	We now specialize the above results to the Drinfeld double quantum group $\tU$ of finite type, viewed as an iquantum group of diagonal type. By the iHopf construction, the dual canonical basis for $\tU$ has a tensor product of dual canonical basis elements of $\ff$ (adjoint with a Cartan algebra factor) as a leading term. This characterization allows us to bridge and compare with the constructions in \cite{BG17a}. Recall the anti-involution $\sigma$ and Chevalley involution $\omega$ on $\tU$ from Lemma \ref{lem:anti-involut-QG}.
	
	\begin{alphatheorem} [Theorem \ref{thm:double=dCB}, Corollaries \ref{cor:dCB:braid}--\ref{cor:dCB:sigma}]
		\label{thm:BG=dCB}
		The dual canonical basis on $\tU$ coincides with the double canonical basis on $\tU$ \`a la Berenstein-Greenstein. Moreover, this basis is preserved by the braid group action, by the Chevalley involution $\omega$, and by the anti-involution $\sigma$.
	\end{alphatheorem}
	
	As explained in Section \ref{sec:doubleCB}, two variants of double canonical bases for $\tU$ were constructed in \cite{BG17a} via two different processes through Heisenberg doubles, and they were conjectured {\em loc. cit.} to coincide. Theorem \ref{thm:double=dCB} shows that this is indeed the case. Our approach bypasses Heisenberg doubles completely. 
	
	The dual canonical bases on $\hU$ and $\hUi$ seem to be more aligned with monoidal categorification or connections to cluster algebras, and it will be interesting to formulate such connections precisely. An intriguing question remains whether there is any direct connection between the canonical bases on modified quantum group $\dot\U$ \cite{Lus93} (and on modified iquantum group $\Uidot$ \cite{BW18b, BW21}) and the dual canonical bases on $\hU$ (and on $\hUi$) constructed here. 
	
	It will also be very interesting but highly nontrivial to generalize our work to iquantum groups beyond quasi-split types; see \cite{BW18b, BK19, BW21, WZ23} for some constructions in such generalities.

	\vspace{2mm}
	
	The paper is organized as follows. 
	In Section \ref{sec:QG and iQG}, we review quantum groups and iquantum groups, including (relative) braid group actions. We also review the iHopf algebra realization of the iquantum group $\hUi$ and several properties arising this way. 
	
	In Section \ref{sec:braid}, via the identification of $\ff$ as subspaces in both $\hB$ and $\hB^\imath_\tau\equiv \hUi$, we establish Theorem \ref{thm:ibraid} relating Lusztig's braid group action to ibraid group action. 
	This result is used in Section \ref{sec:DCBiQG} to establish Theorem \ref{thm:DCB:ibraid}. The dual canonical basis of $\hUi$ is also established in Section \ref{sec:DCBiQG}.
	
	Finally, in Section \ref{sec:doubleCB} we specialize our results on the dual canonical basis to $\tU$, and show they coincide with the double canonical basis \`a la Berenstein-Greenstein. 
	In Appendix \ref{sec:rank1}, recursive formulas for dual canonical basis elements in quasi-split rank one are obtained.

	\vspace{2mm}
	\noindent{\bf Acknowledgments} 
	ML is partially supported by the National Natural Science Foundation of China (No. 12171333). SR is partially supported by
	Fundamental Research Funds for Central Universities of China (No. 20720250059), Fujian Provincial Natural Science Foundation of China (No. 2024J010006) and
	the National Natural Science Foundation of China (Nos. 12271448 and 12471035). WW is partially supported by the NSF grant DMS-2401351, and he thanks National University of Singapore (Department of Mathematics and IMS) for providing an excellent research environment and support during his visit.

	\section{Quantum groups, iquantum groups and iHopf algebras} 
	\label{sec:QG and iQG}
	
	In this preliminary section, we recall quantum groups and iquantum groups in terms of (less standard) dual generators. We also review the realization of iquantum groups via iHopf algebras given in \cite{CLPRW}.
	
	\subsection{Quantum groups}
	\label{subsec:QG}
	
	Let $\I=\{1,\dots,n\}$. 
	Let $C=(c_{ij})_{i,j \in \I}$ be the symmetrizable generalized Cartan matrix (GCM) of a Kac-Moody Lie algebra $\fg$. Let $D=\diag(d_i\mid i\in \I)$ with $d_i\in\Z_{>0}$ be the symmetrizer of $C$, i.e., $DC$ is symmetric. 
	Let $\{\alpha_i\mid i\in\I\}$ be a set of simple roots of $\fg$, and denote the root lattice by $\Z^{\I}:=\Z\alpha_1\oplus\cdots\oplus\Z\alpha_n$. We define a symmetric bilinear form on $\Z^\I$ by setting
	\begin{align} \label{BilForm}
		(\alpha_i,\alpha_j)=d_ic_{ij},\quad \forall i,j\in\I.
	\end{align}
	The simple reflection $s_i:\Z^{\I}\rightarrow\Z^{\I}$ is defined to be $s_i(\alpha_j)=\alpha_j-c_{ij}\alpha_i$, for $i,j\in \I$.
	Denote the Weyl group by $W =\langle s_i\mid i\in \I\rangle$.
	
	Let $v$ be an indeterminate. Let 
	$$v_i=v^{d_i},\qquad \forall i\in\I.$$
	For $A,B$ in a $\Q(v^{\frac12})$-algebra, we write $[A, B]=AB-BA$, and $[A,B]_q=AB-qBA$ for any $q\in\Q(v^{\frac12})^\times$. Denote, for $r\in\N,m \in \Z$,
	\[
	[r]_{v_i}=\frac{v_i^r-v_i^{-r}}{v_i-v_i^{-1}},
	\quad
	[r]_{v_i}^!=\prod_{i=1}^r [i]_{v_i}, \quad \qbinom{m}{r}_{v_i} =\frac{[m]_{v_i}[m-1]_{v_i}\ldots [m-r+1]_{v_i}}{[r]_{v_i}^!}.
	\]
	Following \cite{Dr87, BG17a}, the (Drinfeld double) quantum group $\hU := \hU_v(\fg)$ is defined to be the $\Q(v^{\frac12})$-algebra generated by $E_i,F_i, \tK_i,\tK_i'$, $i\in \I$, subject to the following relations:  for $i, j \in \I$,
	\begin{align}
		[E_i,F_j]= \delta_{ij}(v_i^{-1}-v_i) (\tK_i-\tK_i'),  &\qquad [\tK_i,\tK_j]=[\tK_i,\tK_j']  =[\tK_i',\tK_j']=0,
		\label{eq:KK}
		\\
		\tK_i E_j=v_i^{c_{ij}} E_j \tK_i, & \qquad \tK_i F_j=v_i^{-c_{ij}} F_j \tK_i,
		\label{eq:EK}
		\\
		\tK_i' E_j=v_i^{-c_{ij}} E_j \tK_i', & \qquad \tK_i' F_j=v_i^{c_{ij}} F_j \tK_i',
		\label{eq:K2}
	\end{align}
	and for $i\neq j \in \I$,
	\begin{align}
		& \sum_{r=0}^{1-c_{ij}} (-1)^r \left[ \begin{array}{c} 1-c_{ij} \\r \end{array} \right]_{v_i}  E_i^r E_j  E_i^{1-c_{ij}-r}=0,
		\label{eq:serre1} 
		\\ 
		&\sum_{r=0}^{1-c_{ij}} (-1)^r \left[ \begin{array}{c} 1-c_{ij} \\r \end{array} \right]_{v_i}  F_i^r F_j  F_i^{1-c_{ij}-r}=0.
		\label{eq:serre2}
	\end{align}
	
	We define $\tU=\tU_v(\fg)$ as the $\Q(v^{\frac12})$-algebra with   generators and relations of $\widehat{\U}$ above, but in addition requiring $\tK_i,\tK_i'$ ($i\in\I$) to be invertible. Then $\tU$ and $\hU$ are $\Z^\I$-graded algebras by setting 
	$$\deg E_i=\alpha_i, \qquad\deg F_i=-\alpha_i,\qquad \deg K_i=0=\deg K_i'.$$ 
	Let $\tU_\mu$ be the homogeneous subspace of degree $\mu$.  Then $\tU=\oplus_{\mu\in\Z^\I} \tU_\mu$ and $\hU=\oplus_{\mu\in\Z^\I} \hU_\mu$. 
	
	The Drinfeld-Jimbo quantum group $\bU$ is defined to the $\Q(v^{\frac12})$-algebra generated by $E_i,F_i, K_i, K_i^{-1}$, $i\in \I$, subject to the relations modified from \eqref{eq:KK}--\eqref{eq:serre2} with $\tK_i'$ replaced by $K_i^{-1}$. We can also view $\bU$ as the quotient algebra of $\hU$ (or $\tU$) modulo the ideal generated by $K_iK_i'-1$ ($i\in\I$); see \cite{Dr87}.

	By a slight abuse of notation, let $\U^+$ be the subalgebra of $\hU$ (and also $\tU$, $\U$) generated by $E_i$ $(i\in \I)$, and let $\U^-$ be the subalgebra generated by $F_i$ $(i\in\I)$. Let $\hU^0$ and $\tU^0$ be the subalgebras of $\hU$ and $\tU$ generated by $\tK_i, \tK_i'$ $(i\in \I)$, and $\U^0$ be the subalgebra of $\bU$ generated by $\tK_i^{\pm 1}$ $(i\in \I)$. Then the algebras $\hU$, $\widetilde{\bU}$ and $\bU$ have triangular decompositions:
	\begin{align*}
		\hU=\U^+\otimes\hU^0\otimes \U^-,\qquad 
		\widetilde{\bU} =\U^+\otimes \widetilde{\bU}^0\otimes\U^-,
		\qquad
		\bU &=\bU^+\otimes \bU^0\otimes\bU^-.
	\end{align*}
	For any $\mu=\sum_{i\in\I}m_i\alpha_i\in\Z^\I$, we denote $K_\mu=\prod_{i\in\I} K_i^{m_i}$, $K_\mu'=\prod_{i\in\I} (K_i')^{m_i}$.

	The algebras $\hU$ (and $\tU$, $\U$) are Hopf algebras, with the coproduct $\Delta$ and the counit $\varepsilon$ defined by
	\begin{align}
		\begin{split}
			\Delta(E_i)=E_i\otimes 1+K_i\otimes E_i,\quad &\Delta(F_i)=1\otimes F_i+F_i\otimes K_i',
			\\
			\Delta(K_i)=K_i\otimes K_i,\quad &\Delta(K_i')=K_i'\otimes K_i';
			\\
			\varepsilon(E_i)=0=\varepsilon(F_i),\qquad &\varepsilon(K_i)=1=\varepsilon(K_i');
		\end{split}
	\end{align}
	
	The following two lemmas are either standard or easy to verify.
	\begin{lemma} 
		\label{lem:anti-involut-QG}
		\quad
		\begin{enumerate}
			\item There exists an anti-involution (called the bar-involution) $u\mapsto \ov{u}$ on $\hU$ (and also $\tU$, $\U$) given by $\ov{v^{1/2}}=v^{-1/2}$, $\ov{E_i}=E_i$, $\ov{F_i}=F_i$, and $\ov{K_i}=K_i$, $\ov{K_i'}=K_i'$, for $i\in\I$.
			\item 
			There exists an anti-involution $\sigma$ on $\hU$ (also $\tU$, $\U$) given by $\sigma(E_i)=E_i$, $\sigma(F_i)=F_i$, and $\sigma(K_i)=K_i'$, for $i\in\I$.
			\item 
			There exists a Chevalley involution $\omega$ on $\hU$ (also $\tU$, $\U$) given by $\omega(E_i)=F_i$, $\omega(F_i)=E_i$, and $\omega(K_i)=K_i'$, for $i\in\I$.
		\end{enumerate}
	\end{lemma}
	
	\begin{lemma} \label{lem:twisting Psi}
		Let $\F$ be the algebraic closure of $\Q(v^{\frac{1}{2}})$ and $\F^\times=\F\setminus\{0\}$. For scalars $\ba=(a_i)_{i\in\I}\in (\F^{\times})^\I$, we have an automorphism $\widetilde{\Psi}_{\ba}$ on the $\F$-algebra $\tU$ such that
		\[
		\widetilde{\Psi}_{\ba}:K_i\mapsto a_i^{\frac{1}{2}}K_i,\quad K'_i\mapsto a_i^{\frac{1}{2}}K'_i,\quad E_i\mapsto a_i^{\frac{1}{2}}E_i,\quad F_i\mapsto F_i.
		\]
	\end{lemma}
	
	Let $\Br(W)$ be the braid group associated to the Weyl group $W$, generated by simple reflections $t_i$ ($i\in\I$). Lusztig introduced 4 variants of braid group symmetries on the quantum group $\U$ \cite{Lus90b} \cite[\S37.1.3]{Lus93}. These braid group symmetries can be lifted to the Drinfeld double $\tU$; see, e.g., \cite[Propositions 6.20–6.21]{LW22b}, which are denoted by $\widetilde{T}_{i,e}',\widetilde{T}_{i,-e}''$, $e=\pm1$; also cf. \cite[\S5]{BG17a}. In fact, $\widetilde{T}_{i,e}'=\sigma\circ\widetilde{T}_{i,-e}''\circ\sigma$, which is the inverse of $\widetilde{T}_{i,-e}''$.
	
	\begin{proposition}
		\label{prop:BG2U}
		For $i\in \I$, the automorphisms $\widetilde{T}_{i,e}'$ on $\tU$ satisfy that
		\begin{align*}
			&\widetilde{T}_{i,e}'(K_\mu)= K_{s_i(\mu)},
			\qquad \widetilde{T}_{i,e}'(K'_\mu)= K'_{s_i(\mu)},\;\;\forall \mu\in \Z^\I,\\
			&\widetilde{T}_{i,1}'(E_i)=v_i (K_i')^{-1}F_i,\qquad \widetilde{T}_{i,1}'(F_i)=v_i^{-1} E_iK_i^{-1},\\
			&\widetilde{T}_{i,-1}'(E_i)=v_i^{-1} K_i^{-1}F_i,\qquad \widetilde{T}_{i,-1}'(F_i)=v_i E_i(K_i')^{-1},\\
			&\widetilde{T}_{i,e}'(E_j)
			=\sum_{r+s=-c_{ij}} (-1)^{r}v_i^{-e(r+\frac{1}{2}c_{ij})}(v_i-v_i^{-1})^{c_{ij}} E_i^{(s)}  E_j E_i^{(r)} \quad\forall j\neq i,
			\\
			&\widetilde{T}_{i,e}'(F_j)
			= \sum_{r+s=-c_{ij}} (-1)^{r}v_i^{-e(r+\frac{1}{2}c_{ij})}(v_i-v_i^{-1})^{c_{ij}} F_i^{(s)}  F_j F_i^{(r)} \quad\forall j\neq i.
		\end{align*}
	\end{proposition}

	\begin{lemma}\label{lem:QGbraid-bar}
		The braid group actions $\widetilde{T}_{i,e}'$ and $\widetilde{T}_{i,-e}''$ commute with the bar-involution, i.e., $\ov{\widetilde{T}_{i,e}'(u)}=\widetilde{T}_{i,e}'(\ov{u})$ and $\ov{\widetilde{T}_{i,-e}''(u)}=\widetilde{T}_{i,-e}''(\ov{u})$ for any $u\in\tU$.
	\end{lemma}

	We shall often use the shorthand notation 
	\begin{align}
		\label{braid-shorthand}
		\widetilde{T}_i:=\widetilde{T}_{i,1}',\qquad \widetilde{T}_i^{-1}:=\widetilde{T}_{i,-1}''.
	\end{align}
	The $\widetilde{T}_i$'s satisfy the braid group relations and so 
	$\widetilde{T}_w 
	:= \widetilde{T}_{i_1}\cdots
	\widetilde{T}_{i_r} \in \Aut(\tU)$ is well defined,
	where $w = s_{i_1}\cdots s_{i_r}$ is any reduced expression of $w \in W$.

	\subsection{iQuantum groups}
	
	For a Cartan matrix $C=(c_{ij})_{i,j\in \I}$, let $\text{Inv}(C)$ be the group of permutations $\btau$ of the set $\I$ such that $c_{ij}=c_{\btau i,\btau j}$, for all $i,j$, and $\btau^2=\Id$. Then  $\btau \in \text{Inv}(C)$ can be viewed as an involution (which is allowed to be the identity) of the corresponding Dynkin diagram (which is identified with $\I$ by abuse of notation). We shall refer to the pair $(\I,\btau)$ as a (quasi-split) Satake diagram. 
	
	We denote by $\bs_{i}$ the following element of order 2 in the Weyl group $W$, i.e.,
	\begin{align}
		\label{def:ri}
		\bs_i= \left\{
		\begin{array}{ll}
			s_{i}, & \text{ if } c_{i,\tau i}=2 \, (i.e.,\btau i=i),
			\\
			s_is_{\btau i}, & \text{ if } c_{i,\tau i}=0,
			\\
			s_is_{\btau i}s_i, & \text{ if } c_{i,\tau i}=-1.
		\end{array}
		\right.
	\end{align}
	It is well known that the {\rm restricted Weyl group} associated to $(\I,\tau)$ can be identified with the following subgroup $W_\btau$ of $W$:
	\begin{align}
		\label{eq:Wtau}
		W_{\btau} =\{w\in W\mid \btau w =w \btau\},
	\end{align}
	where $\btau$ is regarded as an automorphism of the root lattice $\Z^\I$. Moreover, the restricted Weyl group $W_{\btau}$ can be identified with a Weyl group with $\bs_i$ ($i\in \I_\btau$) as its simple reflections. 

	Associated with the Satake diagram $(\I,\tau)$, following \cite{LW22a} we define the universal iquantum groups ${\hU}^\imath$ (resp. $\tUi$) to be the $\Q(v^{\frac12})$-subalgebra of $\hU$ (resp. $\tU$) generated by
	\begin{equation}
		\label{eq:Bi}
		B_i= F_i +  E_{\btau i} \tK_i',
		\qquad \tk_i = \tK_i \tK_{\btau i}', 
		\quad \forall i \in \I, 
	\end{equation}
	(with $\tk_i$ invertible in $\tUi$). Let $\hU^{\imath 0}$ be the $\Q(v^{\frac12})$-subalgebra of $\hUi$ generated by $\tk_i$, for $i\in \I$. Similarly, let $\tU^{\imath 0}$ be the $\Q(v^{\frac12})$-subalgebra of $\tUi$ generated by $\tk_i^{\pm1}$, for $i\in \I$. 
	The algebra $\widetilde{\bU}^\imath$ (resp. $\hUi$) is a right coideal subalgebra of $\widetilde{\bU}$ (resp. $\hU$); the pairs $(\widetilde{\bU}, \widetilde{\bU}^\imath)$ and $(\hU,\hUi)$ are called quantum symmetric pairs, and $\hUi$ and $\tUi$ are called the universal {\em (quasi-split) iquantum groups}; they are {\em split} if $\btau =\Id$.
	
	For $i\in\I$, for any $\alpha=\sum_{i\in\I} a_i\alpha_i\in\Z^\I$, we set
	\begin{align}
		\label{eq:bbKi}
		\K_i=v^{\frac{1}{2}c_{i,\tau i}}\tk_i,\qquad \K_\alpha:=\prod_{i\in\I}\K_i^{a_i}.
	\end{align}

	Let $\bvs=(\vs_i)\in(\Q(v^{\frac12})^\times)^\I$ be such that $\vs_i=\vs_{\tau i}$ for each $i\in\I$ which satisfies $c_{i,\tau i}=0$. The iquantum groups \`a la Letzter-Kolb \cite{Let99,Ko14} $\Ui=\Ui_\bvs$ is the $\Q(v^{\frac12})$-subalgebra of $\U$ generated by
	$$B_i=F_i+\vs_i E_{\tau i}K_i^{-1},\quad k_i=K_iK_{\tau i}^{-1},\quad \forall i\in\I.$$ 
	By \cite[Proposition 6.2]{LW22a}, the $\Q(v^{\frac12})$-algebra $\Ui$ is isomorphic to the quotient of $\tUi$ by the ideal generated by $\tk_i-\vs_i$ (for $i=\tau i$) and $\tk_i\tk_{\tau i}-\vs_i\vs_{\tau i}$ (for $i\neq \tau i$).

	\begin{lemma}[see e.g. {\cite[Lemma 3.9]{CLPRW}}] 
		\label{lem:involution-iQG}
		\quad
		\begin{enumerate}
			\item 
			There exists an anti-involution $\sigma^\imath$ on $\hUi$ (and also $\tUi$) given by $\sigma^\imath(B_i)=B_i$, $\sigma^\imath(\tk_i)=\tk_{\tau i}$, for $i\in\I$.
			\item There exists an anti-involution (called bar-involution) $:u\mapsto \ov{u}$ on $\hUi$ (and also $\tUi$) given by $\ov{v^{1/2}}=v^{-1/2}$, $\ov{B_i}=B_i$, and $\ov{\K_i}=\K_i$,  for $i\in\I$. In particular, $\ov{\tk_i}=v_i^{c_{i,\tau i}}\tk_{i}$. 
			\item There exists an involution $\psi^\imath$ of $\tUi$ such that $\psi^\imath(v^{1/2})=v^{-1/2}$, $\psi^\imath(B_i)=B_i$, $\psi^\imath(\tk_i)=v_i^{c_{i,\tau i}}\tk_{\tau i}$, for $i\in\I$.
		\end{enumerate}
	\end{lemma}
	

	\begin{example}
		\label{ex:QGvsiQG}
		{\rm (Quantum groups as iquantum groups of diagonal type)} 
		Consider the $\Q(v^{\frac12})$-subalgebra $\tUUi$ of $\tUU$
		generated by
		\[
		\ck_i:=\tK_{i} \tK_{i^{\diamond}}', \quad
		\ck_i':=\tK_{i^{\diamond}} \tK_{i}',  \quad
		\cb_{i}:= F_{i}+ E_{i^{\diamond}} \tK_{i}', \quad
		\cb_{i^{\diamond}}:=F_{i^{\diamond}}+ E_{i} \tK_{i^{\diamond}}',
		\qquad \forall i\in \I.
		\]
		Here we drop the tensor product notation and use instead $i^\diamond$ to index the generators of the second copy of $\tU$ in $\tUU$. There exists a $\Q(v^{\frac12})$-algebra isomorphism $\widetilde{\phi}: \tU \rightarrow \tUUi$ such that
		\[
		\widetilde{\phi}(E_i)= \cb_{i},\quad \widetilde{\phi}(F_i)= \cb_{i^{\diamond}}, \quad \widetilde{\phi}(\tK_i)= \ck_i', \quad \widetilde{\phi}(\tK_i')= \ck_i, \qquad \forall  i\in \I.
		\]
		In this case, the Satake diagram is $(\I\sqcup\I^\diamond,\swa)$, where $\I^\diamond$ is a copy of $\I$ of $\tU$.  
	\end{example}

	\subsection{Relative braid group symmetries}
	
	Choose one representative for each $\btau$-orbit on $\I$, and let
	\begin{align}\label{eq:ci}
		\ci = \{ \text{the chosen representatives of $\btau$-orbits in $\I$} \}.
	\end{align} 
	The braid group associated to the relative Weyl group $W_\tau$ is denoted
	\begin{equation}
		\label{eq:braidCox}
		\brW =\langle \br_i \mid i\in \I_\btau \rangle
	\end{equation}
	where $\br_i$ satisfy the same braid relations as for $\bs_i$ in $W_{\tau }$. The relative braid (or ibraid) group symmetries $\tTT_{i,e}'$ and $\tTT''_{i,e}$ ($i\in\I$, $e\in\{+1,-1\}$) on $\tUi$ are established in \cite{LW22a, WZ23} (and \cite{Z23}); see \cite{KP11} for earlier conjectures on iquantum groups with specific parameters. In this paper, we shall also use the bar-equivariant versions of these ibraid group symmetries of $\tUi$. 
	
	Quasi K-matrix appeared earlier in different formulations; see \cite{BW18a, BK19,AV22}. We shall need the following.
	
	\begin{proposition}
		[\text{\cite[Theorem~3.6]{WZ23}}]
		\label{prop:Kmatrix}
		There exists a unique element $\widetilde{\Upsilon}=\sum_{\mu\in\N^\I}\widetilde{\Upsilon}^\mu$ (called quasi K-matrix) such that $\widetilde{\Upsilon}^0=1$, $\widetilde{\Upsilon}^\mu\in\U^+_{\mu}$ and the following identities hold:
		\begin{align}
			B_i\widetilde{\Upsilon}&=\widetilde{\Upsilon}B_i^\sigma,\qquad (i\in\I),
			\\
			x\widetilde{\Upsilon}&=\widetilde{\Upsilon}x,\qquad (x\in\tU^{\imath0}),
		\end{align}
		where $B_i^\sigma:=\sigma(B_i)=F_i+K_iE_{\tau i}$.  Moreover, $\widetilde{\Upsilon}^\mu=0$ unless $\tau(\mu)=\mu$.
	\end{proposition}
	
	Denote by $\tU_{i,\tau i}$ the quantum group associated to $\I_i=\{i,\tau i\}$. Let $\widetilde{\Upsilon}_i$ be the rank one quasi K-matrix associated to $\I_i=\{i,\tau i\}$, i.e., $\widetilde{\Upsilon}_i=\sum_{\mu\in\N^\I}\widetilde{\Upsilon}_{i}^{\mu}$ with $\widetilde{\Upsilon}_{i}^{\mu}\in\tU_{i,\tau i}^+$, and $\widetilde{\Upsilon}_{i}^{0}=1$. Define a distinguished parameter $\bm{\varsigma}_\diamond=(\varsigma_{i,\diamond})_{i\in\I}$ by
	\begin{align}
		\label{eq:disting-para}
		\varsigma_{i,\diamond}=v^{-\frac{1}{2}(\alpha_i,\alpha_{\tau i})}.
	\end{align}
	Recall the automorphism $\widetilde{\Psi}_{\varsigma_{\diamond}}$ of $\tU$ from Lemma~\ref{lem:twisting Psi}. We set
	\begin{equation}\label{eq: T_i twisted def}
		\widetilde{\mathscr{T}}_i:=\widetilde{\Psi}_{\bm{\varsigma}_\diamond}^{-1}\circ\widetilde{T}_i\circ \widetilde{\Psi}_{\bm{\varsigma}_\diamond},
		\qquad \widetilde{\mathscr{T}}_i^{-1}:=\widetilde{\Psi}_{\bm{\varsigma}_\diamond}^{-1}\circ\widetilde{T}_i^{-1}\circ \widetilde{\Psi}_{\bm{\varsigma}_\diamond}.
	\end{equation}
	Clearly $\widetilde{\mathscr{T}}_i$ and $\widetilde{\mathscr{T}}_i^{-1}$, for $i\in\I$, are automorphisms of $\tU$ and satisfy the braid group relations. Hence, we can define
	$\widetilde{\mathscr{T}}_w:=\widetilde{\mathscr{T}}_{i_1}\cdots \widetilde{\mathscr{T}}_{i_r},$
	where $w=s_{i_1}\cdots s_{i_r}$ is any reduced expression. 
	
	\begin{theorem}[{cf. \cite{WZ23,Z23}}]
		\label{thm: relative T_i conjugate}
		For $i\in\I$, there are mutually inverse automorphisms $\widetilde{\TT}_i$ and $\widetilde{\TT}_i^{-1}$ on $\tUi$ such that
		\begin{align}
			\label{eq:relbraid1}
			\widetilde{\TT}_i^{-1}(x)\widetilde{\Upsilon}_i &=\widetilde{\Upsilon}_i\widetilde{\mathscr{T}}_{\br_i}^{-1}(x),
			\\
			\label{eq:relbraid2}
			\widetilde{\TT}_i(x)\widetilde{\mathscr{T}}_i^{-1}(\widetilde{\Upsilon}_i)^{-1} &=\widetilde{\mathscr{T}}_i^{-1}(\widetilde{\Upsilon}_i)^{-1}\widetilde{\mathscr{T}}_{\br_i}(x).
		\end{align}
		Moreover, we have $\widetilde{\TT}_i^{-1}=\sigma^\imath\circ \widetilde{\TT}_i\circ \sigma^\imath$, and there exists a group homomorphism $\Br(W_\tau)\rightarrow \Aut(\tUi)$, $\br_i\mapsto \tTT_i$ for $i\in\I$.
	\end{theorem}

	Denote by $\tau_i$ the diagram involution of $\I_i:=\{i,\tau i\}$ defined by 
	\begin{align} \label{taui}
		r_i(\alpha_i)=-\alpha_{\tau_i(i)}, 
		\qquad
		r_i(\alpha_{\tau i})=-\alpha_{\tau_i(\tau i)}.
	\end{align}
	
	\begin{proposition} [{\cite[Proposition 4.11, Theorem 4.14]{WZ23}}]
		\label{prop:TiBi}
		For $i,j\in\I$, we have $\widetilde{\TT}_i(\K_j)=\K_{r_i(\alpha_j)}$ and
		\[\widetilde{\TT}_i(B_i)=v^{\frac{1}{2}(\alpha_i-\alpha_{\tau i},\alpha_i)}\K_{\tau_i(i)}^{-1}B_{\tau_i(\tau i)},\quad \widetilde{\TT}_i(B_{\tau i})=v^{\frac{1}{2}(\alpha_i-\alpha_{\tau i},\alpha_i)}\K_{\tau_i(\tau i)}^{-1}B_{\tau_i(i)}.
		\]
	\end{proposition}

	\begin{lemma} [cf. \cite{CLPRW}]
		\label{lem:bar-invar-braid-Ui}
		The braid group actions $\tTT_{i}$ commute with the bar-involution, i.e., $\ov{\tTT_i(u)}=\tTT_i(\ov{u})$ for any $u\in\tUi$.
	\end{lemma}
	
	Corresponding to Lusztig's braid group symmetries $\widetilde{T}_{i,e}'$, $\widetilde{T}_{i,e}''$ on $\tU$, as in \cite{LW22b,WZ23,Z23}, we define
	\begin{align}
		\tTT_{i,1}'&=\tTT_{i},\quad \tTT_{i,-1}''=\tTT_{i}^{-1},
		\\
		\tTT_{i,-1}'&=\psi^\imath\circ \tTT_i\circ \psi^\imath,
		\quad \tTT_{i,1}''=\psi^\imath\circ \tTT_i^{-1}\circ \psi^\imath.
	\end{align}
	Moreover, we have 
	$$\tTT_{i,e}'=\sigma^\imath\circ \TT_{i,-e}''\circ \sigma^\imath,\quad e\in\{+1,-1\}, i\in\I.$$
	Then all the braid group actions $\tTT_{i,e}',\tTT_{i,e}''$ commute with the bar-involution since $\sigma^\imath,\psi^\imath$ commute with the bar-involution.

	\subsection{iHopf algebra defined on $\tB$} 
	\label{subsec:iHopf:ff}
	
	Recall the Cartan matrix $C=(c_{ij})$ and $D=\diag(d_i\mid i\in\I)$. Let $'\ff$ be the free associative $\Q(v^{\frac12})$-algebra with generators $\theta_i$ ($i\in\I$); see \cite[Chap.~ 1]{Lus93}. We denote a rescaled version of $\theta_i$ by 
	\begin{align} \label{rescale:theta}
		\vartheta_i=(v_i -v_i^{-1}) \theta_i.
	\end{align}
	
	Let $'\hB$  be the $\Q(v^{\frac12})$-algebra generated by
	$\vartheta_i,\h_i$ $(i\in\I)$ subject to 
	\begin{align*}
		[\h_i,\h_j]=0, \qquad h_i\vartheta_j=v_i^{c_{ij}} \vartheta_j h_i.  
	\end{align*}
	
	Let $\ff$ (resp. $\hB$) be the quotient algebra of $'\ff$ (resp. $'\hB$) by the ideal generated by
	\begin{align}
		& \sum_{r=0}^{1-c_{ij}} (-1)^r \left[ \begin{array}{c} 1-c_{ij} \\r \end{array} \right]_{v_i}  \vartheta_i^r \vartheta_j  \vartheta_i^{1-c_{ij}-r},\qquad \forall i\neq j.
	\end{align}
	We endow $\ff$ with an $\N^\I$-grading by setting $\wt(\vartheta_i)=\alpha_i$. Let $\ff_\mu$ be the homogeneous subspace of degree $\mu$. Then $\ff=\bigoplus_{\mu\in\N^\I}\ff_\mu$. 
	
	Let $r:\ff\to\ff\otimes\ff$ be the homomorphism defined by Lusztig \cite[1.2.6]{Lus93}. In Sweedler notation, we write $r(x)=\sum x_{(1)}\otimes x_{(2)}$ for any $x\in\ff$. Then the coproduct of $\hB$ satisfies
	\begin{align*}
		\Delta(x)=\sum x_{(1)}h_{\wt(x_{(2)})}\otimes x_{(2)}, \qquad \Delta^2(x)=\sum x_{(1)}h_{\wt(x_{(2)})}h_{\wt(x_{(3)})}\otimes x_{(2)}h_{\wt(x_{(3)})}\otimes x_{(3)}.
	\end{align*}
	This convention greatly improves the clarity of the computation and will be adopted throughout this paper.
	
	We identify 
	\begin{align} \label{ffUU}
		\ff \stackrel{\cong}{\longrightarrow}\U^+, \,  \vartheta_i \mapsto \vartheta_i^+ :=E_i,  \qquad
		\ff \stackrel{\cong}{\longrightarrow}\U^-, \,\vartheta_i \mapsto \vartheta_i^- := F_i. 
	\end{align}
	Then we have  $\hB\cong \U^+\otimes \Q(v^{\frac12})[K_i\mid i\in\I]$. 
	
	Let $\tB$ (resp.  $'\tB$) be the algebra constructed from $\hB$ (resp. $'\hB$) with $h_i$ invertible for $i\in\I$. Define the coproduct and counit  
	\begin{align}
		\Delta(\vartheta_i)= \vartheta_i\otimes 1 +  h_i\otimes \vartheta_i, \qquad
		\Delta(\h_{i}) = \h_{i} \otimes \h_{i},\quad \forall i\in\I;
		\\
		\varepsilon(\vartheta_i)= 0, \qquad
		\varepsilon(\h_{i}) = 1=\varepsilon(\h_i^{-1}),\quad \forall i\in\I.
	\end{align}
	In this way, ${'\hB},\hB,{'\tB},\tB$ are all Hopf algebras. Define 
	\begin{align}
		\label{eq:hopf-pairing}
		\varphi(\vartheta_i,\vartheta_j)=\delta_{ij}(v_i-v_i^{-1}),\quad \varphi(\h_i,\h_j)=v_i^{c_{ij}},\quad\varphi(\vartheta_i,\h_j)=0,\quad \forall i,j\in\I,x,y\in\ff.
	\end{align}
	Then it gives (symmetric) Hopf pairings on the Hopf algebras $'\hB,\hB,'\tB,$ and $\tB$. Moreover, these Hopf pairings are non-degenerate on $\hB$ and $\tB$.
	
	Let  $(\hB\otimes\hB )^\imath$ and $(\tB\otimes\tB )^\imath$ be the iHopf algebras of diagonal type, associated to $(\hB,\varphi)$ and $(\tB,\varphi)$, respectively; see \cite[\S4.1]{CLPRW}. In fact, $(\hB\otimes\hB )^\imath$ is defined on the same vector space as $\hB\otimes \hB$ equipped with a new multiplication
	\[(a\otimes b)\ast (c\otimes d)=\sum\varphi (a_{(1)},d_{(2)})\cdot \varphi(c_{(2)},b_{(1)})\cdot a_{(2)}c_{(1)}\otimes b_{(2)}d_{(1)},\qquad\forall a,b,c,d\in \hB.
	\]
	Similar for $(\tB\otimes\tB )^\imath$. 
	
	\begin{lemma}[{\cite[Lemma 4.1]{CLPRW}}]\label{lem:U=iHopfBB}
		We have Hopf algebra isomorphisms 
		\begin{align*}
			\widehat{\Phi}_{\sharp}:&\hU\longrightarrow (\hB\otimes\hB )^\imath ,\qquad \widetilde{\Phi}_{\sharp}:\tU\longrightarrow (\tB\otimes\tB )^\imath,\\
			E_i\mapsto\vartheta_i\otimes 1,&\quad F_i\mapsto 1\otimes \vartheta_i,\quad K_i \mapsto  h_i\otimes1,\quad K_i'\mapsto 1\otimes h_i,\quad \forall i\in\I.
		\end{align*}
	\end{lemma}

	Let $\tau$ be an involution in $\Inv(C)$. Clearly $\tau$ preserves the Hopf pairing $\varphi$. By the construction of iHopf algebras in \cite{CLPRW}, we denote by
	\[
	\hB^\imath_\tau = \text{iHopf}\, \big(\hB, \varphi \circ (\tau \otimes 1)\big), \qquad
	\tB^\imath_\tau =\text{iHopf}\, \big(\tB, \varphi\circ (\tau \otimes 1)\big)
	\]
	the iHopf algebras defined on $\big(\hB, \varphi \circ (\tau \otimes 1)\big)$ and $\big(\tB, \varphi\circ (\tau \otimes 1)\big)$, respectively. In fact, $\hB^\imath_\tau$ is the same vector space as $\hB$ equipped with a new multiplication: 
	\begin{align}\label{star product}
		a\ast b:=\sum \varphi(\tau b_{(2)},a_{(1)})\cdot a_{(2)}b_{(1)},\quad \forall a,b\in \hB,
	\end{align}
	where $\Delta(a)=\sum a_{(1)}\otimes a_{(2)}$, $\Delta(b)=\sum b_{(1)}\otimes b_{(2)}$. Similar for $\tB^\imath_\tau.$
	
	\begin{theorem}[{\cite[Theorem 4.4]{CLPRW}}]\label{thm:iH=iQG}
		We have algebra isomorphisms 
		\begin{align*}
			\widehat{\Phi}^\imath: \hB^\imath_\tau&\longrightarrow \hUi,\qquad \widetilde{\Phi}^\imath:\tB^\imath_\tau\longrightarrow \tUi, \\
			\vartheta_i &\mapsto B_i,\qquad h_i\mapsto \tk_{\tau i},\quad \forall i\in\I.
		\end{align*}
	\end{theorem}
	
	In the following, we always identify $\hUi \equiv\hB^\imath_\tau$, and $\tUi \equiv\tB^\imath_\tau$.
	
	Denote by $\chi: \hB\rightarrow \Q(v^{\frac{1}{2}})$ (respectively, $\chi: \tB\rightarrow \Q(v^{\frac{1}{2}})$) the $\tau$-twisted compatible linear map given in \cite[Lemma 4.2]{CLPRW}. That is, $\chi: \hB\rightarrow \Q(v^{\frac{1}{2}})$ is the linear map such that $\chi(1)=1$, 
	$\chi(\vartheta_i)=0$, $\chi(h_i)=\varphi(h_i,h_{\tau i})$, for $i\in\I$, and  $\chi(ab)=\sum\chi(a_{(1)})\chi(b_{(2)})\varphi(\tau (a_{(2)}),b_{(1)})$ holds for all $a,b\in \hB$. (The same statement holds  when replacing $\hB$ by $\tB$.)

	\begin{lemma}[{\cite[Lemma 4.3, Theorem 4.4]{CLPRW}}] \label{lem:iHembed}
		There are algebra homomorphisms 
		\begin{equation} \label{eq:xi_theta embedding of iHopf}
			\widehat{\xi}_\tau: \hB^\imath_\tau\longrightarrow(\hB\otimes\hB)^\imath,
			\qquad 
			\widetilde{\xi}_\tau: \tB^\imath_\tau\longrightarrow(\tB\otimes\tB)^\imath
		\end{equation}
		which send $a \mapsto \sum \chi(a_{(2)})\cdot \tau(a_{(3)})\otimes a_{(1)}$. In particular, we have the following commutative diagrams
		\begin{equation} \label{CD:fU}
			\begin{tikzcd}
				\hB^\imath_\tau\ar[r,"\widehat{\xi}_\tau"]\ar[d,dashed,swap,"\widehat{\Phi}^\imath"]&(\hB\otimes\hB)^\imath\ar[d,"\widehat{\Phi}_\sharp^{-1}"]\\
				\hUi\ar[r,hook]&\hU
			\end{tikzcd}
			\qquad \qquad 
			\begin{tikzcd}
				\tB^\imath_\tau\ar[r,"\widetilde{\xi}_\tau"]\ar[d,dashed,swap,"\widetilde{\Phi}^\imath"]&(\tB\otimes\tB)^\imath\ar[d,"\widetilde{\Phi}_\sharp^{-1}"]\\
				\tUi\ar[r,hook]&\tU
			\end{tikzcd}
		\end{equation}
	\end{lemma}
	
	\subsection{A recursive formula and $\diamond$-action}
	
	For the algebra $\ff$, there exist linear maps known as skew-derivations (cf. \cite{Lus93}) 
	\[
	\partial_i^R:\ff\longrightarrow \ff, \qquad
	\partial_i^L:\ff\longrightarrow \ff
	\]
	such that $\partial_i^R(1)=\partial_i^L(1)=0$, $\partial_i^R(\vartheta_j)=\delta_{ij}=\partial_i^L(\vartheta_j)$, and
	\begin{align*}
		&\quad    
		\partial_i^R(fg)=\partial_i^R(f)g+v^{(\alpha_i,\mu)}f\partial_i^R(g),
		\\
		&\quad \partial_i^L(fg)=v^{(\alpha_i,\nu)}\partial_i^L(f)g+f\partial_i^L(g),
	\end{align*}
	for any $j\in\I, f\in\ff_\mu, g\in\ff_\nu$.

	Recall the two algebras $(\widehat\BB, \cdot)$ and $(\hB^\imath_\tau, *)$ have the same underlying vector space (which contains $\ff$ as a subspace). 
	
	\begin{lemma}[{\cite[Lemma 4.5]{CLPRW}}]
		\label{lem:recursive}
		In $\hB^\imath_\tau$ (and $\tB^\imath_\tau$), for $x\in\ff$ and $i\in \I$, we have
		\begin{align*}
			\vartheta_i*x&=\vartheta_i\cdot x+(v_i-v_i^{-1})\partial_{\tau i}^L(x)\cdot h_{\tau i},
			\\
			x*\vartheta_i&=x\cdot \vartheta_i+(v_i-v_i^{-1})\partial_{\tau i}^R(x)\cdot h_{i}.
		\end{align*}
	\end{lemma}

	\begin{lemma}
		\label{lem:bar}
		There exists a bar involution $\ov{\phantom{x}}$ on $\hB^\imath_\tau$ (also on $\tB^\imath_\tau$), which is an anti-involution of $\Q$-algebra such that
		\[
		\ov{v^{1/2}}=v^{-1/2}, \quad
		\ov{\vartheta_i}=\vartheta_i,
		\quad    \ov{h_\alpha}=v^{(\alpha,\tau\alpha)}h_\alpha, 
		\qquad \text{ for } i\in\I, \alpha\in\N^\I.
		\]
	\end{lemma}
	
	\begin{proof}
		Follows from Lemma \ref{lem:involution-iQG}~(2) by using the isomorphisms in Theorem~\ref{thm:iH=iQG}.
	\end{proof}

	Let $\widetilde{\mathcal{T}}$ be the subalgebra of $\tB^\imath_\tau$ generated by $h_\alpha$, $\alpha\in\Z^\I$, which is a Laurent polynomial algebra in $h_i$, for $i\in \I$. Similarly, one can define the subalgebra $\widehat{\mathcal{T}}$ of $\hB^\imath_\tau$, which is a polynomial algebra in $h_i$, for $i\in \I$. 
	We define a $\diamond$-action of $\widehat{\mathcal{T}}$ on $\hB^\imath_\tau$ by letting
	\begin{align}  \label{eq:diamond}
		h_\alpha\diamond x :=v^{\frac{1}{2}(\tau\alpha-\alpha,\wt(x))}\cdot h_\alpha\ast x,
	\end{align}
	for $\alpha\in\N^\I$ and $x\in\hB$. 
	The $\diamond$-action of $\widetilde{\mathcal{T}}$ on $\tB^\imath_\tau$ is defined similarly. 
	
	\begin{lemma}
		We have $\ov{h_\alpha\diamond x}=\ov{h_\alpha}\diamond\ov{x}$, for $\alpha\in\N^\I$ and $x\in\hB$.
	\end{lemma}
	
	\begin{proof}
		Assume that $x$ is homogeneous, and note that $\wt(x)=\wt(\ov{x})$. Applying the bar involution defined in Lemma~\ref{lem:bar}, we have 
		\begin{align*}
			\ov{h_\alpha\diamond x} &= v^{-\frac{1}{2}(\tau\alpha-\alpha,\wt(x))}v^{(\alpha,\tau \alpha)}\ov{x}*h_\alpha
			\\
			&= v^{-\frac{1}{2}(\tau\alpha-\alpha,\wt(x))}v^{(\alpha,\tau \alpha)} v^{(\tau \alpha-\alpha,\wt(x))} h_\alpha*\ov{x}
			\\
			&= v^{\frac{1}{2}(\tau\alpha-\alpha,\wt(x))}v^{(\alpha,\tau \alpha)} h_\alpha*\ov{x}
			=\ov{h_\alpha}\diamond \ov{x}.
		\end{align*}
		The lemma is proved.    
	\end{proof}
	
	For $\alpha\in\N^\I$, note that $$\K_\alpha=v^{\frac{1}{2}(\alpha,\tau\alpha)}h_{\tau\alpha}$$
	is the unique bar-invariant element in $\widehat{\mathcal{T}}$ of the form $\K_\alpha=\lambda h_{\tau\alpha}$ where $\lambda\in v^{\frac12 \Z}$.

	\section{$\mathrm{i}$Braid group symmetries on $\tUi$ via $\mathrm{i}$Hopf algebra} 
	\label{sec:braid}
	
	In this section, we shall give an iHopf algebra interpretation of the ibraid group action on iquantum groups, providing a new connection to braid group action of quantum groups. 
	
	\subsection{Connecting 2 braid group actions via iHopf}
	
	Let $j\neq i,\tau i$ in $\I$ in this section. Define the root vectors in $\ff$:
	\begin{align}
		f_{i,j;m}=\sum_{r+s=m}(-1)^rv_i^{r(c_{ij}+m-1)+\frac{1}{2}m}(v_i-v_i^{-1})^{-m}\vartheta_i^{(r)}\vartheta_j\vartheta_i^{(s)},\label{eq:root vector i=taui-1}
		\\
		f'_{i,j;m}=\sum_{r+s=m}(-1)^rv_i^{r(c_{ij}+m-1)+\frac{1}{2}m}(v_i-v_i^{-1})^{-m}\vartheta_i^{(s)}\vartheta_j\vartheta_i^{(r)},\label{eq:root vector i=taui-2}
	\end{align}
	which are slightly normalized versions of Lusztig's definition \cite{Lus93}. By the same proof of \cite[Proposition 37.2.5]{Lus93}, the rescaled braid group symmetry $\widetilde{T}_i$ satisfies
	\begin{align}\label{eq: tU T_i action on f_m}
		\widetilde{T}_i(f_{i,j;m}^+)=f_{i,j,-c_{ij}-m}'^+,\qquad \forall m,n\in\Z.
	\end{align}
	
	Let $\ad:\ff\to\ff$ be the adjoint action via the identification $\ff\cong\U^+$; it is given by
	\begin{align}\label{eq: ad action def}
		\ad(\vartheta_i)(x)=\vartheta_ix-h_ixh_i^{-1}\vartheta_i.
	\end{align}
	Recalling the anti-involution $\sigma$ from Lemma~\ref{lem:anti-involut-QG}, we have
	\begin{align} \label{rootvector2}
		f_{i,j;m}=v_i^{\frac{1}{2}m}(v_i-v_i^{-1})^{-m}\sigma\big(\ad(\vartheta_i^{(m)})(\vartheta_j)\big),
		\quad 
		f'_{i,j;m}=v_i^{\frac{1}{2}m}(v_i-v_i^{-1})^{-m}\ad(\vartheta_i^{(m)})(\vartheta_j).
	\end{align}
	
	The ibraid group symmetry $\tTT_i$ of $\tUi$ can be transported to $\tB^\imath_\tau$ via the isomorphism $\widehat{\Phi}^\imath:\tB^\imath_\tau \stackrel{\cong}{\rightarrow} \tUi$ in Theorem \ref{thm:iH=iQG}. Recalling $\K_\alpha\in \tUi$ from \eqref{eq:bbKi}, we define
	\begin{align}
		\K_\alpha:=v^{\frac{1}{2}(\alpha,\tau\alpha)}h_{\tau\alpha} \in \tB^\imath_\tau
	\end{align}
	so that $\widehat{\Phi}^\imath(\K_\alpha)=\K_\alpha$. Then $\tTT_i(\K_\alpha)=\K_{r_i\alpha}$ in $\tB^\imath_\tau$; see  \eqref{def:ri} for definition of $r_i$. Denote by 
	\begin{align} \label{embed:f}
		\iota:\ff \longrightarrow \tB^\imath_\tau
	\end{align}
	the canonical embedding. 
	

	\subsection{Subalgebras $\ff[i,\tau i]$ and ${^\sigma\ff}[i,\tau i]$}
	
	For $i\in\I$, we let $\hU_{i,\tau i}$ be the subalgebra of $\hU$ generated by $K_j,K_j',E_j,F_j$, for $j\in\{i,\tau i\}$. Also, let $\ff_{i,\tau i}$ be the subalgebra of $\ff$ generated by $\vartheta_i,\vartheta_{\tau i}$. That is,
	\begin{align} \label{eq:firank1}
		\hU_{i,\tau i} =\Q(v^{\frac12}) \langle K_j,K_j',E_j,F_j \mid j\in\{i,\tau i\} \rangle
		;
		\qquad
		\ff_{i,\tau i} =\Q(v^{\frac12}) \langle \vartheta_i,\vartheta_{\tau i} \rangle.
	\end{align}
	We further set
	\begin{align} \label{eq:fitaui}
		\ff[i,\tau i] :=\{x\in\ff\mid \widetilde{T}_{r_i}(x^+)\in\U^+\},
		\qquad 
		{^\sigma\ff}[i,\tau i] :=\{x\in\ff\mid \widetilde{T}_{r_i}^{-1}(x^+)\in\U^+\},
	\end{align}
	which are subalgebras of $\ff$ since $\widetilde{T}_{r_i}$ is an algebra homomorphism. 
	Recall the anti-involution $\sigma$ on $\tU$ from Lemma \ref{lem:anti-involut-QG}. Since $\widetilde{T}_{r_i}^{-1}=\sigma\circ\widetilde{T}_{r_i}\circ\sigma$, we see that ${^\sigma\ff}[i,\tau i]=\sigma(\ff[i,\tau i])$ and $\sigma$ induces anti-isomorphisms between $\ff[i,\tau i]$ and ${^\sigma\ff}[i,\tau i]$. 
	
	In the following subsections, we will focus on the subalgebra ${^\sigma\ff}[i,\tau i]$ and study the relative braid group action on it. To that end, we first specify a generating set of ${^\sigma\ff}[i,\tau i]$, which is provided by certain root vectors which are studied in depth in \cite[\S4]{CLPRW}. 
	\begin{itemize}
		\item If $c_{i,\tau i}=2$, then the root vectors $f_{i,j;m}$, $f'_{i,j;m}$ are given by \eqref{eq:root vector i=taui-1}--\eqref{eq:root vector i=taui-2} or \eqref{rootvector2}.
		\item 
		If $c_{i,\tau i}=0$, then for $m,n\in\Z$ we set 
		\begin{align}\label{rootvector0}
			\begin{split}
				f_{i,\tau i,j;m,n}&=v_i^{\frac{1}{2}(m+n)}(v_i-v_i^{-1})^{-(m+n)}\sigma(\ad(\vartheta_i^{(m)}\vartheta_{\tau i}^{(n)})(\vartheta_j)),\\
				f'_{i,\tau i,j;m,n}&=v_i^{\frac{1}{2}(m+n)}(v_i-v_i^{-1})^{-(m+n)}\ad(\vartheta_i^{(m)}\vartheta_{\tau i}^{(n)})(\vartheta_j).
			\end{split}
		\end{align}
		\item 
		If $c_{i,\tau i}=-1$, then for $a,b,c\in\Z$ we set 
		\begin{align} \label{rootvector-1}
			\begin{split}
				f_{i,\tau i,j;a,b,c}&=v_i^{\frac{1}{2}(a+b+c)}(v_i-v_i^{-1})^{-(a+b+c)}\sigma(\ad(\vartheta_i^{(a)}\vartheta_{\tau i}^{(b)}\vartheta_i^{(c)})(\vartheta_j)),\\
				f'_{i,\tau i,j;a,b,c}&=v_i^{\frac{1}{2}(a+b+c)}(v_i-v_i^{-1})^{-(a+b+c)}\ad(\vartheta_i^{(a)}\vartheta_{\tau i}^{(b)}\vartheta_i^{(c)})(\vartheta_j).
			\end{split}
		\end{align}
	\end{itemize}
	
	Recall $\I_i=\{i,\tau i\}$. Now define a subset ${^\sigma\!R_{i,\tau i}}$ of ${^\sigma\ff}[i,\tau i]$ by 
	\[
	{^\sigma\!R_{i,\tau i}}:=
	\begin{dcases}
		\{f'_{i,j;m}\mid m\in\Z,j\notin\I_i\}&\text{if $i=\tau i$},\\
		\{f'_{i,\tau i,j;m,n}\mid m,n\in\Z,j\notin\I_i\}&\text{if $c_{i,\tau i}=0$},\\
		\{f'_{i,\tau i,j;a,b,c},f'_{\tau i,i,j;a,b,c}\mid a,b,c\in\Z,j\notin\I_i\}&\text{if $c_{i,\tau i}=-1$}.
	\end{dcases}
	\]
	
	The following statement seems well known; in case when $i=\tau i$ this can be found in \cite[38.1.2, 38.1.6]{Lus93}. It can be obtained from \cite[Conjecture~ 5.3]{BG17a}, which was proved independently in \cite[Proposition 2.10]{Tan17} and \cite[Theorem 1.1]{Kim17}. 
	It is also given in \cite[Corollary 2.3, (2.12)]{KY21} who refers back to \cite{Rad85}.  
	
	\begin{proposition}
		\label{prop:fract-U}
		The algebra ${^\sigma\ff}[i,\tau i]$ is generated by the set ${^\sigma\!R_{i,\tau i}}$. Moreover, the following multiplication maps are linear isomorphisms
		\begin{align}\label{eq: f[i,tau i] linear iso}
			{^\sigma\ff}[i,\tau i]\otimes\ff_{i,\tau i}\stackrel{\cong}{\longrightarrow} \ff,\qquad 
			\ff_{i,\tau i}\otimes {^\sigma\ff}[i,\tau i]\stackrel{\cong}{\longrightarrow} \ff.
		\end{align}
	\end{proposition}

	\subsection{Relative braid group action on ${^\sigma\ff}[i,\tau i]$}
	
	Recall from \eqref{eq:xi_theta embedding of iHopf} the embedding $\widetilde{\xi}_\tau:\tB^\imath_\tau\to (\tB\otimes\tB)^\imath \equiv \tU$, where $(\tB\otimes\tB)^\imath$ is identified with $\tU$ via the isomorphism in Lemma~\ref{lem:U=iHopfBB}.

	\begin{lemma}\label{lem: xi_theta leading term}
		For $x\in\ff$, we have 
		\[\widetilde{\xi}_\tau(x)\in x^-+\sum_{\alpha\in\N^\I\setminus\{0\}}\U^+\U^-\cdot K_\alpha'.\]
	\end{lemma}
	\begin{proof}
		For $x\in\ff$, the definition of $\widetilde{\xi}_\tau(x)$ becomes
		\[
		\widetilde{\xi}_\tau(x)=\sum\chi(x_{(2)}h_{\wt(x_{(3)})})\tau (x_{(3)})\otimes x_{(1)}h_{\wt(x_{(2)})}h_{\wt(x_{(3)})}.
		\]
		Recall that we identify $(\tB\otimes\tB)^\imath$ with $\tU$.
		We note that for $a,b\in\ff$, the element $a\otimes b$ lies in $\sum_{\alpha\in\N^\I}\U^+\U^-\cdot K_\alpha'$ by induction on $\wt(a)$ and $\wt(b)$ since 
		\[(a\otimes 1)*(1\otimes b)=\sum\varphi(a_{(1)},b_{(2)})a_{(2)}\otimes b_{(1)}=a\otimes b+\sum_{b_{(2)}\neq 1}\varphi(a_{(1)},b_{(2)})a_{(2)}\otimes b_{(1)}.\]
		Now if $h_{\wt(x_{(2)})}h_{\wt(x_{(3)})}\neq 1$ then $\tau (x_{(3)})\otimes x_{(1)}h_{\wt(x_{(2)})}h_{\wt(x_{(3)})}\in \sum_{\alpha\in\N^\I\setminus\{0\}}\U^+\U^-\cdot K_\alpha'$. So the claim follows.
	\end{proof}
	
	Recall the braid group symmetry $\widetilde{\mathscr{T}}_{r_i}$ on $\tU$, cf. \eqref{eq: T_i twisted def}. Note that $\widetilde{\mathscr{T}}_{r_i}$ coincides with $\widetilde{T}_{r_i}$ on $\U^-$. Denote 
	\begin{gather*}
		\Lambda_{i,\tau i}   :=\N\alpha_i+\N\alpha_{\tau i},
		\qquad 
		\Lambda[i,\tau i] :=\{\beta\in\N^\I\mid r_i(\beta)\in\N^\I\},
		\\
		\Xi_i := \{r_i(\beta)-\eta\mid \beta\in\Lambda[i,\tau i],\eta\in\Lambda_{i,\tau i}\}.
	\end{gather*}
	In the following, we identify $\tU$ with $(\tB\otimes\tB)^\imath$ via the isomorphism  $\widetilde{\Phi}_\sharp$ in Lemma~ \ref{lem:U=iHopfBB}. For $x,y\in\tB$, we view $x\otimes y\in\tB\otimes \tB$ as an element in $(\tB\otimes\tB)^\imath=\tU$.
	
	\begin{lemma}\label{lem:T_ri apply on tensor of x}
		For any $\alpha_0\in\N^\I$, $\beta_0\in\Lambda[i,\tau i]$, $x\in{^\sigma\ff}[i,\tau i]$, $y\in\ff$, by viewing $x\otimes y\in\tU$ we have
		\[
		\widetilde{\mathscr{T}}_{r_i}^{-1}((x\otimes y)K'_{\alpha_0+\beta_0})\in\sum_{\gamma\in\Xi_i}\U^+\U^-\cdot K'_\gamma.
		\]
	\end{lemma}
	
	\begin{proof}
		For $y=1$ the claim is clear. Now we prove the general case by induction on $\wt(y)$. Note that
		\begin{align*}
			(x\otimes y)K_{\alpha_0+\beta_0}'&=(x\otimes 1)*(1\otimes y)K'_{\alpha+\beta}-\sum\varphi(x_{(1)},y_{(2)})(x_{(2)}\otimes y_{(1)}h_{\wt(y_{(2)})})K'_{\alpha_0+\beta_0}\\
			&\in(x\otimes 1)*(1\otimes y)K'_{\alpha_0+\beta_0}+\sum\Q(v)(x_{(2)}\otimes y_{(1)})K'_{\wt(y_{(2)})+\alpha_0+\beta_0}.
		\end{align*}
		For the leading term $(x\otimes 1)*(1\otimes y)K'_{\alpha_0+\beta_0}$, we may assume that $y=uy'$ where $u\in\ff_{i,\tau i}$ and $y'\in{^\sigma\ff}[i,\tau i]$. Then
		\[\widetilde{\mathscr{T}}_{r_i}^{-1}((x\otimes 1)*(1\otimes y)K'_{\alpha_0+\beta_0})=\widetilde{\mathscr{T}}_{r_i}^{-1}(x\otimes 1)* \widetilde{\mathscr{T}}_{r_i}^{-1}(1\otimes u)*\widetilde{\mathscr{T}}_{r_i}^{-1}(1\otimes y')K'_{r_i(\alpha_0+\beta_0)}\]
		where $\widetilde{\mathscr{T}}_{r_i}^{-1}(x\otimes 1)\in\U^+$, $\widetilde{\mathscr{T}}_{r_i}^{-1}(1\otimes y')\in\U^-$, and $\widetilde{\mathscr{T}}_{r_i}^{-1}(1\otimes u)\in\sum_{\alpha\in\Lambda_{i,\tau i}}\tU_{i,\tau i}^+K_\alpha'^{-1}$. Thus the left-hand side belongs to $\sum_{\gamma\in\Xi_i}\U^+\U^-\cdot K'_\gamma$. Note that 
		\[
		r({^\sigma\ff}[i,\tau i])\subseteq\ff\otimes{^\sigma\ff}[i,\tau i]; 
		\]
		this follows from the identity ${^\sigma\ff}[i,\tau i]=\{x\in\ff:\widetilde{T}_i^{-1}(x^+)\in\U^+,\widetilde{T}_{\tau i}^{-1}(x^+)\in\U^+\}$ (see \cite[Lemma 2.4]{KY21}) and then applying \cite[38.1.8]{Lus93}. By the induction hypothesis, the other terms $\widetilde{\mathscr{T}}_{r_i}^{-1}((x_{(2)}\otimes y_{(2)})K'_{\wt(y_{(1)})+\alpha_0+\beta_0})$ also belong to $\sum_{\gamma\in\Xi_i}\U^+\U^-\cdot K'_\gamma$. This completes the induction step.
	\end{proof}
	
	\begin{corollary}
		For any $x\in{^\sigma\ff}[i,\tau i]$, we have
		\begin{align}
			\label{eq:T_i on f[i,tau i] K-matrix conjugate-1}
			\begin{split}
				\widetilde{\xi}_\tau(\tTT_i^{-1}(x))\widetilde{\Upsilon}_i
				&=\widetilde{\Upsilon}_i\widetilde{\mathscr{T}}_{r_i}^{-1}(\widetilde{\xi}_\tau(x))
				\\
				& \in \widetilde{\Upsilon}_i\widetilde{T}_{r_i}^{-1}(x^-)+\sum_{\gamma\in\Xi_i\setminus\{0\}}\U^+\U^-\cdot K'_\gamma.
			\end{split}
		\end{align}
	\end{corollary}
	
	\begin{proof}
		Note that for $\gamma=r_i(\beta)-\eta\in\Xi_i$, we have $\gamma=0$ if and only if $\beta=\eta=0$. Therefore Lemma~\ref{lem:T_ri apply on tensor of x} implies
		\begin{equation}\label{eq:T_ri on f[i,tau i] leading term-1}
			\widetilde{\mathscr{T}}_{r_i}^{-1}(\widetilde{\xi}_\tau(x))\in \widetilde{T}_{r_i}^{-1}(x^-)+\sum_{\gamma\in\Xi_i\setminus\{0\}}\U^+\U^-\cdot K'_\gamma,\quad \forall x\in{^\sigma\ff}[i,\tau i].
		\end{equation}
		
		Now the statement \eqref{eq:T_i on f[i,tau i] K-matrix conjugate-1} follows by identifying $\tUi$ with $\tB^\imath_\tau$ via the embedding $\widetilde{\xi}_\tau$ and  applying Theorem~\ref{thm: relative T_i conjugate}. 
	\end{proof}
	
	\begin{lemma}\label{lem: T_i on f[i,taui] cartan part}
		For any $x\in{^\sigma\ff}[i,\tau i]$, we have
		\[
		\tTT_i^{-1}(x)\in \sum_{\alpha\in\N^\I-\Lambda_{i,\tau i}} h_{\alpha}\ast \ff.
		\]
	\end{lemma}
	
	\begin{proof}
		The statement is clear for the generators of ${^\sigma\ff}[i,\tau i]$ from the results in \cite{CLPRW} (see \cite[Propositions 5.5,5.10,5.15]{CLPRW}). Now assuming that the statement holds for $x,y\in{^\sigma\ff}[i,\tau i]$, we shall prove that it holds for the product $xy$. 
		
		To that end, we have 
		\begin{align*}
			\tTT_i^{-1}(x\ast y)=\tTT_i^{-1}(x)\ast\tTT_i^{-1}(y)\in\sum_{\alpha\in \N^\I-\Lambda_{i,\tau i}}h_\alpha\ast \ff.
		\end{align*}
		On the other hand, we have by \eqref{star product}  
		\begin{align*}
			\tTT_i^{-1}(x\ast y)&=\sum\varphi(x_{(1)},y_{(2)})\tTT_i^{-1}(x_{(2)}y_{(1)}h_{\wt(y_{(2)})})\\
			&=\tTT_i^{-1}(xy)+\sum_{y_{(2)}\neq 1}\varphi(x_{(1)},y_{(2)})\tTT_i^{-1}(x_{(2)}y_{(1)}h_{\wt(y_{(2)})}).
		\end{align*}
		Note that $x_{(2)},y_{(2)}\in{^\sigma\ff}[i,\tau i]$, and thus $\wt(y_{(2)})\in\Lambda[i,\tau i]$. 
		Since $\tTT_i^{-1}(x_{(2)}y_{(1)})\in\sum_{\alpha\in\N^\I-\Lambda_{i,\tau i}}h_{\alpha}\ast \ff$, we have $\tTT_i^{-1}(x_{(2)}y_{(1)}h_{\wt(y_{(2)})})\in \sum_{\alpha\in \N^\I-\Lambda_{i,\tau i}}h_{\alpha}\ast \ff$. Therefore, the statement in the lemma holds for $xy$.
	\end{proof}
	
	\begin{corollary}
		For any $x\in{^\sigma\ff}[i,\tau i]$, we have
		\begin{align}
			\label{eq:T_i on f[i,tau i] leading term}
			\begin{split}
				\widetilde{\xi}_\tau\big(\tTT_i^{-1}(x)\big) &=\widetilde{\Upsilon}_i\widetilde{\mathscr{T}}_{r_i}^{-1}\big(\widetilde{\xi}_\tau(x)\big) \widetilde{\Upsilon}^{-1}_i
				\\
				&\in \widetilde{\Upsilon}_i\widetilde{T}_{r_i}^{-1}(x^-)\widetilde{\Upsilon}^{-1}_i+\sum_{\delta,\eta\in\Lambda_{i,\tau i},\gamma\in\Xi_i\setminus(-\Lambda_{i,\tau i})}\U^+\U^-\cdot K_\delta K'_{\gamma+\eta}.
			\end{split}
		\end{align}
	\end{corollary}
	
	\begin{proof}
		By Lemma~\ref{lem: T_i on f[i,taui] cartan part}, for any $x\in{^\sigma\ff}[i,\tau i]$, we can write
		\[\tTT_i^{-1}(x)\in\sum_{\alpha\in \N^\I-\Lambda_{i,\tau i},y\in\ff}a_{\alpha,y}h_\alpha\ast y.\]
		Then Lemma~\ref{lem: xi_theta leading term} gives us
		\begin{align*}
			\widetilde{\xi}_\tau(\tTT_i^{-1}(x))\widetilde{\Upsilon}_i&=\Big(\sum_{\alpha\in \N^\I-\Lambda_{i,\tau i},y\in\ff}a_{\alpha,y}K_{\tau\alpha}K_{\alpha}'\cdot \widetilde{\xi}_\tau(y)\Big)\widetilde{\Upsilon}_i\\
			&\in \sum_{\alpha\in \N^\I-\Lambda_{i,\tau i}}\Big(\U^+\U^-\hU_{i,\tau i}^{0} K_{\tau\alpha}K_\alpha'+\sum_{\beta\in\N^\I\setminus\{0\}}\U^+\U^-\hU_{i,\tau i}^{0} K_{\tau\alpha}K_{\alpha+\beta}'\Big),
		\end{align*}
		where $\hU^{0}_{i,\tau i}$ is the subalgebra of $\hU_{i,\tau i}$ generated by $K_i,K_{\tau i}, K_{i}', K_{\tau i}'$. Comparing this with \eqref{eq:T_i on f[i,tau i] K-matrix conjugate-1}, we find that $a_{\alpha,y}=0$ for $\alpha\neq 0$, i.e. $\tTT_i^{-1}(x)\in\ff$, and
		\begin{equation}
			\label{eq:T_ri on f[i,tau i] leading term-2}
			\widetilde{\mathscr{T}}_{r_i}^{-1}\big(\widetilde{\xi}_\tau(x)\big) \in \widetilde{T}_{r_i}^{-1}(x^-)+\sum_{\gamma\in\Xi_i\setminus(-\Lambda_{i,\tau i})}\U^+\U^-\cdot K'_\gamma.
		\end{equation}
		The desired statement \eqref{eq:T_i on f[i,tau i] leading term}
		follows from conjugating $\widetilde{\mathscr{T}}_{r_i}^{-1}\big(\widetilde{\xi}_\tau(x)\big)$ in \eqref{eq:T_ri on f[i,tau i] leading term-2} by $\widetilde{\Upsilon}_i$.
	\end{proof}
	
	The following is the main result of this subsection.
	\begin{theorem}
		\label{thm:Tifi}
		We have $\tTT_i^{-1}(x)=\widetilde{T}_{r_i}^{-1}(x)$, for any $x\in{^\sigma\ff}[i,\tau i]$.
	\end{theorem}
	
	\begin{proof}
		The statement is verified when $x$ is one of the generators of ${^\sigma\ff}[i,\tau i]$ in \cite{CLPRW} (see \cite[Propositions 5.5,5.10,5.15]{CLPRW}). Now assuming the statement holds for $x,y\in{^\sigma\ff}[i,\tau i]$ we shall prove that it holds for $xy$. 
		
		To that end, applying \eqref{eq:T_i on f[i,tau i] leading term} to $xy$ gives us
		\begin{equation*}
			\begin{aligned}
				\widetilde{\Upsilon}_i\widetilde{T}_{r_i}^{-1}(x^-)\widetilde{\Upsilon}_i^{-1}\cdot \widetilde{\Upsilon}_i&\widetilde{T}_{r_i}^{-1}(y^-)\widetilde{\Upsilon}_i^{-1}=\widetilde{\Upsilon}_i\widetilde{T}_{r_i}^{-1}(x^-y^-)\widetilde{\Upsilon}_i^{-1}\\
				&\in \widetilde{\xi}_\tau\big(\tTT_i^{-1}(xy))+\sum_{\delta,\eta\in\Lambda_{i,\tau i},\gamma\in\Xi_i\setminus(-\Lambda_{i,\tau i})}\U^+\U^-\cdot K_\delta K'_{\gamma+\eta}.
			\end{aligned}
		\end{equation*}
		Together with Lemma \ref{lem: xi_theta leading term}, we know that
		\begin{align} 
			\label{eq: conjugate leading term}   
			\begin{split} &\widetilde{\Upsilon}_i\widetilde{T}_{r_i}^{-1}(x^-)\widetilde{\Upsilon}_i^{-1}\cdot \widetilde{\Upsilon}_i\widetilde{T}_{r_i}^{-1}(y^-)\widetilde{\Upsilon}_i^{-1}
				\\
				&\in \tTT_i^{-1}(xy)^-+
				\sum_{\alpha\in\N^\I\setminus\{0\}} \U^+\U^-\cdot K_\alpha'+
				\sum_{\delta,\eta\in\Lambda_{i,\tau i},\gamma\in\Xi_i\setminus(-\Lambda_{i,\tau i})}\U^+\U^-\cdot K_\delta K'_{\gamma+\eta}.
			\end{split}
		\end{align}
		
		On the other hand, as $\widetilde{\xi}_\vartheta$ and $\widetilde{\TT}_i$ are algebra homomorphisms, we have
		\begin{equation}\label{eq: T_i ast product leading term}
			\begin{aligned}
				&\widetilde{\xi}_\tau\big(\tTT_i^{-1}(x\ast y)\big)=\widetilde{\xi}_\tau(\tTT_i^{-1}(x))\cdot \widetilde{\xi}_\tau(\tTT_i^{-1}(y))\\
				&\in\Big(\widetilde{\Upsilon}_i\widetilde{T}_{r_i}^{-1}(x^-)\widetilde{\Upsilon}_i^{-1}+\sum_{\delta,\eta\in\Lambda_{i,\tau i},\gamma\in\Xi_i\setminus(-\Lambda_{i,\tau i})}\U^+\U^-\cdot K_\delta K'_{\gamma+\eta}\Big)\\
				&\quad \cdot\Big(\widetilde{\Upsilon}_i\widetilde{T}_{r_i}^{-1}(y^-)\widetilde{\Upsilon}_i^{-1}+\sum_{\delta',\eta'\in\Lambda_{i,\tau i},\gamma'\in\Xi_i\setminus(-\Lambda_{i,\tau i})}\U^+\U^-\cdot K_{\delta'} K'_{\gamma'+\eta'}\Big)\\
				&\subseteq \widetilde{\Upsilon}_i\widetilde{T}_{r_i}^{-1}(x^-)\widetilde{\Upsilon}_i^{-1}\cdot \widetilde{\Upsilon}_i\widetilde{T}_{r_i}^{-1}(y^-)\widetilde{\Upsilon}_i^{-1}+\sum_{\delta'',\eta''\in\Lambda_{i,\tau i},\gamma''\in\Xi_i\setminus(-\Lambda_{i,\tau i})}\U^+\U^-\cdot K_{\delta''}K'_{\gamma''+\eta''}.
			\end{aligned}
		\end{equation}
		Combining \eqref{eq: conjugate leading term} and \eqref{eq: T_i ast product leading term} then gives us
		\begin{equation}\label{eq:leading term of T_i^-1(x ast y)}
			\widetilde{\xi}_\tau(\tTT_i^{-1}(x\ast y))\in \tTT_i^{-1}(xy)^-+
			\sum_{\alpha\in\N^\I\setminus\{0\}} \U^+\U^-\cdot K_\alpha'+\sum_{\delta,\eta\in\Lambda_{i,\tau i},\gamma\in\Xi_i\setminus(-\Lambda_{i,\tau i})}\U^+\U^-\cdot K_\delta K'_{\gamma+\eta}.
		\end{equation}
		Now by the induction hypothesis $\tTT_i^{-1}(x\ast y)=\widetilde{T}_{r_i}^{-1}(x)\ast \widetilde{T}_{r_i}^{-1}(y)$, we have
		\begin{align*}
			&\widetilde{\xi}_\tau\big(\tTT_i^{-1}(x\ast y)\big) \\
			&=\widetilde{\xi}_\tau(\widetilde{T}_{r_i}^{-1}(x)\ast \widetilde{T}_{r_i}^{-1}(y))\\
			&=\widetilde{\xi}_\tau\Big(\sum\varphi(\widetilde{T}_{r_i}^{-1}(x)_{(1)},\widetilde{T}_{r_i}^{-1}(y)_{(2)})\widetilde{T}_{r_i}^{-1}(x)_{(2)}\widetilde{T}_{r_i}^{-1}(y)_{(1)}h_{\widetilde{T}_{r_i}^{-1}(y)_{(2)}}\Big)\\
			&=\widetilde{\xi}_\tau(\widetilde{T}_{r_i}^{-1}(x)\widetilde{T}_{r_i}^{-1}(y))+\widetilde{\xi}_\tau\Big(\sum_{\widetilde{T}_{r_i}^{-1}(y)_{(2)}\neq 1}\varphi(\widetilde{T}_{r_i}^{-1}(x)_{(1)},\widetilde{T}_{r_i}^{-1}(y)_{(2)})\widetilde{T}_{r_i}^{-1}(x)_{(2)}\widetilde{T}_{r_i}^{-1}(y)_{(1)}h_{\widetilde{T}_{r_i}^{-1}(y)_{(2)}}\Big)\\
			&\in \widetilde{T}_{r_i}^{-1}(x)^-\widetilde{T}_{r_i}^{-1}(y)^-+\sum_{\beta\in\N^\I\setminus\{0\}}\U^+\U^-\cdot K_{\beta}'+\sum_{\alpha\in\N^\I\setminus\{0\},\beta\in\N^\I\setminus\{0\}}\U^+\U^-\cdot K_{\tau\alpha}K_{\alpha+\beta}'
		\end{align*}
		where we have used Lemma~\ref{lem: xi_theta leading term}. Comparing the leading term of the right-hand side with the leading term of \eqref{eq:leading term of T_i^-1(x ast y)}, we obtain
		\[
		\tTT_i^{-1}(xy)=\widetilde{T}_{r_i}^{-1}(x)\widetilde{T}_{r_i}^{-1}(y)=\widetilde{T}_{r_i}^{-1}(xy).
		\]
		This completes the induction step and the theorem is proved.
	\end{proof}

	\subsection{Relative braid group action on $\ff$}
	
	We now consider the action of $\widetilde{\TT}_i$ on $\ff$. 
	
	\begin{lemma}\label{lem: T_bri apply on tensor of xu}
		For any $x\in{^\sigma\ff}[i,\tau i]$, $y\in\ff$, $u,w\in\ff_{i,\tau i}$, $\alpha_0\in\N^\I$, $\beta_0\in\Lambda_{i,\tau i}$, we have 
		\[
		\widetilde{\mathscr{T}}_{r_i}^{-1} \big((\tau (xu)\otimes yw)K'_{\alpha_0+\wt(x)}K'_{\beta_0+\wt(u)}\big)\in \tilde{k}_{\tau\tau_i\wt(u)}^{-1}\sum_{\gamma_1\in\Lambda_{i,\tau i},\gamma_2\in\Xi_i}\U^+\U^-\cdot K_{\gamma_1}K_{\gamma_2}'.\]
		Moreover, the left-hand side has its leading term in $\tilde{k}_{\tau\tau_i\wt(u)}^{-1}\U^+\U^-$ only if $\alpha_0=\beta_0=0$ and $x=1$.
	\end{lemma}
	\begin{proof}
		For $y=w=1$ the claim is clear. Now we prove by induction on $\wt(yw)$. Note that 
		\begin{align*}
			(\tau (xu)\otimes yw)&=(\tau (xu)\otimes 1)(1\otimes yw)\\
			&\quad -\sum_{y_{(2)}w_{(2)}\neq 1}\varphi(\tau(x_{(1)}u_{(1)}),y_{(2)}w_{(2)})(\tau (x_{(2)}u_{(2)})\otimes y_{(1)}w_{(1)}h_{\wt(w_{(2)}w_{(2)})})\\
			&\in(\tau (xu)\otimes 1)(1\otimes yw)-\sum_{y_{(2)}w_{(2)}\neq 1}\Q(v)(\tau (x_{(2)}u_{(2)})\otimes y_{(1)}w_{(1)})K_{\wt(y_{(2)}w_{(2)})}'.
		\end{align*}
		Since $x_{(2)}\in{^\sigma\ff}[i,\tau i]$ and $\wt(y_{(1)}w_{(1)})<\wt(yw)$, by the induction hypothesis we find that 
		\begin{align*}
			&(\tau (x_{(2)}u_{(2)})\otimes y_{(1)}w_{(1)})K_{\wt(y_{(2)}w_{(2)})}'K'_{\alpha_0+\wt(x)}K'_{\beta_0+\wt(u)}\\
			&=(\tau (x_{(2)}u_{(2)})\otimes y_{(1)}w_{(1)})K'_{\wt(y_{(2)})+\alpha_0+\wt(x)}K'_{\wt(w_{(2)})+\beta_0+\wt(u_{(1)})+\wt(u_{(2)})}\\
			&\in \tilde{k}_{\tau\tau_i\wt(u_{(2)})}^{-1}\sum_{\substack{\gamma_1\in\Lambda_{i,\tau i},\gamma_2\in\Xi_i\\(\gamma_1,\gamma_2)\neq (0,0)}}\U^+\U^-\cdot K_{\gamma_1}K_{\gamma_2}'\\
			&\subseteq \tilde{k}_{\tau\tau_i\wt(u)}^{-1}\sum_{\substack{\gamma_1\in\Lambda_{i,\tau i},\gamma_2\in\Xi_i\\(\gamma_1,\gamma_2)\neq (0,0)}}\U^+\U^-\cdot K_{\gamma_1}K_{\gamma_2}'.
		\end{align*}
		It therefore remains to consider $\widetilde{\mathscr{T}}_{r_i}^{-1}((\tau (xu)\otimes 1)(1\otimes yw)K'_{\alpha_0+\wt(x)}K'_{\beta_0+\wt(u)})$, for which we may assume that $y=w'y'$ with $w'\in\ff_{i,\tau i},y'\in\ff[i,\tau i]$. Then
		\begin{equation}\label{eq: T_ri on tensor leading}
			\begin{aligned}
				&\widetilde{\mathscr{T}}_{r_i}^{-1}((\tau (xu)\otimes 1)(1\otimes yw)K'_{\alpha_0+\wt(x)}K'_{\beta_0+\wt(u)})\\
				&\in\Q(v)\cdot\widetilde{\mathscr{T}}_{r_i}^{-1}(\tau x^+)\tau\tau_i(u^-)\tau_i(w'^+)\widetilde{\mathscr{T}}_{r_i}^{-1}(y'^-)\tau_i(w^+)\tilde{k}_{\tau\tau_i\wt(u)}^{-1}\cdot K'_{r_i(\alpha_0+\wt(x))-\tau_i(\beta_0+\wt(w)+\wt(w'))}\\
				&\subseteq \tilde{k}_{\tau\tau_i\wt(u)}^{-1}\sum_{\gamma_1,\gamma_2\in\Lambda_{i,\tau i}}\U^+\U^-\cdot K_{\gamma_1}K_{r_i(\alpha_0+\wt(x))-\tau_i(\beta_0)-\gamma_2}'
			\end{aligned}
		\end{equation}
		where we have used the fact
		\[\tau\tau_i(u^-)\tau_i(w'^+)\widetilde{\mathscr{T}}_{r_i}^{-1}(y'^-)\tau_i(w^+)\in\sum_{\gamma_1,\gamma_2\leq\tau_i(\wt(w)+\wt(w'))}\U^+\U^-K_{\gamma_1}K'_{\gamma_2}.\]
		Note that $r_i(\alpha_i+\wt(x))-\tau_i(\beta_0)-\gamma_2=0$ only if $x=1$ and $\alpha_0=\beta_0=\gamma_2=0$, so this completes the induction step.
	\end{proof}
	
	\begin{corollary}
		For any $x\in{^\sigma\ff}[i,\tau i]$ and $u\in\ff_{i,\tau i}$, we have 
		\begin{equation}
			\label{eq:T_ri on product conjugate-1}
			\begin{aligned}
				&\widetilde{\xi}_\tau\big(\tTT_i^{-1}(xu)\big)\widetilde{\Upsilon}_i =\widetilde{\Upsilon}_i\widetilde{\mathscr{T}}_{r_i}^{-1}\big(\widetilde{\xi}_\tau(xu)\big)\\
				&\qquad\in v^{\frac{1}{2}\Z}\tilde{k}_{\tau\tau_i\wt(u)}^{-1}\cdot\widetilde{\Upsilon}_i \tau\tau_i(u^-)\widetilde{\mathscr{T}}_{r_i}^{-1}(x^-)+\tilde{k}_{\tau\tau_i\wt(u)}^{-1}\sum_{\substack{\gamma_1\in\Lambda_{i,\tau i},\gamma_2\in\Xi_i\\(\gamma_1,\gamma_2)\neq (0,0)}}\U^+\U^-\cdot K_{\gamma_1}K_{\gamma_2}'.
			\end{aligned}
		\end{equation}
	\end{corollary}
	
	\begin{proof}
		Note that 
		\begin{equation*}
			\widetilde{\xi}_\tau(xu) =\sum\chi(x_{(2)}h_{\wt(x_{(3)})}u_{(2)}h_{u_{(3)}})\tau (x_{(3)} u_{(3)})\otimes x_{(1)}h_{\wt(x_{(2)})+\wt(x_{(3)})}u_{(1)}h_{\wt(u_{(2)})+\wt(u_{(3)})}.
		\end{equation*}
		
		Combining this expression with Lemma~\ref{lem: T_bri apply on tensor of xu} gives us  
		\begin{equation*}
			\widetilde{\mathscr{T}}_{r_i}^{-1}\big(\widetilde{\xi}_\tau(xu)\big)\in v^{\frac{1}{2}\Z}\tilde{k}_{\tau\tau_i\wt(u)}^{-1}\cdot \tau\tau_i(u^-)\widetilde{\mathscr{T}}_{r_i}^{-1}(x^-)+\tilde{k}_{\tau\tau_i(u)}^{-1}\sum_{\substack{\gamma_1\in\Lambda_{i,\tau i},\gamma_2\in\Xi_i\\(\gamma_1,\gamma_2)\neq (0,0)}}\U^+\U^-\cdot K_{\gamma_1}K_{\gamma_2}'
		\end{equation*}
		where the leading term comes from $x_{(2)}=x_{(3)}=1$, $u_{(1)}=u_{(2)}=1$, and we also use \eqref{eq: T_ri on tensor leading}. The claim \eqref{eq:T_ri on product conjugate-1} now follows from this by Theorem~\ref{thm: relative T_i conjugate}.
	\end{proof}
	
	\begin{lemma}\label{lem: T_i on product cartan}
		For $x\in{^\sigma\ff}[i,\tau i]$ and $u\in\ff_{i,\tau i}$, we have 
		\[\widetilde{\TT}_i^{-1}(xu)\in \K_{\tau\tau_i\wt(u)}^{-1}\ast\ff.\]
	\end{lemma}
	
	\begin{proof}
		Recall from Proposition \ref{prop:TiBi} that
		\begin{align} \label{Tifi}
			\widetilde{\TT}_i^{-1}(u)=\K_{\tau\tau_i\wt(u)}^{-1}\diamond \tau\tau_i(u),
			\qquad \text{ for } u\in\ff_{i,\tau i}.
		\end{align}
		
		Let $x'=\widetilde{T}_{r_i}^{-1}(x)$ and $u'=\tau\tau_i(u)$. Clearly we have $x'u' \in \ff$. Applying Theorem~\ref{thm:Tifi} and \eqref{Tifi}, we have 
		\[
		\widetilde{\TT}_i^{-1}(x\ast u)=\widetilde{T}_i^{-1}(x)\ast \K_{\tau\tau_i\wt(u)}^{-1}\diamond u'\in v^{\frac{1}{2}\Z}\K_{\tau\tau_i\wt(u)}^{-1}\ast(x'u').\]
		On the other hand, we have 
		\begin{align*}
			\widetilde{\TT}_i^{-1}(x\ast u)&=\widetilde{\TT}_i^{-1}\sum\varphi(x_{(1)},u_{(2)})x_{(2)}u_{(1)}h_{\wt(u_{(2)})}\\
			&=\widetilde{\TT}_i^{-1}(xu)+\sum\varphi(x_{(1)},u_{(2)})\widetilde{\TT}_i^{-1}(x_{(2)}u_{(1)}h_{\wt(u_{(2)})})\\
			&\in\widetilde{\TT}_i^{-1}(xu)+\sum\Q(v)\K_{\tau\tau_i\wt(u_{(1)})}^{-1}\K_{\tau\tau_i\wt(u_{(2)})}^{-1}\ast(\tau\tau_i(u_{(1)})\widetilde{\TT}_i^{-1}(x_{(2)}))\\
			&\subseteq \widetilde{\TT}_i^{-1}(xu)+\K_{\tau\tau_i\wt(u)}^{-1}\ast\ff.
		\end{align*}
		The lemma now follows by comparing the two statements above.
	\end{proof}
	
	We can now prove the second main result of this section, generalizing the formula in Theorem \ref{thm:Tifi}. This result seems new even in the context of $\tU$ viewed as iquantum group of diagonal type. 
	
	\begin{theorem}
		\label{thm:Tifif_i}
		\begin{enumerate}
			\item For any $x\in{^\sigma\ff}[i,\tau i]$ and $u\in\ff_{i,\tau i}$, we have 
			\begin{align*}
				\widetilde{\TT}_i^{-1}(xu)=v^{-\frac{1}{2}(\tau\wt(u)+\wt(u),\wt(x))}\K_{\tau\tau_i\wt(u)}^{-1}\diamond \big(\tau\tau_i(u)\widetilde{T}_{r_i}^{-1}(x)\big).
			\end{align*}
			\item 
			For any $x\in\ff[i,\tau i]$ and $u\in\ff_{i,\tau i}$, we have 
			\begin{align*}
				\widetilde{\TT}_i(ux)=v^{-\frac{1}{2}(\tau\wt(u)+\wt(u),\wt(x))}\K_{\tau_i\wt(u)}^{-1}\diamond \big(\widetilde{T}_{r_i}(x)\tau\tau_i(u)\big).
			\end{align*}
		\end{enumerate}
	\end{theorem}
	
	\begin{proof}
		Since $\widetilde{T}_{r_i}$ induces an isomorphism from $\ff[i,\tau i]$ to ${^\sigma\ff}[i,\tau i]$, the second statement follows immediately from the first one. 
		
		Now we prove (1). First, by applying Lemma~\ref{lem: T_i on product cartan} and Lemma~\ref{lem: xi_theta leading term}, we have that 
		\begin{align*}
			\widetilde{\xi}_\tau(\tTT_i^{-1}(xu))\widetilde{\Upsilon}_i&\in \tilde{k}_{\tau\tau_i\wt(u)}^{-1}\Big(\U^-+\sum_{\alpha\in\N^\I\setminus\{0\}}\U^+\U^-\cdot K_\alpha'\Big)\widetilde{\Upsilon}_i\\
			&\subseteq\tilde{k}_{\tau\tau_i\wt(u)}^{-1}\Big(\U^+\U^-\hU_{i,\tau i}^0+\sum_{\alpha\in\N^\I\setminus\{0\}}\U^+\U^-\hU_{i,\tau i}^0\cdot K_\alpha'\Big).
		\end{align*}
		Comparing this with \eqref{eq:T_ri on product conjugate-1} gives us
		\begin{equation}\label{eq:T_ri on product leading term-2}
			\widetilde{\mathscr{T}}_{r_i}^{-1}(\widetilde{\xi}_\tau(xu))
			\in v^{\frac{1}{2}\Z}\tilde{k}_{\tau\tau_i\wt(u)}^{-1}\cdot \tau\tau_i(u^-)\widetilde{\mathscr{T}}_{r_i}^{-1}(x^-)+\tilde{k}_{\tau\tau_i\wt(u)}^{-1}\sum_{\substack{\gamma_1\in\Lambda_{i,\tau i},\gamma_2\in\Xi_i\cap\N^\I\\(\gamma_1,\gamma_2)\neq (0,0)}}\U^+\U^-\cdot K_{\gamma_1}K_{\gamma_2}'.
		\end{equation}
		Now conjugating $\widetilde{\mathscr{T}}_{r_i}^{-1}(\widetilde{\xi}_\tau(xu))$ in \eqref{eq:T_ri on product leading term-2} by $\widetilde{\Upsilon}_i$, we finally get
		
		\begin{equation}\label{eq:T_i on product leading term}
			\begin{aligned}
				&\widetilde{\xi}_\tau(\tTT_i^{-1}(xu))=\widetilde{\Upsilon}_i\widetilde{\mathscr{T}}_{r_i}^{-1}(\widetilde{\xi}_\tau(xu))\widetilde{\Upsilon}_i^{-1}\\
				&\in v^{\frac{1}{2}n(x,u)}\tilde{k}_{\tau\tau_i\wt(u)}^{-1}\cdot\widetilde{\Upsilon}_i \tau\tau_i(u^-)\widetilde{\mathscr{T}}_{r_i}^{-1}(x^-)\widetilde{\Upsilon}_i^{-1}+\tilde{k}_{\tau\tau_i\wt(u)}^{-1}\sum_{\substack{\gamma_1\in\Lambda_{i,\tau i},\gamma_2\in\Xi_i\cap\N^\I\\(\gamma_1,\gamma_2)\neq (0,0)}}\U^+\U^-\cdot K_{\gamma_1}K_{\gamma_2}'
			\end{aligned}
		\end{equation}
		for some $n(x,u)\in\Z$. 
		
		To prove (1), from Lemma~\ref{lem: T_i on product cartan} we see that $y:=\K_{\tau\tau_i\wt(u)}\ast(\tTT_i^{-1}(xu))\in\ff$, so \eqref{eq:T_i on product leading term} and Lemma~\ref{lem: xi_theta leading term} gives
		\[y^-\in v^{\frac{1}{2}n(x,u)}\widetilde{\Upsilon}_i \tau\tau_i(u^-)\widetilde{\mathscr{T}}_{r_i}^{-1}(x^-)\widetilde{\Upsilon}_i^{-1}+\sum_{\gamma\in\N^\I\setminus\{0\}}\U^+\U^-\cdot K_{\gamma}'.\]
		Since $\widetilde{\TT}_i^{-1}(u)=\K_{\tau\tau_i\wt(u)}^{-1}\diamond \tau\tau_i u$, by similar argument as in the proof of Theorem~\ref{thm:Tifi} we have 
		\[
		y=v^{\frac{1}{2}n(x,u)} \tau\tau_i(u)\widetilde{T}_{r_i}^{-1}(x).
		\]
		To determine the integer $n(x,u)$ we must compute the leading term of $\widetilde{\xi}_\tau(\tTT_i^{-1}(xu))$. In view of the proof of Lemma~\ref{lem: T_bri apply on tensor of xu}, this is given by the leading term of $\widetilde{\mathscr{T}}_{r_i}^{-1}(\chi(h_{\wt(u)})\tau u\otimes xh_{\wt(u)})$, which is equal to
		\begin{align*}
			&\chi(h_{\wt(u)})\widetilde{\mathscr{T}}_{r_i}^{-1}(v^{-(\wt(u),\tau\wt(u))}\tau (u^+)x^-K_u')\\
			&=v^{-\frac{1}{2}(\wt(u),\wt(u))-\frac{1}{2}(\wt(u),\tau\wt(u))}K_{\tau\tau_i\wt(u)}^{-1}\tau\tau_i(u^-)\cdot\widetilde{T}_{r_i}^{-1}(x^-)\cdot v^{-\frac{1}{2}(\wt(u),\tau\wt(u))}K_{\tau_i\wt(u)}'^{-1}\\
			&=v^{\frac{1}{2}(\tau \wt(u)-\wt(u),\wt(u))-(\wt(u),\wt(x))}\widetilde{\xi}_\tau(\K_{\tau\tau_i\wt(u)}^{-1})(\tau\tau_i(u^-)\cdot\widetilde{T}_{r_i}^{-1}(x^-))\\
			&=v^{-\frac{1}{2}(\tau\wt(u)+\wt(u),\wt(x))}\widetilde{\xi}_\tau(\K_{\tau\tau_i\wt(u)}^{-1})\diamond(\tau\tau_i(u^-)\cdot\widetilde{T}_{r_i}^{-1}(x^-)).
		\end{align*}
		This completes the proof.
	\end{proof}
	
	The following corollary is a special case of \cite[Theorem 7.13]{WZ23} for quasi-split iquantum groups, whose original proof is completely different and difficult.
	
	\begin{corollary}
		\label{cor:braid-simple}
		Suppose that $wi\in\I$, for $w \in W_\tau$ and $i \in \I$. Then  $\tTT_{w}(B_i) = B_{wi}$ in $\widetilde{\mathbf{B}}^\imath_\tau\equiv\tUi$.
	\end{corollary}
	
	\begin{proof}
		We assume that $w=r_{i_1}r_{i_2}\cdots r_{i_t}$ is a reduced expression of $w\in W_\tau$. Denote by $w_k:=r_{i_k}r_{i_{k+1}}\cdots r_{i_t}$ for $1\leq k\leq t$. Then we have 
		$\widetilde{T}_{w_k}(E_i)\in\U^+$ by using the assumption $wi\in\I$; see \cite[Proposition~ 8.20]{Jan96}. Therefore, we have $\widetilde{T}_{w_{k+1}}(\vartheta_i)\in\ff[i_{k},\tau i_{k}]$ for $1\leq k< t$, and it follows from Theorem~ \ref{thm:Tifif_i}(2) that 
		$\widetilde{T}_{w_{k}}(\vartheta_i)=\tTT_{w_k}(\vartheta_i)$ for $1\leq k\leq t$  by induction. 
		In particular $\tTT_{{w}}(\vartheta_i)=\widetilde{T}_{{w}}(\vartheta_i)$, which equals to $\vartheta_{wi}$ by using \cite[Proposition 8.20]{Jan96} again.
	\end{proof}
	

	\section{Dual canonical bases for iquantum groups}
	\label{sec:DCBiQG}
	
	In this section, we restrict ourselves to quantum groups and quasi-split iquantum groups $\tUi$ of arbitrary finite type. We shall construct the dual canonical basis on $\tUi$. 
	
	\subsection{Dual canonical basis on $\ff$}
	\label{subsec:dCB:f}
	
	A canonical basis $\mathbf{C}^{\mathrm{can}}$ of $\ff$ was constructed by Lusztig \cite{Lus90a} in ADE type and it is now available for all (finite) types \cite{Ka91, Lus93}. We denote by $\{b^*\mid b\in\mathbf{C}^{\mathrm{can}}\}$ the dual basis of $\mathbf{C}^{\mathrm{can}}$ with respect to the bilinear form $\varphi(\cdot,\cdot)$ on $\ff$ defined in \eqref{eq:hopf-pairing}, i.e. $\varphi(b^*, b') =\delta_{b,b'}$, for $b,b'\in\mathbf{C}^{\mathrm{can}}$. Recalling the bilinear form $(\cdot,\cdot)$ from \eqref{BilForm}, we define a norm function $\texttt{N}:\Z^\I\rightarrow\Z$ by 
	\[
	\texttt{N}(\alpha)=\frac{1}{2}(\alpha,\alpha)-\text{ht}(\alpha),
	\]
	where the height function $\text{ht}:\Z^\I\rightarrow\Z$ is given by $\text{ht}(\sum_ia_i\alpha_i)=\sum_ia_i$. The rescaled dual canonical basis of $\ff$ is then defined to be $\delta_b:=v^{\frac{1}{2}\texttt{N}(\wt(b))}b^*$ and
	\[
	\mathbf{C}:=\{v^{\frac{1}{2}\texttt{N}(\wt(b))}b^*\mid b\in\mathbf{C}^{\mathrm{can}}\}.
	\]
	
	\begin{example}
		For $\tU=\tU_v(\mathfrak{sl}_2)$, we have $\mathbf{C}=\{\vartheta_1^n\mid n\in\N\}$. 
	\end{example}
	Set $\mathcal{Z}=\Z[v^{\frac12},v^{-\frac12}]$. We define the integral form $\ff_{\mathcal{Z}}$ to be the free $\mathcal{Z}$-submodule of $\ff$ generated by $\mathbf{C}$. It is known \cite[Theorem 14.4.13]{Lus93} that $\ff_{\mathcal{Z}}$ is an algebra over $\mathcal{Z}$; and for any $b',b'',c\in\mathbf{C}$ we have
	\begin{align}
		\label{eqn:coprod-dCB}
		b'b''=\sum_{b\in\mathbf{C}}g_{b',b''}^bb,
		\qquad\Delta(c)=\sum_{c',c''\in\mathbf{C}} f_{c',c''}^c c'h_{\wt(c'')}\otimes c'',
	\end{align}
	where $g_{b',b''}^b,f_{c',c''}^c\in\mathcal{Z}$.
	
	By a slight abuse of notation, we also denote by $\mathbf{C}^{\mathrm{can}}$ the canonical basis for $\U^+$ via the isomorphism $\ff \rightarrow \U^+$ from \eqref{ffUU}. Because our $E_i$ are dual generators (see \cite[Remark 3.1]{CLPRW}), we have $(v_i -v_i^{-1})^{-1} E_i \in \mathbf{C}^{\mathrm{can}}$ and $E_i \in \mathbf{C}$.
	
	\subsection{Integral form on $\hUi$}
	\label{subsec:integralform}
	
	Let $C$ be a Cartan matrix of finite type. Satake diagrams of finite type with $\tau\neq \Id$ can  be found in \cite[Table 3.1]{CLPRW}. Let $\tau$ be an involution in $\mathrm{Inv}(C)$ and consider the iHopf algebra $\hB^\imath_\tau$. Recall from \eqref{embed:f} the canonical embedding $\iota:\ff \to \hB^\imath_\tau$. Usually, we view $x\in\hB^\imath_\tau$ for any $x\in\ff$ by omitting $\iota$ if there is no confusion.

	We denote by $_{\mathcal{Z}}\hB^\imath_\tau$ the free $\mathcal{Z}$-submodule of $\hB^\imath_\tau$ generated by $\K_\alpha\diamond \iota(b)$, where $\alpha\in\N^\I$ and $b\in\mathbf{C}$. Similarly, let $_{\mathcal{Z}}\tB^\imath_\tau$ be the free $\mathcal{Z}$-submodule of $\tB^\imath_\tau$ generated by $\K_\alpha\diamond \iota(b)$, where $\alpha\in\Z^\I$ and $b\in\mathbf{C}$. 
	By definition, we have $\iota(\ff_\mathcal{Z})\subseteq {}_\mathcal{Z} \tB^\imath_\tau$.
	
	\begin{lemma}\label{lem: iHopf graded}
		The submodule $_{\mathcal{Z}}\hB^\imath_\tau$ is an $\N^\I$-graded $\mathcal{Z}$-algebra with weights given by \[
		\wt(h_\alpha)=\alpha+\tau\alpha, \quad
		\wt\big(\iota(x)\big)=\wt(x), \quad \text{ for } \alpha\in\N^\I, \, x\in\ff.
		\]
	\end{lemma}
	
	\begin{proof}
		By \cite[Theorem 3.11]{BG17a} we have $\varphi(b,b')\in\mathcal{Z}$ for any $b,b'\in\mathbf{C}$. Since $\mathbf{C}$ is an integral basis for the Hopf algebra $\hB$, it follows from the definition of $\ast$ that $_{\mathcal{Z}}\hB^\imath_\tau$ is an algebra over $\mathcal{Z}$. 
		
		To prove that $_{\mathcal{Z}}\hB^\imath_\tau$ is $\N^\I$-graded, we note that for $x,y\in\ff$,
		\begin{equation}\label{eq: monimial ast}
			x\ast y=\sum \varphi(x_{(1)}h_{\wt(x_{(2)})},y_{(2)})x_{(2)}y_{(1)}h_{\wt(y_{(2)})}.
		\end{equation}
		Since $\varphi(x,y)\neq 0$ only if $\wt(x)=\wt(y)$, any nonzero term on the right-hand side of \eqref{eq: monimial ast} satisfies $\wt(x_{(1)})=\wt(y_{(2)})$, and in this case
		\begin{align*}
			\wt(x_{(2)}y_{(1)}h_{\wt(y_{(2)})})&=\wt(x_{(2)})+\wt(y_{(1)})+\wt(y_{(2)})+\tau(\wt(y_{(2)}))\\
			&=\wt(x_{(2)})+\wt(y_{(1)})+\wt(y_{(2)})+\wt(x_{(1)})\\
			&=\wt(x)+\wt(y).
		\end{align*}
		This proves our assertion.
	\end{proof}
	
	Via the isomorphisms in Theorem \ref{thm:iH=iQG}, we define the integral forms  $\hUi_\mathcal{Z}$ (resp. $\tUi_\mathcal{Z}$) of $\hUi$ (resp. $\tUi$) such that the following $\mathcal{Z}$-algebra isomorphisms hold:
	\begin{align}
		\label{eq:iso-integral}
		\widehat{\Phi}^\imath: {}_{\mathcal{Z}}\hB^\imath_\tau \stackrel{\cong}{\longrightarrow} \hUi_\mathcal{Z}, \qquad \widetilde{\Phi}^\imath: {}_{\mathcal{Z}}\tB^\imath_\tau \stackrel{\cong}{\longrightarrow} \tUi_\mathcal{Z}.
	\end{align}
	
	\begin{remark}
		For Cartan matrix $C$ not of finite type, the $\mathcal{Z}$-module $_{\mathcal{Z}}\hB^\imath_\tau$ may not be an algebra, since $\varphi(b,b')$ may not be in $\mathcal{Z}$ for $b,b'\in\mathbf{C}$; see  \cite[Proposition 3.9]{BG17a}.
	\end{remark}
	
	Let $w_0$ be the longest element of the Weyl group $W$ with a reduced expression $w_0=s_{i_1}s_{i_2}\cdots s_{i_N}$. Set $\bi=(i_1,\dots,i_N)$, and define 
	\begin{align}  \label{roots:beta}
		\beta_{\bi,k} =s_{i_1}\cdots s_{i_{k-1}}(\alpha_{i_k}),\quad \forall 1\leq k\leq N.
	\end{align}
	Then $\{\beta_{\bi,1},\dots,\beta_{\bi,N}\}$ is the set of positive roots. Following \cite[Proposition 40.1.3]{Lus93}, we define $\vartheta_{\bi,k} \in \ff$ such that 
	\begin{align}  \label{Eik}
		\vartheta_{\bi,k}^+=E_{\bi,k},
		\qquad
		\text{ where }
		E_{\bi,k}=\widetilde{T}_{i_1}^{-1}\cdots\widetilde{T}^{-1}_{i_{k-1}}(E_{i_k}),\quad\forall 1\leq k\leq N.
	\end{align}
	Here $\widetilde{T}_i$ are the braid group symmetries on $\tU$ from Proposition \ref{prop:BG2U}. 
	
	For any $\ba=(a_1,\dots,a_N)\in\N^N$, we set
	\begin{align} \label{theta:ia}
		\vartheta_{\bi,\ba}=v^{\frac{1}{2}n_{\bi,\mathbf{a}}}\prod_{k=1}^{N}\vartheta_{\bi,k}^{a_k}
	\end{align}
	where 
	\[n_{\bi,\ba}=\sum_{1\leq k<l\leq N}(\beta_{\bi,k},\beta_{\bi,l})a_ka_l.\]
	Then $\{\vartheta_{\bi,\ba}\mid \ba\in\N^N\}$ forms a basis of $\ff_\mathcal{Z}$, called the dual PBW basis; see \cite{Lus93, BG17b}.
	A direct construction of dual canonical basis $\mathbf{C}$ was given in \cite[Theorem 1.1]{BG17b} via the dual PBW basis: for each $\ba=(a_1,\dots,a_N)\in\N^N$, there is a unique element $b_{\bi,\ba}\in\mathbf{C}$ such that 
	\begin{equation} 
		\label{eq: QG dPBW and dCB}
		b_{\bi,\ba} \in \vartheta_{\bi,\ba} +\sum_{\ba'\prec\ba}v^{-1}\Z[v^{-1}]\vartheta_{\bi,\ba'}
	\end{equation}
	where $\preceq$ is the partial order on $\N^N$ defined in \cite[Section 4.2]{BG17b}. In spite of the notation $b_{\bi,\ba}$, we emphasize that the dual canonical basis $\mathbf{C}$ is independent of the choice of the reduced expression $\bi$. 
	
	\begin{corollary}
		$\{\K_\alpha\diamond \vartheta_{\bi,\ba}\mid \alpha\in\N^\I,\ba\in\N^N\}$ forms a $\mathcal Z$-basis of $_{\mathcal{Z}}\hB^\imath_\tau$.
	\end{corollary}
	This basis will be called the (dual) PBW basis for $_{\mathcal{Z}}\hB^\imath_\tau$ and $\hB^\imath_\tau$. 
	
	\begin{proof}
		Follows by definition of the integral $\mathcal Z$-form $_{\mathcal{Z}}\hB^\imath_\tau$ and \eqref{eq: QG dPBW and dCB}. 
	\end{proof}
	
	\begin{proposition}\label{lem:bar inv on iHopf integral}
		The bar-involution on $\hB^\imath_\tau$ (and also $\tB^\imath_\tau$) preserves the integral forms $_{\mathcal{Z}}\hB^\imath_\tau$ (and also $_{\mathcal{Z}}\tB^\imath_\tau$).
	\end{proposition}

	\begin{proof}
		We focus on $_{\mathcal{Z}}\hB^\imath_\tau$. Since the $\mathcal Z$-algebra $_{\mathcal{Z}}\hB^\imath_\tau$ is generated by $\{\K_\alpha,\iota(\vartheta_{\bi,k})\mid 1\leq k\leq N,\alpha\in\N^\I\}$, it is enough to prove that $\ov{\iota(\vartheta_{i,k})}\in{_{\mathcal{Z}}\hB}^\imath_\tau$ for $1\leq k\leq N$.
		
		It is known that the longest element $w_0\in W$ belongs to $W_\tau$, and let
		$w_0=r_{j_1}r_{j_2}\cdots r_{j_m}$ be a reduced expression in $W_\tau$; see \eqref{def:ri} for the definition of $r_i$. Let
		$\bi=(i_1,\dots,i_N)$ be a (fixed) sequence constructed from $(j_1,\dots,j_m)$ by replacing $j_k$ by $(t_1,\dots,t_p)$ if we have reduced expressions $r_{j_k}=s_{t_1}s_{t_2}\cdots s_{t_p}$ for $1\leq k\leq m$. 
		
		Corresponding to \eqref{def:ri}, we define
		\begin{align}
			\Phi^+(r_i)=\begin{cases}
				\{\alpha_i\},& \text{ if }c_{i,\tau i}=2,
				\\
				\{\alpha_i,\alpha_{\tau i}\}, &\text{ if }c_{i,\tau i}=0,
				\\
				\{\alpha_i,\alpha_{\tau i},\alpha_i+\alpha_{\tau i}\},&\text{ if }c_{i,\tau i}=-1.
			\end{cases}
		\end{align}
		Then $\Phi^+=\{r_{j_1}\cdots r_{j_{k-1}}(\beta)\mid 1\leq k\leq m, \beta\in\Phi^+(r_{j_{k}})\}$. 
		
		For $1\leq k\leq N$, if $s_{i_1}\cdots s_{i_{k-1}}(\alpha_{i_k})=r_{j_1}\cdots r_{j_{t-1}}(\alpha_{l})$ for some $1\leq t\leq m$ and $\alpha_l\in \Phi^+(r_{j_{t}})$, then we have 
		$\iota(\vartheta_{\bi,k})=\tTT_{j_1}^{-1}\cdots \tTT_{j_{t-1}}^{-1}(B_l)$ by Theorem \ref{thm:Tifi}. In this case, we have $\ov{\iota(\vartheta_{\bi,k})}=\iota(\vartheta_{\bi,k})$ by Lemma \ref{lem:bar-invar-braid-Ui}. Otherwise, $s_{i_1}\cdots s_{i_{k-1}}(\alpha_{i_k})=r_{j_1}\cdots r_{j_{t-1}}(\alpha_{j_t}+\alpha_{\tau j_t})$ for some $1\leq t\leq m$ (in this case, $c_{j_t,\tau j_t}=-1$). By our assumption,  $s_{i_{k-1}}(\alpha_{i_k})=\alpha_{j_t}+\alpha_{\tau j_t}$, and we further assume that $i_k=j_t$, and then $i_{k-1}=\tau j_{t}$. In this way, we have
		$\iota(\vartheta_{\bi,k})=\tTT_{j_1}^{-1}\cdots \tTT_{j_{t-1}}^{-1}(\iota (\widetilde{T}_{j_t}^{-1}(\vartheta_{\tau j_t})))$. By definition, 
		\begin{align*}
			\widetilde{T}_{j_t}^{-1}(\vartheta_{\tau j_t}))&=(v_{j_t}-v_{j_t}^{-1})^{-1}\big(v_{j_t}^{1/2}\vartheta_{\tau j_t}\vartheta_{j_t}-v_{j_t}^{-1/2}\vartheta_{j_t}\vartheta_{\tau j_t}\big)
			\\
			&=(v_{j_t}-v_{j_t}^{-1})^{-1}\big(v_{j_t}^{1/2}\vartheta_{\tau j_t}*\vartheta_{j_t}-v_{j_t}^{-1/2}\vartheta_{j_t}*\vartheta_{\tau j_t}\big)+\K_{j_t}-v_{j_t}\K_{\tau j_t},
		\end{align*}
		so we have
		\begin{align*}
			\ov{\iota(\widetilde{T}_{j_t}^{-1}(\vartheta_{\tau j_t})))}=\iota(\widetilde{T}_{j_t}^{-1}(\vartheta_{\tau j_t})))+(v_{j_t}-v_{j_t}^{-1})\K_{\tau j_t}.
		\end{align*}
		Therefore,
		\begin{align*}
			\ov{\iota(\vartheta_{\bi,k})}=&\iota(\vartheta_{\bi,k})+(v_{j_t}-v_{j_t}^{-1})\tTT_{j_1}^{-1}\cdots \tTT_{j_{t-1}}^{-1}(\K_{\tau j_t})\in {}_{\mathcal{Z}}\hB^\imath_\tau.
		\end{align*}
		The proof is completed.
	\end{proof}

	\subsection{Dual canonical basis on $\hUi$}
	\label{subsec: dCB from CB}
	
	We now define a partial order $\preceq$ on the set $\N^\I\times \mathbf{C}$ as follows: $(\alpha,b)\preceq(\beta,b')$ if 
	\begin{itemize}
		\item[(1)] $\alpha+\tau\alpha+\wt(\iota(b))=\beta+\tau\beta+\wt(\iota(b'))$, and
		\item[(2)] $0\neq \beta-\alpha\in\N^\I$ or $(\alpha,b)=(\beta,b')$.
	\end{itemize}
	We denote $(\alpha,b)\prec(\beta,b')$ if $(\alpha,b)\preceq(\beta,b')$ and $(\alpha,b)\neq (\beta,b')$.
	
	\begin{lemma}
		For any $\alpha\in\N^\I$ and $b\in \mathbf{C}$ we have
		\begin{equation}
			\label{eq:bar_dCB}
			\ov{\K_\alpha\diamond \iota(b)} \in \K_\alpha\diamond\iota(b) +\sum_{(\alpha,b)\prec(\beta,b')}\Z[v,v^{-1}]\cdot \K_\beta\diamond \iota(b').
		\end{equation}
	\end{lemma}
	
	\begin{proof}
		Let us consider the quotient map $\pi:\hB^\imath_\tau\to\ff$ defined by the ideal generated by $h_i$, $i\in\I$. From \eqref{star product} we see that $\pi\circ\iota=\Id$ and $\pi$ is a homomorphism of algebras. Moreover, if we define an anti-involution of $\ff$ by $\ov{\vartheta_i}=\vartheta_i$, $\ov{v}=v^{-1}$, then we have $\pi(\ov{x})=\ov{x}$ for any $x\in\hB^\imath_\tau$.
		
		Now it suffices to consider the lemma for $\alpha=0$, in which case we have for $b\in\mathbf{C}^{\mathrm{can}}$,
		\[
		\pi(\ov{\iota(b^*)})-v^{\texttt{N}(\wt(b))}\pi(\iota(b^*))=\ov{\pi(\iota(b^*))}-v^{\texttt{N}(\wt(b))}\pi(\iota(b^*))=\ov{b^*}-v^{\texttt{N}(\wt(b))}b^*=0.
		\]
		Since $\iota$ is weight-preserving, this implies by Proposition~\ref{lem:bar inv on iHopf integral} that 
		\[
		\ov{\iota(b^*)}-v^{\texttt{N}(\wt(b))}\iota(b^*)\in\sum \Z[v,v^{-1}]\cdot h_{\tau\beta}\ast\iota(b'^*).
		\]
		From the definition of $\delta_b$, we then get
		\[
		\ov{\iota(\delta_b)}-\iota(\delta_b)\in\sum \Z[v,v^{-1}]\cdot v^{d(\beta,b')} \K_{\beta}\diamond\iota(\delta_{b'})
		\]
		where
		\[
		d(\beta,b')=-\frac{1}{2}\texttt{N}(\wt(b))-\frac{1}{2}\texttt{N}(\wt(b'))-\frac{1}{2}(\beta,\tau\beta)-\frac{1}{2}(\beta-\tau\beta,\wt(b')).
		\]
		Write $\eta:=\wt(b)\in\N^\I$, then we have $\wt(b')=\eta-\beta-\tau\beta$, so
		\begin{align*}
			d(\beta,b')&=-\frac{1}{4}(\eta,\eta)+\frac{1}{2}\text{ht}(\eta)-\frac{1}{4}(\eta-\beta-\tau\beta,\eta-\beta-\tau\beta)+\frac{1}{2}\text{ht}(\eta-\beta-\tau\beta)\\
			&-\frac{1}{2}(\beta,\tau\beta)-\frac{1}{2}(\beta-\tau\beta,\eta-\beta-\tau\beta)\\
			&=-\frac{1}{2}(\beta,\beta)-\frac{1}{2}(\eta,\eta)+(\eta-\beta,\tau\beta)+\text{ht}(\eta-\beta)\in\Z.
		\end{align*}
		This completes the proof.
	\end{proof}
	
	For a given pair $(\alpha,b)$, there are only finitely many pairs $(\beta,b')$ such that $(\alpha,b)\prec(\beta,b')$. Hence Lusztig's Lemma is applicable thanks to \eqref{eq:bar_dCB} and we obtain a bar-invariant basis for $\hB^\imath_\tau$. We also note the $\diamond$-action preserves our basis, more precisely, we have established the following. 
	
	\begin{theorem}
		\label{thm:dCB}
		For each $\alpha\in\N^\I$ and $b\in \mathbf{C}$, there exists a unique element $C_{\alpha,b}\in \hB^\imath_\tau$ such that $\ov{C_{\alpha,b}}=C_{\alpha,b}$ and
		\begin{align}
			\label{eq:dualCB}
			C_{\alpha,b}\in \K_\alpha\diamond \iota(b) +\sum_{(\beta,b')}v^{-1}\Z[v^{-1}]\cdot \K_\beta\diamond \iota(b').
		\end{align}
		Then $\{C_{\alpha,b} \mid \alpha\in\N^\I,b\in \mathbf{C}\}$ forms a basis for $\hB^\imath_\tau$. 
		Moreover, $C_{\alpha,b}$ satisfies 
		\[
		C_{\alpha,b}\in\K_\alpha\diamond \iota(b) +\sum_{(\alpha,b)\prec(\beta,b')}v^{-1}\Z[v^{-1}]\cdot \K_\beta\diamond \iota(b'),
		\]
		and $C_{\alpha,b}=\K_\alpha\diamond C_{0,b}$.
	\end{theorem}
	
	Writing $C_b:=C_{0,b}$, we use $\{\K_\alpha\diamond C_b\mid \alpha\in\N^\I,b\in \mathbf{C}\}$ to denote the basis in Theorem~\ref{thm:dCB}; this is called the dual canonical basis of $\hB^\imath_\tau$. The dual canonical basis of $\tB^\imath_\tau$ is defined to be $\{\K_\alpha\diamond C_b\mid \alpha\in\Z^\I,b\in\mathbf{C}\}$. 
	
	From the integral properties of Lusztig's canonical basis, we have the following.
	
	\begin{corollary}
		$\{\K_\alpha\diamond C_b\mid \alpha\in\N^\I,b\in\mathbf{C}\}$ forms a basis of the $\mathcal{Z}$-algebra ${}_\mathcal{Z}\hB_\tau^\imath$, and $\{\K_\alpha\diamond C_b\mid \alpha\in\Z^\I,b\in\mathbf{C}\}$ forms a basis of the $\mathcal{Z}$-algebra ${}_\mathcal{Z}\tB_\tau^\imath$.
	\end{corollary}
	
	These dual canonical bases can be transfered to the ones of the integral forms $\hUi_{\mathcal{Z}}$ and $\tUi_{\mathcal{Z}}$ of iquantum groups via the isomorphisms in \eqref{eq:iso-integral}.

	
	Due to Lemma \ref{lem:involution-iQG} and the isomorphisms given in Theorem~\ref{thm:iH=iQG}, there exists an anti-involution $\sigma^\imath$ on $\hB^\imath_\tau$ (also $\tB^\imath_\tau$) given by  $\sigma^\imath(\vartheta_i)=\vartheta_i$, $i\in\I$ and $\sigma^\imath(h_\alpha)=h_{\tau(\alpha)}$, $\alpha\in\N^\I$.  We use $\sigma$ to denote the anti-involution on $\ff$ by sending $\vartheta_i\mapsto \vartheta_i$. Recall the natural inclusion $\iota:\ff\rightarrow \hB^\imath_\tau$.
	\begin{lemma}
		\label{lem:com-anti-involution}
		For any $x\in\ff$, we have 
		$\sigma^\imath(\iota(x))=\iota(\sigma(x))$.
	\end{lemma}
	
	\begin{proof}
		We prove by induction on the weight of $x$. For $x=1$ or $x=\vartheta_i$ this is obvious.
		Now let $\mu \in \N^\I \backslash\{0\}$ and assume the result holds for $x$ with $\wt(x)<\mu$. It suffices to prove the statement for any element in $\ff_\mu$ of the form $\vartheta_i\cdot x$.
		
		It follows from Lemma \ref{lem:recursive} that
		\begin{align*}
			\iota(\sigma(\vartheta_i\cdot x))&=\iota(\sigma(x))\cdot\iota(\vartheta_i)
			\\
			&=\sigma(x)*\vartheta_i-(v_i-v_i^{-1})\partial_{\tau i}^R(\sigma(x))\cdot h_i
			\\
			&=\sigma(x)*\vartheta_i-(v_i-v_i^{-1})\sigma(\partial_{\tau i}^L(x))\cdot h_i,
		\end{align*}
		where the last equality follows from $\sigma\circ \partial_{\tau i}^R=\partial_{\tau i}^L\circ \sigma$ by definition.
		
		Therefore, by using the induction hypothesis, we have 
		\begin{align*}
			\sigma^\imath(\iota(\vartheta_i\cdot x))&=\sigma^\imath(\vartheta_i*x)-(v_i-v_i^{-1})\sigma^\imath(\partial_{\tau i}^L(x)\cdot h_{\tau i})
			\\
			&=\sigma^\imath(x)*\vartheta_i-(v_i-v_i^{-1})v^{-(\alpha_{\tau i},\mu-\alpha_{\tau i})}\sigma^\imath(h_{\tau i}*\partial_{\tau i}^L(x))
			\\
			&=\sigma(x)*\vartheta_i-(v_i-v_i^{-1})v^{-(\alpha_{\tau i},\mu-\alpha_{\tau i})}\sigma(\partial_{\tau i}^L(x))*h_i
			\\
			&=\sigma(x)*\vartheta_i-(v_i-v_i^{-1})v^{-(\alpha_{i},\mu-\alpha_{\tau i})}h_i\sigma(\partial_{\tau i}^L(x))
			\\
			&=\sigma(x)*\vartheta_i-(v_i-v_i^{-1})\sigma(\partial_{\tau i}^L(x))\cdot h_i
			\\
			&= \iota(\sigma(\vartheta_i\cdot x)).
		\end{align*}
		The proof is completed.   
	\end{proof}
	
	\begin{proposition}
		\label{prop:dualCB-anti-inv}
		The dual canonical basis $\{\K_\alpha\diamond C_b\mid \alpha\in\N^\I,b\in\mathbf{C}\}$ of $\hB^\imath_\tau$ is preserved by the anti-involution $\sigma^\imath$. Moreover, $$ \sigma^{\imath}(C_{\alpha,b})= C_{\tau\alpha,\sigma(b)},\quad \forall \alpha\in\N^\I,b\in\mathbf{C}.
		$$
	\end{proposition}
	Similar results hold for $\tB_\tau^\imath$. 
	
	\begin{proof}
		By \cite[Lemma 3.5]{BG17a}, we have $\sigma(b)\in\mathbf{C}$ for  $b\in\mathbf{C}$. Then by Lemma \ref{lem:com-anti-involution}, we have 
		\begin{align*}
			\sigma^\imath(\K_\alpha\diamond \iota(b))&=v^{\frac{1}{2}(\tau \alpha-\alpha,\wt(b))}\sigma^\imath(\K_\alpha*\iota(b))
			\\
			&=v^{\frac{1}{2}(\tau \alpha-\alpha,\wt(b))}\sigma(\iota(b))*\K_{\tau \alpha}
			\\
			&=v^{\frac{1}{2}( \alpha-\tau \alpha,\wt(b))}\K_{\tau \alpha}*\sigma(\iota(b))
			\\
			&=\K_{\tau \alpha}\diamond\sigma(\iota(b)).
		\end{align*}
		Applying $\sigma^{\imath}$ to \eqref{eq:dualCB} and using the above identity, we obtain
		\begin{align*}
			\sigma^{\imath}(C_{\alpha,b})-\K_{\tau \alpha}\diamond \sigma(\iota(b))\in\sum_{(\beta,b')}v^{-1}\Z[v^{-1}]\cdot \K_\beta\diamond \iota(b').
		\end{align*}
		Since $\sigma^\imath$ commutes with the bar involution, $\sigma^\imath(C_{\alpha,b})$ is bar-invariant.
		By the characterization of the dual canonical basis element $C_{\tau\alpha,\sigma(b)}$ in Theorem \ref{thm:dCB}, we must have $\sigma^{\imath}(C_{\alpha,b})= C_{\tau\alpha,\sigma(b)}$. 
	\end{proof}
	
	\begin{remark}
		We conjecture the following positivity property:
		\begin{align}
			\label{conj:CB-dualCB}
			\K_\alpha\diamond\iota(b)\in C_{\alpha,b}+\sum_{(\alpha,b)\prec(\beta,b')}v^{-1}\N[v^{-1}]\cdot C_{\beta,b'},
			\qquad\text{ for } b\in \mathbf{C}, \alpha \in \Z^\I,
		\end{align}
		in any iquantum group $\hB^\imath_\tau$ of finite type; see \cite[Conjecture 1.21]{BG17a} for a similar conjecture on quantum groups. This conjecture is proved in \cite[Lemma 7.14]{LP25} for quantum groups and iquantum groups of type ADE (except type ${\rm AIII}_{2r}$).
	\end{remark}
	
	\subsection{Dual canonical basis via dual PBW}

	Define a partial order $\preceq$ on $\N^\I\times\N^N$ by declaring $(\alpha,\ba)\preceq(\beta,\ba')$ if 
	\begin{itemize}
		\item[(1)] $\alpha+\tau\alpha+\wt(\iota(\vartheta_{\bi,\ba}))=\beta+\tau\beta+\wt(\iota(\vartheta_{\bi,\ba'}))$, and
		\item[(2)] $0\neq\beta-\alpha\in\N^\I$ or $\beta=\alpha$ and $\ba'\preceq\ba$ (here $\preceq$ is the partial order on $\N^N$ defined in \cite[Section 4.2]{BG17b}).
	\end{itemize}
	
	We shall use the dual PBW basis $\{\K_\alpha\diamond \iota(\vartheta_{\bi,\ba})\mid \alpha\in\N^\I,\ba\in\N^N\}$ to give a second construction of the dual canonical basis $\{\K_\alpha\diamond C_{b_{\bi,\ba}}\mid \alpha\in\N^\I,\ba\in\N^N\}$ of $_{\mathcal{Z}}\hB^\imath_\tau$ (see Theorem~ \ref{thm:dCB}).
	
	\begin{proposition} \label{prop:dCB:PBW}
		For each $\alpha\in\N^\I$ and $\ba\in\N^N$, there exists a unique element $\K_\alpha\diamond C_{\vartheta_{\bi,\ba}}$ in $_{\mathcal{Z}}\hB^\imath_\tau$ satisfying $\ov{\K_\alpha\diamond C_{\vartheta_{\bi,\ba}}}=\K_\alpha\diamond C_{\vartheta_{\bi,\ba}}$ and
		\begin{align}
			\label{eqn:dCBviaPBW}
			\K_\alpha\diamond C_{\vartheta_{\bi,\ba}}-\K_\alpha\diamond \iota(\vartheta_{\bi,\ba})\in \sum_{(\alpha,\ba)\prec(\beta,\ba')}v^{-1}\Z[v^{-1}]\cdot \K_\beta\diamond \iota(\vartheta_{\bi,\ba'}).
		\end{align}
		Moreover, we have
		$\K_\alpha\diamond C_{\vartheta_{\bi,\ba}} =\K_\alpha\diamond C_{b_{\bi,\ba}}$.
	\end{proposition}
	
	\begin{proof}
		Plugging \eqref{eq: QG dPBW and dCB} into \eqref{eq:bar_dCB}, we have 
		\begin{align}
			\label{eq:PBW-bar}
			\ov{\K_\alpha\diamond \iota(\vartheta_{\bi,\ba})}-\K_\alpha\diamond\iota(\vartheta_{\bi,\ba})\in \sum_{(\alpha,\ba)\prec(\beta,\ba')}\Z[v,v^{-1}]\cdot \K_\beta\diamond \iota(\vartheta_{\bi,\ba'}).
		\end{align}
		Hence the existence and uniqueness of the desired element $\K_\alpha\diamond C_{\vartheta_{\bi,\ba}}$ follows from Lusztig's Lemma.
		
		On the other hand, it follows from \eqref{eq: QG dPBW and dCB} and Theorem~\ref{thm:dCB} that the element $\K_\alpha\diamond C_{b_{\bi,\ba}}$ from Theorem \ref{thm:dCB} also satisfies the requirement \eqref{eq:PBW-bar}. Hence, $\K_\alpha\diamond C_{\vartheta_{\bi,\ba}}=\K_\alpha\diamond C_{b_{\bi,\ba}}$ by uniqueness.
	\end{proof}
	
	We can view $\tau\in\Inv(C)$ as an involution of $\hB_\tau^\imath$ (and also $\hB$), which maps $\vartheta_i\mapsto \vartheta_{\tau i}$, $h_i\mapsto h_{\tau i}$ for $i\in\I$. Similarly to Lemma \ref{lem:com-anti-involution}, one can see $\iota\circ \tau=\tau\circ\iota$.
	
	\begin{corollary}
		\label{cor:dCB-tau}
		The dual canonical basis $\{\K_\alpha\diamond C_b\mid \alpha\in\N^\I,b\in\mathbf{C}\}$ of $\hB^\imath_\tau$ is preserved by the involution $\tau$. 
	\end{corollary}
	
	\begin{proof}
		Keep the notation as in Proposition~ \ref{prop:dCB:PBW}. Set $\tau\bi:=(\tau i_1,\dots,\tau i_N)$ for $\bi=(i_1,\dots,i_N)$.
		Note that $\tau(\widetilde{T}_i(\vartheta_j))=\widetilde{T}_{\tau i}(\vartheta_{\tau j})$ for any $i\neq j\in\I$. 
		Then we have $\tau(\vartheta_{\bi,\ba})=\vartheta_{\tau\bi,\ba}$ for any $\ba\in\N^N$.
		Obviously, $\{\vartheta_{\tau \bi,\ba}\mid \ba\in\N^N \}$ is a basis of $\ff_\mathcal{Z}$. 
		We denote by $\{\K_\alpha\diamond C_{b_{\bi,\ba}}\mid \ba\in\N^N,\alpha\in\N^\I\}$ the dual canonical basis constructed from Proposition~ \ref{prop:dCB:PBW} by using the PBW basis $\{\vartheta_{\tau \bi,\ba}\mid \ba\in\N^N \}$. Applying $\tau$ to \eqref{eqn:dCBviaPBW}, we obtain $\tau(\K_\alpha\diamond C_{\vartheta_{\bi,\ba}})=\K_{\tau \alpha}\diamond C_{\vartheta_{\tau\bi,\ba}}$ by the uniqueness.
	\end{proof}

	\begin{remark}
		A new construction of dual canonical bases (cf. \cite{LW21b}) for the universal iquantum groups $\hUi$ and $\tUi$ of type ADE (except type ${\rm AIII}_{2r}$) was given in \cite{LP25} using the (dual) Hall bases of Hall algebras. As the Hall bases coincide with special PBW bases (see \cite{Rin96}), the dual canonical bases constructed in this paper recover those constructed in ADE type geometrically {\em loc. cit.}, and extend to the ${\rm AIII}_{2r}$ type and all non-ADE types. 
	\end{remark}

	\subsection{iBraid group symmetries}
	
	The main result of this subsection is the following.
	
	\begin{theorem}
		\label{thm:dCB-braid}
		The dual canonical basis of $\tB_\tau^\imath$ is preserved by the ibraid group symmetries.
	\end{theorem}
	
	\begin{proof}
		It suffices to show that the dual canonical basis $\{\K_\alpha \diamond C_b \}$ of $\tB_\tau^\imath$ in Theorem \ref{thm:dCB} is preserved by the action of $\widetilde{\TT}_i$, for any given $i\in \I$. 
		By the explicit action of $\widetilde{\TT}_i$ on $\K_\alpha$, this reduces to checking that $\widetilde{\TT}_i(C_b)$ is a dual canonical basis element.
		
		Recall $r_i$ from \eqref{def:ri}. Let $\ell$ be the length of $r_i$ in $W$. We can take a reduced expression $\bi=(i_1,\dots,i_N)$ of $w_0$ such that $r_i =s_{i_1}\cdots s_{i_\ell}$ (which clearly is a reduced expression of $r_i$). Recalling $\vartheta_{\bi,k}^+$ from \eqref{Eik} and noting also $r_i =s_{i_\ell} \cdots s_{i_1}$, we have
		\[\widetilde{T}_{r_i}(\vartheta_{\bi,k}^+)=\widetilde{T}_{r_i}\widetilde{T}_{i_1}^{-1}\cdots \widetilde{T}_{i_\ell}^{-1}\widetilde{T}_{i_{\ell+1}}^{-1}\cdots \widetilde{T}_{i_{k-1}}^{-1}(E_{i_k})=\widetilde{T}_{i_{\ell+1}}^{-1}\cdots \widetilde{T}_{i_{k-1}}^{-1}(E_{i_k})\in\U^+,\quad \forall \ell< k\leq N,\]
		and in particular, $\vartheta_{\bi,k}\in\ff[i,\tau i]$. Denoting $w'=s_{i_{\ell+1}}\cdots s_{i_N}$ (reduced of length $N-\ell$), we have $w_0=r_iw'=w'\cdot w'^{-1}r_iw'$. It follows that the length of $w'^{-1}r_iw'$ is $\ell$. 
		
		Recall from \eqref{taui} the involution $\tau_i$ on $\I_i$; note that $\I_i =\{i_1,\dots,i_\ell\}$ as a set. Denote $j_a=\tau_0\tau_i(i_a)$, for $1\le a \le \ell$. We have 
		\[
		w'^{-1}s_{i_a}w'
		=w_0^{-1}(r_i^{-1}s_{i_a}r_i)w_0 =w_0^{-1}s_{\tau_i(i_a)}w_0=s_{\tau_0\tau_i(i_a)}, 
		\]
		and thus, $w'^{-1}r_iw'=s_{j_1}\cdots s_{j_\ell}=r_{\tau_0(i)}$, which is reduced for length reason. Set 
		\[
		\bi'=(i_{\ell+1},\dots,i_N,j_1,\dots,j_\ell). 
		\]
		Then it is clear that 
		\[
		E_{\bi',k}=\widetilde{T}_{r_i}(E_{\bi,k+\ell}),\quad\text{ for } 1\leq k\leq N-\ell.
		\] 
		Moreover, for $N-\ell+1\leq k\leq N$, we have
		\begin{align*}
			E_{\bi',k}&=\widetilde{T}_{i_{\ell+1}}^{-1}\cdots\widetilde{T}_{i_{N}}^{-1}(\widetilde{T}_{j_1}^{-1}\cdots \widetilde{T}_{j_{k-N+\ell-1}}^{-1})(E_{j_{k-N+\ell}})\\
			&=\widetilde{T}_{i_{\ell+1}}^{-1}\cdots\widetilde{T}_{i_{N}}^{-1}(\widetilde{T}_{j_1}^{-1}\cdots \widetilde{T}_{j_{k-N+\ell-1}}^{-1})\widetilde{T}_{w'^{-1}}(E_{i_{k-N+\ell}})\\
			&=\widetilde{T}_{i_1}^{-1}\cdots\widetilde{T}_{i_{k-N+\ell-1}}^{-1}(\widetilde{T}_{i_{\ell+1}}^{-1}\cdots\widetilde{T}_{i_{N}}^{-1})\widetilde{T}_{w'^{-1}}(E_{i_{k-N+\ell}})\\
			&=\widetilde{T}_{i_1}^{-1}\cdots\widetilde{T}_{i_{k-N+\ell-1}}^{-1}\widetilde{T}_{w'^{-1}}^{-1}\widetilde{T}_{w'^{-1}}(E_{i_{k-N+\ell}})\\
			&=\widetilde{T}_{i_1}^{-1}\cdots\widetilde{T}_{i_{k-N+\ell-1}}^{-1}(E_{i_{k-N+\ell}})\\
			&=E_{\bi,k-N+\ell}.
		\end{align*}
		
		Write $\ba\in\N^N$ as $\ba =(\ba_1,\ba_2)$ where $\ba_1=(a_1,\dots,a_\ell)$, $\ba_2=(a_{\ell+1},\dots,a_N)$. Note that 
		$\vartheta_{\bi,\ba}=v^{\frac{1}{2}(\wt(\vartheta_{\bi,\ba_1}),\wt(x))}\vartheta_{\bi,\ba_1}x$, where $\vartheta_{\bi,\ba_1}\in\ff_{i,\tau i}$ and  $x:=\vartheta_{\bi,\ba_2} \in\ff[i,\tau i]$. Applying Theorem~ \ref{thm:Tifif_i} gives us
		\begin{align*}
			\widetilde{\TT}_i(\vartheta_{\bi,\ba})&=\widetilde{\TT}_i(v^{\frac{1}{2}(\wt(\vartheta_{\bi,\ba_1}),\wt(x))}\vartheta_{\bi,\ba_1}x)\\
			&=v^{-\frac{1}{2}(\tau\wt(\vartheta_{\bi,\ba_1}),\wt(x))}\K_{\tau_i\wt(\vartheta_{\bi,\ba_1})}^{-1}\diamond \big(\widetilde{T}_{r_i}(x)\tau\tau_i(\vartheta_{\bi,\ba_1})\big)\\
			&=v^{\frac{1}{2}(\tau\tau_i\wt(\vartheta_{\bi,\ba_1}),r_i(\wt(x)))}\K_{\tau_i\wt(\vartheta_{\bi,\ba_1})}^{-1}\diamond \big(\widetilde{T}_{r_i}(x)\tau\tau_i(\vartheta_{\bi,\ba_1})\big)\\
			&=v^{\frac{1}{2}(\tau\tau_i\wt(\vartheta_{\bi',\ba_1}),\wt(\vartheta_{\bi',\ba_2}))}\K_{\tau_i\wt(\vartheta_{\bi,\ba_1})}^{-1}\diamond \big(\vartheta_{\bi',\ba_2}\tau\tau_i(\vartheta_{\bi',\ba_1})\big).
		\end{align*}
		Note that $\tau\tau_i$ is either the identity on $\I_i$, or exchanges $i$ with $\tau i$ exactly when $c_{i,\tau i}=0$. Therefore, we can always write $\tau\tau_i(\vartheta_{\bi',\ba_1})=\vartheta_{\bi',\ba_1'}$ for some $\ba_1'=(a_1',\dots,a_\ell')$. This allows us to simplify the right-hand side above as 
		\begin{align}
			\widetilde{\TT}_i(\vartheta_{\bi,\ba})&=v^{\frac{1}{2}(\tau\tau_i\wt(\vartheta_{\bi',\ba_1}),\wt(\vartheta_{\bi',\ba_2}))}\K_{\tau_i\wt(\vartheta_{\bi,\ba_1})}^{-1}\diamond \big(\vartheta_{\bi',\ba_2}\tau\tau_i(\vartheta_{\bi',\ba_1})\big)
			\notag \\
			&=v^{\frac{1}{2}(\wt(\vartheta_{\bi',\ba_1'}),\wt(\vartheta_{\bi',\ba_2}))}\K_{\tau_i\wt(\vartheta_{\bi,\ba_1})}^{-1}\diamond(\vartheta_{\bi',\ba_2}\vartheta_{\bi',\ba_1'})
			\notag \\
			&=\K_{\tau_i\wt(\vartheta_{\bi,\ba_1})}^{-1}\diamond \vartheta_{\bi',\ba'},
			\label{Ti:PBW}
		\end{align}
		where $\ba'= (\ba_2,\ba_1') =(a_{\ell+1},\dots,a_N,a_1',\dots,a_{\ell}')$. 
		
		Recall from Theorem \ref{thm:dCB} and Proposition \ref{prop:dCB:PBW} that the dual canonical basis element $C_{\vartheta_{\bi,\ba}}$ is the unique bar invariant element such that 
		$C_{\vartheta_{\bi,\ba}} \in \vartheta_{\bi,\ba} + \sum_{(\beta,\tilde\ba)}v^{-1}\Z[v^{-1}]\cdot \K_\beta\diamond \vartheta_{\bi,\tilde\ba}$; write $\tilde{\ba} =(\tilde{\ba}_1,\tilde{\ba}_2)$ similarly as above for $\ba$. 
		Applying $\widetilde{\TT}_i$ (which commutes with the bar involution) to $C_{\vartheta_{\bi,\ba}}$ and using \eqref{Ti:PBW}, we see that $\widetilde{\TT}_i(C_{\vartheta_{\bi,\ba}})$ is bar invariant and 
		\[
		\widetilde{\TT}_i(C_{\vartheta_{\bi,\ba}}) \in \K_{\tau_i\wt(\vartheta_{\bi,\ba_1})}^{-1} \diamond \vartheta_{\bi',\ba'} + \sum_{(\beta,\tilde\ba)}v^{-1}\Z[v^{-1}]\cdot \K_{r_i\beta}\K_{\tau_i\wt(\vartheta_{\bi,\tilde\ba_1})}^{-1} \diamond \vartheta_{\bi,\tilde\ba}.
		\]
		It follows by the uniqueness that $\widetilde{\TT}_i(C_{\vartheta_{\bi,\ba}})$ is a dual canonical basis element and 
		\[
		\widetilde{\TT}_i(C_{\vartheta_{\bi,\ba}})=\K_{\tau_i\wt(\vartheta_{\bi,\ba_1})}^{-1}\diamond C_{\vartheta_{\bi',\ba'}}.
		\]
		The theorem is proved.
	\end{proof}
	
	\begin{remark}
		Using iHall algebras developed in \cite{LW22a,LW23}, Lu and Pan \cite[Theorem 7.19]{LP25} proved that the dual canonical basis is preserved under the braid group action for universal iquantum groups of type ADE (except for the type ${\rm AIII}_{2r}$ listed in \cite[Table 3.1]{CLPRW}). Moreover, it was shown in \cite[Theorem C]{LP25} that the structure constants of the dual canonical bases are positive (i.e., belonging in $\N[v^{\pm1/2}]$). We conjecture the positivity also holds for quasi-split type ${\rm AIII}_{2r}$.  
	\end{remark}

	\subsection{Dual canonical basis on iquantum groups} 
	
	Recall the distinguished parameter $\bvs_{\diamond}=(\vs_{i,\diamond})$ from \eqref{eq:disting-para}. By \cite[Proposition 6.2(1)]{LW22a}, 
	there exists an algebra epimorphism $\pi_{\bvs_\diamond}:\tB_\tau^\imath\rightarrow \Ui_{\bvs_\diamond}$
	by sending
	\begin{align}
		\vartheta_i\mapsto B_i,\qquad &\K_j\mapsto \begin{cases}
			k_j&\text{ if }\tau j\neq j
			\\
			1&\text{ if }\tau j=j,
		\end{cases}
		\qquad \K_{\tau j}\mapsto \begin{cases}
			k_j^{-1}&\text{ if }\tau j\neq j
			\\
			1&\text{ if }\tau j=j,
		\end{cases}
	\end{align}
	for any $i\in\I$, $j\in\I\setminus\I_\tau$. The kernel of $\pi_{\bvs_\diamond}$ is generated by $\K_i\K_{\tau i}-1$, $i\in\I$. 
	
	\begin{lemma}
		\label{lem:bar-Ui}
		$\Ui_{\bvs_\diamond}$ admits an anti-involution (called bar-involution) such that $\ov{v^{1/2}}=v^{-1/2}$, $\ov{B_i}=B_i$, and $\ov{k_i}=k_i$,  for $i\in\I\setminus\I_\tau$.
	\end{lemma}
	
	\begin{proof}
		Follows by noting that the bar-involution of $\tB^\imath_\tau$ preserves the kernel of $\pi_{\bvs_\diamond}$.
	\end{proof}
	
	Let $i\in\I$. As the ibraid group symmetry $\tTT_i$ of $\tB^\imath_\tau=\tUi$ preserves the kernel of $\pi_{\bvs_\diamond}$, it induces an automorphism $\TT_i$ of $\Ui_{\bvs_\diamond}$; cf. \cite[Proposition 7.2]{LW21a}; that is, we have the following commutative diagram:
	\begin{equation}
		\label{eq:TTi} 
		\xymatrix{
			\tUi\ar[rr]^{\tTT_i} \ar[d]^{\pi_{\bvs_\diamond}}&& \tUi\ar[d]^{\pi_{\bvs_\diamond}}\\
			\Ui_{\bvs_\diamond}\ar[rr]^{\TT_i}&&\Ui_{\bvs_\diamond}
		}
	\end{equation}
	
	For any $\alpha=\sum_ia_i\alpha_i\in\Z^{\I\setminus\I_\tau}$ and $x\in\Ui_{\bvs_\diamond}$, denote
	\[
	k_\alpha=\prod_i k_i^{a_i},
	\qquad
	k_\alpha\diamond x=v^{\frac{1}{2}(\tau\alpha-\alpha,\wt(x'))}\cdot k_\alpha  x,
	\]
	where $x'$ is any preimage of $x$ in $\tUi$ under $\pi_{\bvs_\diamond}$.
	
	\begin{proposition}
		\label{prop:dCB:parameter}
		For each $\alpha\in\Z^{\I\setminus\I_\tau}$ and $b\in \mathbf{C}$, there exists a unique element $C^\diamond_{\alpha,b}\in\Ui_{\bvs_\diamond}$ such that $\ov{C^\diamond_{\alpha,b}}=C^\diamond_{\alpha,b}$ and
		\begin{align}
			\label{eq:dualCB-distinguished}
			C^\diamond_{\alpha,b}-k_\alpha\diamond \pi_{\bvs_\diamond}(\iota(b))\in\sum_{(\gamma,b') \in \Z^{\I\setminus\I_\tau} \times \mathbf{C}}v^{-1}\Z[v^{-1}]\cdot k_\gamma\diamond \pi_{\bvs_\diamond}(\iota(b')).
		\end{align}
		Moreover, 
		$\{C_{\alpha,b}^\diamond\mid \alpha\in\Z^{\I\setminus\I_\tau},b\in\mathbf{C}\}$ forms a basis of $\Ui_{\bvs_{\diamond}}$, and it is preserved by the braid group symmetries $\TT_i$, for $i\in \I_\tau$.   
	\end{proposition}
	
	\begin{proof}
		Recall the dual canonical basis $\{C_{\tilde\alpha,b}\}$ from Theorem \ref{thm:dCB}. We set $C_{\alpha,b}^\diamond:=\pi_{\bvs_\diamond}(C_{\alpha,b})$, for $\alpha\in\Z^{\I\setminus\I_\tau}$ and $b\in \mathbf{C}$, which clearly satisfies $\ov{C^\diamond_{\alpha,b}}=C^\diamond_{\alpha,b}$. For $\beta=\sum_{i\in\I} a_i\alpha_i$, we have $\pi_{\bvs_\diamond}(\K_\beta)=k_\gamma$, where $\gamma=\sum_{i\in\I\setminus\I_\tau}(a_i-a_{\tau i})\alpha_i$. Then \eqref{eq:dualCB-distinguished} follows from \eqref{eq:dualCB}. The uniqueness of $C_{\alpha,b}^\diamond$ follows by a standard argument. 
		
		It is clear that 
		$\{C_{\alpha,b}^\diamond\mid \alpha\in\Z^{\I\setminus\I_\tau},b\in\mathbf{C}\}$ is a basis of $\Ui_{\bvs_{\diamond}}$.  Finally it follows from Theorem \ref{thm:dCB-braid} and the commutative diagram \eqref{eq:TTi} that this basis is preserved by the braid group symmetry $\TT_i$.
	\end{proof}

	\section{Dual canonical basis for Drinfeld doubles} \label{sec:doubleCB}
	
	In case of iHopf algebra on Borel of diagonal type, our construction yields a dual canonical basis of Drinfeld double quantum group $\hU$. Berenstein and Greenstein earlier defined two bases for $\hU$, called positive/negative double canonical bases. In this section we prove the coincidence of their bases as well as the dual canonical basis of $\hU$ constructed in this paper. Several conjectures from their work quickly follow from such identification. 
	
	\subsection{Results of Berenstein-Greenstein}
	First, following \cite{BG17a} we define the quantum Heisenberg algebras $\mathcal{H}^{\pm}$ by 
	\[\mathcal{H}^+=\hU/\langle K_i'\mid i\in\I\rangle,\quad \mathcal{H}^-=\hU/\langle K_i\mid i\in\I\rangle.\]
	Let $\mathbf{K}^+$ (resp. $\mathbf{K}^-$) be the submonoid of $\hU$ generated by the $K_i$ (resp. the $K_i'$), $i\in\I$. Then we have triangular decompositions
	\[\mathcal{H}^+=\mathbf{K}^+\otimes\U^-\otimes \U^+,\quad \mathcal{H}^-=\mathbf{K}^-\otimes\U^+\otimes \U^-.\]
	The induced natural embeddings of vector spaces
	\begin{align*}
		\iota_+&:\mathcal{H}^+=\mathbf{K}^+\otimes\U^-\otimes \U^+\hookrightarrow \hU=\mathbf{K}^-\otimes(\mathbf{K}^+\otimes\U^-\otimes \U^+),\\
		\iota_-&:\mathcal{H}^-=\mathbf{K}^-\otimes\U^+\otimes \U^-\hookrightarrow \hU=\mathbf{K}^+\otimes(\mathbf{K}^-\otimes\U^+\otimes \U^-)
	\end{align*}
	split the canonical projections $\hU\rightarrow\mathcal{H}^+$ and $\hU\rightarrow\mathcal{H}^-$, respectively.
	
	Let $\mathbf{C}^{\pm}$ be the rescaled dual canonical basis of $\hU^{\pm}$ defined in \S\ref{subsec:dCB:f}. Recall from Lemma~\ref{lem:anti-involut-QG} that there is a bar-involution on $\hU$ defined by $\ov{v}=v^{-1}$, $\ov{E}_i=E_i$, $\ov{F}_i=F_i$ and $\ov{K}_i=K_i$, $\ov{K}'_i=K_i'$ for $i\in\I$.
	
	Let $\alpha_{+i}=(\alpha_i,0)$, $\alpha_{-i}=(0,\alpha_i)$. We define a weight function on $\hU$ by setting 
	\[
	\wt^{2}(E_i)=\alpha_{+i},\quad \wt^{2}(F_i)=\alpha_{-i},\quad \wt^{2}(K_i)=\wt^{2}(K_i')=\alpha_{+i}+\alpha_{-i}.
	\]
	It is easily seen that $\hU$ becomes a $\N^{\I^2}$-graded algebra. Moreover, we have a partial order $\prec$ on $\N^{\I^2}$ defined by $\alpha\prec\beta$ if and only if $\beta-\alpha\in\N^\I\times\N^\I$. Using this degree function we define an action $\diamond$ of the algebra $\hU^0$ on $\hU$ via
	\begin{align*}
		K_i\diamond x=v^{-\frac{1}{2}\check{\alpha}_i(\wt^{2}(x))}K_ix,\quad K_i'\diamond x=v^{\frac{1}{2}\check{\alpha}_i(\wt^{2}(x))}K_i'x
	\end{align*}
	where $\check{\alpha}_i\in\Hom_\Z(\N^{\I^2},\Z)$ is defined by $\check{\alpha_i}(\alpha_{\pm i})=\pm c_{ij}$ and $x\in\hU$ is homogeneous. This action is characterized by the following property: 
	\begin{equation}\label{QG diamond action char}
		\ov{K\diamond x}=K \diamond \ov{x},\quad  K\in\hU^0,x\in\hU.
	\end{equation}
	Note that the $\diamond$-action as well as the bar-involution factors through to a $\mathbf{K}^{\pm}$-action and an anti-involution on $\mathcal{H}^{\pm}$ via the canonical projection $\hU\rightarrow\mathcal{H}^{\pm}$, and (\ref{QG diamond action char}) still holds. 
	
	The following results are due to \cite{BG17a}.
	
	\begin{proposition}[\text{\cite[Theorem 1.3]{BG17a}}]
		\label{prop:doubleCB-H+}
		For any $(b_+,b_-)\in\mathbf{C}^+\times\mathbf{C}^-$, there is a unique element $b_-\circ b_+\in\mathcal{H}^+$ fixed by $\bar{\cdot}$ and satisfying
		\[b_-\circ b_+-b_-b_+\in\sum v\Z[v]K\diamond(b'_-b'_+)\]
		where the sum is taken over $K\in\mathbf{K}^+\setminus\{1\}$ and $b'_{\pm}\in\mathbf{C}^{\pm}$ such that $\wt^{2}(b_-b_+)=\wt^{2}(K)+\wt^{2}(b'_-b'_+)$. The basis $\{K\diamond(b_-\circ b_+)\mid K\in\mathbf{K}^+,b_{\pm}\in\mathbf{C}^{\pm}\}$ is called the double canonical basis of $\mathcal{H}^+$.
	\end{proposition} 
	
	\begin{theorem}[\text{\cite[Theorem 1.5]{BG17a}}]
		\label{thm:doubleCB-U +}
		For any $(b_+,b_-)\in\mathbf{C}^+\times\mathbf{C}^-$, there is a unique element $b_-\bullet b_+\in\hU$ fixed by $\bar{\cdot}$ and satisfying
		\[b_-\bullet b_+-\iota_+(b_-\circ b_+)\in \sum v^{-1}\Z[v^{-1}]K\diamond\iota_+(b'_-\circ b'_+)\]
		where the sum is taken over $K\in\hU^0\setminus\mathbf{K}^+$ and $b'_{\pm}\in\mathbf{C}^{\pm}$ such that $\wt^{2}(b_-b_+)=\wt^{2}(K)+\wt^{2}(b'_-b'_+)$. The basis $\{K\diamond(b_-\bullet b_+)\mid K\in\hU^0,b_{\pm}\in\mathbf{C}^{\pm}\}$ is called the positive double canonical basis of $\hU$.
	\end{theorem}
	
	The following are variants of Proposition~\ref{prop:doubleCB-H+} and Theorem \ref{thm:doubleCB-U +}.
	
	\begin{proposition}
		[cf. \text{\cite{BG17a}}]
		\label{prop:doubleCB-H-}
		For any $(b_+,b_-)\in\mathbf{C}^+\times\mathbf{C}^-$, there is a unique element $b_+\circ b_-\in\mathcal{H}^-$ fixed by $\bar{\cdot}$ and satisfying
		\[b_+\circ b_--b_+b_-\in\sum v\Z[v]K\diamond(b'_+b'_-)\]
		where the sum is taken over $K\in\mathbf{K}^-\setminus\{1\}$ and $b'_{\pm}\in\mathbf{C}^{\pm}$ such that $\wt^{2}(b_+b_-)=\wt^{2}(K)+\wt^{2}(b'_+b'_-)$. The basis $\{K\diamond(b_+\circ b_-)\mid K\in\mathbf{K}^-,b_{\pm}\in\mathbf{C}^{\pm}\}$ is called the double canonical basis of $\mathcal{H}^-$.
	\end{proposition} 
	
	\begin{theorem}[cf. \text{\cite{BG17a}}]
		\label{thm:doubleCB-U -}
		For any $(b_+,b_-)\in\mathbf{C}^+\times\mathbf{C}^-$, there is a unique element $b_+\bullet b_-\in\hU$ fixed by $\bar{\cdot}$ and satisfying
		\[b_+\bullet b_--\iota_-(b_+\circ b_-)\in \sum v^{-1}\Z[v^{-1}]K\diamond\iota_-(b'_+\circ b'_-)\]
		where the sum is taken over $K\in\hU^0\setminus\mathbf{K}^-$ and $b'_{\pm}\in\mathbf{C}^{\pm}$ such that $\wt^{2}(b_+b_-)=\wt^{2}(K)+\wt^{2}(b'_+b'_-)$. The basis $\{K\diamond(b_+\bullet b_-)\mid K\in\hU^0,b_{\pm}\in\mathbf{C}^{\pm}\}$ is called the negative double canonical basis of $\hU$.
	\end{theorem}
	It was conjectured by Bernstein-Greenstein (see \cite[Conjecture 1.11]{BG17a} and \cite[Remark 1.12]{BG17a}) that positive and negative canonical bases coincide.

	\subsection{Coincidence of bases on $\hU$}
	
	By identifying $\hU$ with the iHopf algebra $(\hB\otimes\hB )^\imath$ of \S\ref{subsec:iHopf:ff}, we have constructed the dual canonical basis of $\hU$. For $(b_+,b_-)\in \mathbf{C}^+\times\mathbf{C}^-$, we denote by $C_{b_+,b_-}$ the corresponding dual canonical basis, characterized by the following properties:
	\[\ov{C_{b_+,b_-}}=C_{b_+,b_-},\quad C_{b_+,b_-}\in b_+\otimes b_-+\sum v^{-1}\Z[v^{-1}]K\diamond(b_+'\otimes b_-'),\]
	where $K\in \tU^0$ and $(b_+',b_-')\in \mathbf{C}^+\times\mathbf{C}^-$. We need the following lemma:
	
	\begin{lemma}\label{lem: H^pm embedding on iHopf}
		For any elements $x,y\in\hB$, we denote by $x\otimes^{\pm}y$ the image of $x\otimes y\in\hU$ in $\mathcal{H}^{\pm}$. Then
		\[\iota_+(x\otimes^+ y)=\iota_-(x\otimes^- y)=x\otimes y.\]
	\end{lemma}
	\begin{proof}
		Let us denote by $\ast^\pm$ the multiplication of $\mathcal{H}^\pm$, and by $\ast$ the multiplication of $\hU$. Then the definition of $\iota_{\pm}$ spells out as
		\begin{equation}\label{eq: H^pm embedding def}
			\iota_-(x^+\ast^-y^-)=x^+\ast y^-,\quad \iota_+(y^-\ast^+x^+)=y^-\ast x^+.
		\end{equation}
		We now prove the lemma by induction on $\wt^{2}(x\otimes y)$. The lemma is clear if $x=1$ or $y=1$. Now assume that claim for $x',y'\in\ff$ with $\wt^{2}(x'\otimes y')\prec\wt^{2}(x\otimes y)$. The definition of $\mathcal{H}^\pm$ implies
		\begin{align*}
			x^+\ast^-y^-&=x\otimes^- y+\sum_{y_{(2)}\neq 1}\varphi(x_{(1)},y_{(2)}) x_{(2)}\otimes^- y_{(1)}h_{\wt(y_{(2)})},\\
			y^-\ast^+ x^+&=x\otimes^+ y+\sum_{y_{(1)}\neq 1}\varphi(x_{(2)},y_{(1)}) x_{(1)}h_{\wt(x_{(2)})}\otimes^+ y_{(2)}.
		\end{align*}
		From \eqref{eq: H^pm embedding def}, we then get
		\begin{equation}\label{eq: H^pm embedding on iHopf}
			\begin{aligned}
				x^+\ast y^-&=\iota_-(x^+\ast^-y^-)=\iota_-(x\otimes^- y)+\sum_{y_{(2)}\neq 1}\varphi(x_{(1)},y_{(2)}) \iota_-(x_{(2)}\otimes^- y_{(1)}h_{\wt(y_{(2)})}),\\
				y^-\ast x^+&=\iota_+(y^-\ast^+ x^+)=\iota_+(x\otimes^+ y)+\sum_{y_{(1)}\neq 1}\varphi(x_{(2)},y_{(1)}) \iota_+(x_{(1)}h_{\wt(x_{(2)})}\otimes^+ y_{(2)}).
			\end{aligned}
		\end{equation}
		The induction hypothesis implies 
		\begin{align*}
			\iota_-(x_{(2)}\otimes^- y_{(1)}h_{\wt(y_{(2)})})=x_{(2)}\otimes^- y_{(1)}h_{\wt(y_{(2)})},\quad\forall y_{(2)}\neq 1,\\
			\iota_+(x_{(1)}h_{\wt(x_{(2)})}\otimes^+ y_{(2)})=x_{(1)}h_{\wt(x_{(2)})}\otimes^+ y_{(2)},\quad\forall y_{(1)}\neq 1.
		\end{align*}
		Comparing \eqref{eq: H^pm embedding on iHopf} with the definitions of $x^+\ast y^-$ and $y^-\ast x^+$ then gives $\iota_+(x\otimes^+y)=\iota_-(x\otimes^-y)=x\otimes y$. This completes the induction step.
	\end{proof}
	
	Because of Lemma~\ref{lem: H^pm embedding on iHopf}, we will not distinguish $x\otimes y$ with its images in $\mathcal{H}^\pm$. The multiplication formulas of $\mathcal{H}^\pm$ then become
	\begin{align}
		x^+\ast^-y^-&=x\otimes y+\sum_{y_{(2)}\neq 1}\varphi(x_{(1)},y_{(2)}) x_{(2)}\otimes y_{(1)}h_{\wt(y_{(2)})},\label{eq: H^pm multiplication-1}\\
		y^-\ast^-x^+&=x\otimes y,\label{eq: H^pm multiplication-2}\\
		x^+\ast^+y^-&=x\otimes y,\label{eq: H^pm multiplication-3}\\
		y^-\ast^+x^+&=x\otimes y+\sum_{y_{(1)}\neq 1}\varphi(x_{(2)},y_{(1)}) x_{(1)}h_{\wt(x_{(2)})}\otimes y_{(2)}.\label{eq: H^pm multiplication-4}
	\end{align}
	
	We can now formulate the main result of this section. 
	
	\begin{theorem}
		\label{thm:double=dCB}
		The positive double canonical basis, the negative double canonical basis, and the dual canonical basis coincide with each other. More explicitly, for any $(b_+,b_-)\in\mathbf{C}^+\times\mathbf{C}^-$, we have $b_-\bullet b_+=b_+\bullet b_-=C_{b_+,b_-}$.
	\end{theorem}
	\begin{proof}
		By comparing the definition of $C_{b_+,b_-}$ with Theorems~\ref{thm:doubleCB-U +} and \ref{thm:doubleCB-U -}, it suffices to establish the properties \eqref{eq: double CB triangular -}--\eqref{eq: double CB triangular +} below:
		\begin{align}
			\iota_-(b_+\circ b_-)\in b_+\otimes b_-+\sum_{\alpha\in\N^\I,(b_+',b_-')\in\mathbf{C}^+\times\mathbf{C}^-} v^{-1}\Z[v^{-1}]K_\alpha'\diamond(b_+'\otimes b_-'),\label{eq: double CB triangular -}\\
			\iota_+(b_-\circ b_+)\in b_+\otimes b_-+\sum_{\alpha\in\N^\I,(b_+',b_-')\in\mathbf{C}^+\times\mathbf{C}^-} v^{-1}\Z[v^{-1}]K_\alpha\diamond(b_+'\otimes b_-').\label{eq: double CB triangular +}
		\end{align}
		These clearly hold for $b_+=1$ or $b_-=1$. Now assume that \eqref{eq: double CB triangular -}--\eqref{eq: double CB triangular +} are valid for any $(b_+',b_-')\in\mathbf{C}^+\times\mathbf{C}^-$, where $\wt^{2}(b'_+)+\wt^{2}(b'_-)\prec\wt^{2}(b_+)+\wt^{2}(b_-)$. We note that
		\begin{equation}\label{eq: bar of dCB product-1}
			\begin{aligned}
				\ov{b_+\ast^- b_-}-b_+\ast^-b_-&=b_+\otimes b_--b_+\ast^-b_-\\
				&=-\sum_{(b_-)_{(2)}\neq 1}\varphi\big((b_+)_{(1)},(b_-)_{(2)}\big) \,(b_+)_{(2)}\otimes (b_-)_{(1)}h_{\wt((b_-)_{(2)})}\\
				&=\sum_{\wt^{2}(b'_+\otimes b'_-)+\wt^{2}(K_\alpha')=\wt^{2}(b_+\otimes b_-)} a_{\alpha,b_+',b_-'}K_\alpha'\diamond (b_+'\circ b_-').
			\end{aligned}
		\end{equation}
		where $a_{\alpha,b_+',b_-'}\in\Z[v,v^{-1}]$ can be written uniquely in the form $a_{\alpha,b_+',b_-'}=a^+_{\alpha,b_+',b_-'}-a^-_{\alpha,b_+',b_-'}$ for  $a^+_{\alpha,b_+',b_-'}\in v\Z[v]$ and $a^-_{\alpha,b_+',b_-'} =\ov{a^+_{\alpha,b_+',b_-'}}$. The following element
		\begin{align*}
			b_+\ast^-b_-+\sum a^+_{\alpha,b_+',b_-'}K_\alpha'\diamond (b_+'\circ b_-')
		\end{align*}
		is bar invariant, and it has $b_+\ast^-b_-$ as a leading term and other terms are in $\sum v\Z[v]K_\alpha'\diamond(b_+'\circ b_-')$, by induction. Using Proposition~\ref{prop:doubleCB-H-}, we conclude that 
		\begin{equation}\label{eq: b_+ circ b_- expression}
			b_+\circ b_-=b_+\ast^-b_-+\sum a^+_{\alpha,b_+',b_-'}K_\alpha'\diamond (b_+'\circ b_-').
		\end{equation}
		
		Similarly, we have 
		\begin{equation*} 
			\begin{aligned}
				\ov{b_-\ast^+ b_+}-b_-\ast^+b_+&=b_+\otimes b_--b_-\ast^+b_+\\
				&=-\sum_{(b_-)_{(1)}\neq 1}\varphi((b_+)_{(2)},(b_-)_{(1)})(b_+)_{(1)}\otimes (b_-)_{(2)}h_{\wt((b_-)_{(1)})}\\
				&=\sum_{\wt^{2}(b'_+\otimes b'_-)+\wt^{2}(K_\alpha)=\wt^{2}(b_+\otimes b_-)} c_{\alpha,b_+',b_-'}K_\alpha\diamond (b_+'\circ b_-').
			\end{aligned}
		\end{equation*}
		where $c_{\alpha,b_+',b_-'}\in\Z[v,v^{-1}]$ can be uniquely written as $c_{\alpha,b_+',b_-'}=c^+_{\alpha,b_+',b_-'}-c^-_{\alpha,b_+',b_-'}$ for  $c^+_{\alpha,b_+',b_-'}\in v\Z[v]$ and $c^-_{\alpha,b_+',b_-'} =\ov{c^+_{\alpha,b_+',b_-'}}$. Using Proposition~ \ref{prop:doubleCB-H+}, we conclude that 
		\begin{equation}\label{eq: b_- circ b_+ expression}
			b_-\circ b_+=b_-\ast^+b_++\sum c^+_{\alpha,b_+',b_-'}K_\alpha\diamond (b_+'\circ b_-').
		\end{equation}
		
		Using \eqref{eq: b_+ circ b_- expression} and \eqref{eq: b_- circ b_+ expression}, we can now compute that
		\begin{align*}
			\iota_-(b_+\circ b_-)&=\iota_-(b_+\ast^-b_-)+\sum a^+_{\alpha,b_+',b_-'}\iota_-(K_\alpha'\diamond (b_+'\circ b_-'))\\
			&=\iota_-(b_+\otimes b_-)+\iota_-\Big(\sum_{(b_-)_{(2)}\neq 1}\varphi((b_+)_{(1)},(b_-)_{(2)})(b_+)_{(2)}\otimes (b_-)_{(1)}h_{\wt((b_-)_{(2)})}\Big)\\
			&\quad+\sum a^+_{\alpha,b_+',b_-'}K_\alpha'\diamond \iota_-(b_+'\circ b_-')\\
			&=b_+\otimes b_--\sum a_{\alpha,b_+',b_-'}K_\alpha'\diamond \iota_-(b_+'\circ b_-')+\sum a^+_{\alpha,b_+',b_-'}K_\alpha'\diamond \iota_-(b_+'\circ b_-')\\
			&=b_+\otimes b_-+\sum a^-_{\alpha,b_+',b_-'}K_\alpha'\diamond \iota_-(b_+'\circ b_-'),
		\end{align*}
		and
		\begin{align*}
			\iota_+(b_-\circ b_+)&=\iota_+(b_-\ast^+b_+)+\sum c^+_{\alpha,b_+',b_-'}\iota_-(K_\alpha\diamond (b_+'\circ b_-'))\\
			&=\iota_+(b_-\otimes b_+)+\iota_+\Big(\sum_{(b_-)_{(1)}\neq 1}\varphi((b_+)_{(2)},(b_-)_{(1)})(b_+)_{(1)}h_{\wt((b_+)_{(2)})}\otimes (b_-)_{(2)}\Big)\\
			&\quad+\sum c^+_{\alpha,b_+',b_-'}K_\alpha\diamond \iota_+(b_-'\circ b_+')\\
			&=b_+\otimes b_--\sum c_{\alpha,b_+',b_-'}K_\alpha\diamond \iota_+(b_-'\circ b_+')+\sum c^+_{\alpha,b_+',b_-'}K_\alpha\diamond \iota_+(b_-'\circ b_+')\\
			&=b_+\otimes b_-+\sum c^-_{\alpha,b_+',b_-'}K_\alpha\diamond \iota_+(b_-'\circ b_+').
		\end{align*}
		By induction hypothesis we have
		\[\iota_-(b_+'\circ b_-')\in \sum \Z[v^{-1}]K'_\beta\diamond(b_+''\otimes b_-''),\quad  \iota_+(b_-'\circ b_+')\in \sum \Z[v^{-1}]K_\beta\diamond(b_+''\otimes b_-'').\]
		Since $a^-_{\alpha,b_+',b_-'},c^-_{\alpha,b_+',b_-'}\in v^{-1}\Z[v^{-1}]$, we have established \eqref{eq: double CB triangular -}--\eqref{eq: double CB triangular +}, and hence proved the theorem.
	\end{proof}
	
	Thanks to Theorem \ref{thm:double=dCB}, we do not need to distinguish positive and negative double canonical bases, and will refer to them as double canonical basis of $\tU$. 
	Together with Theorem \ref{thm:dCB-braid}, we have the following corollary, which proves \cite[Conjecture 1.15]{BG17a} for all quantum groups of finite type.
	\begin{corollary}
		\label{cor:dCB:braid}
		The  double (= dual) canonical basis of $\tU$ is preserved by the braid group action.
	\end{corollary}
	
	The Chevalley involution $\omega$ of $\tU$ is the algebra automorphism such that $\omega(E_i)=F_i$, $\omega(F_i)=E_i$, $\omega(K_i)=K_i'$, $\omega(K_i')=K_i$, for $i\in\I$. We can view quantum groups as iquantum groups of diagonal type, where the involution $\tau$ is $\swa$; Example \ref{ex:QGvsiQG}. Since the Chevalley involution $\omega$ coincides with $\swa$, we have the following variant of Corollary~ \ref{cor:dCB-tau} which follows by the same argument.
	
	\begin{corollary}
		\label{cor:dCB-chevalley}
		The double ($=$ dual) canonical basis of $\tU$ is preserved by the Chevalley involution $\omega$. 
	\end{corollary}
	
	From the construction, we see that $\omega$ maps positive double canonical basis to the negative one, and vice versa. Consequently, Corollary \ref{cor:dCB-chevalley} confirms again the coincidence of positive and negative double canonical bases.
	
	Using Proposition \ref{prop:dualCB-anti-inv}, we have the following corollary, which proves \cite[Conjecture 1.11]{BG17a} for quantum groups of finite type.
	
	\begin{corollary}[{\cite[Conjecture 1.11]{BG17a}}]
		\label{cor:dCB:sigma}
		The double (= dual) canonical basis of $\tU$ is preserved by the anti-involution $\sigma$. More precisely, we have
		\[\sigma(K\diamond(b_-\bullet b_+))=\sigma(K)\diamond (\sigma(b_+)\bullet \sigma(b_-))=\sigma(K)\diamond (\sigma(b_-)\bullet \sigma(b_+)).\]
	\end{corollary}
	
	\begin{remark}
		For Drinfeld double quantum groups of type ADE, Theorem \ref{thm:double=dCB}, Corollaries~\ref{cor:dCB:braid} and \ref{cor:dCB:sigma} were established earlier by the second and third authors in \cite{LP25} by entirely different approaches. More explicitly, in the framework of generalized quiver varieties and perverse sheaves, it is proved in \cite[Theorem 8.15]{LP25} that double canonical bases coincide with dual canonical bases for quantum groups $\tU$. Using Hall algebras, it is proved in \cite[Corollary 8.17]{LP25} that the dual canonical basis is preserved under the braid group action. The statement that the dual canonical basis is preserved by $\sigma$ appears in \cite[Proposition 8.20]{LP25}. 
	\end{remark}
	
	\appendix
	\section{Dual canonical bases in rank one}
	\label{sec:rank1}
	
	Among 3 rank one quasi-split (universal) iquantum groups, closed formulas for dual canonical bases were known in 2 rank one cases. In this section, we obtain explicit recursive formulas in the remaining most involved rank one case.

	\subsection{Split and diagonal rank one cases}
	
	The rank of the Satake diagram is the number of the $\tau$-orbits. In this way, we can define the iquantum groups of rank one. The split (universal) iquantum group $\tUi_v(\mathfrak{sl}_2)$ of rank one is associated to $\I$ which consists of a single vertex, and is a commutative algebra. The dual canonical basis of $\tUi_v(\mathfrak{sl}_2)$ is obtained in \cite[Section 9]{LP25}. 
	
	A second iquantum group of rank one is associated to the Satake diagram
	\begin{center}\setlength{\unitlength}{0.7mm}
		\vspace{-.4cm}
		\begin{equation}
			\label{eq:satakerank2}
			\begin{picture}(50,13)(0,0)
				
				\put(-0.5,-6){\small $1$}
				\put(20,-6){\small $2$}
				\put(-0.5,-2){$\circ$}
				\put(20,-2){$\circ$}
				
				%
				\color{purple}
				\qbezier(11,4)(15,3.7)(19.5,1)
				\qbezier(11,4)(7,3.7)(2.5,1)
				\put(19,1.1){\vector(2,-1){0.5}}
				\put(2,1.1){\vector(-2,-1){0.5}}
				\put(10,3){\small $^{\tau}$}
			\end{picture}
		\end{equation}
		\vspace{0.1cm}
	\end{center}
	This iquantum group is isomorphic to the Drinfeld double $\tU_v(\mathfrak{sl}_2)$; see Example \ref{ex:QGvsiQG}. Its dual canonical basis is obtained in \cite[Section 10]{LP25}, cf. \cite[Section 4.1]{BG17a}.
	
	\subsection{Quasi-split rank one}
	
	The remaining iquantum group of rank one is the $\tUi_v(\mathfrak{sl}_3)$, which is associated to the Satake diagram
	\begin{center}\setlength{\unitlength}{0.7mm}
		\vspace{-.4cm}
		\begin{equation}
			\label{eq:satakerank1}
			\begin{picture}(50,13)(0,0)
				\put(-0.5,-6){\small $1$}
				\put(20,-6){\small $2$}
				\put(-0.5,-2){$\circ$}
				\put(20,-2){$\circ$}	

				\put(3,-.5){\line(1,0){16}}
				\color{purple}
				\qbezier(11,4)(15,3.7)(19.5,1)
				\qbezier(11,4)(7,3.7)(2.5,1)
				\put(19,1.1){\vector(2,-1){0.5}}
				\put(2,1.1){\vector(-2,-1){0.5}}
				\put(10,3){\small $^{\tau}$}
			\end{picture}
		\end{equation}
		\vspace{0.1cm}
	\end{center}
	In this section, we shall obtain closed formulas for its dual canonical basis.
	
	Let $\ff$ be the algebra of type $A_2$. Denote
	\[
	\vartheta_{12}=\frac{v^{\frac{1}{2}}\vartheta_1\vartheta_2-v^{-\frac{1}{2}}\vartheta_2\vartheta_1}{v-v^{-1}},\quad \vartheta_{21}=\frac{v^{\frac{1}{2}}\vartheta_2\vartheta_1-v^{-\frac{1}{2}}\vartheta_1\vartheta_2}{v-v^{-1}}.
	\]
	We adopt the convention that $\vartheta_i^a=0$ for $a<0$; similar for $\vartheta_{12}$ and $\vartheta_{21}$. By \cite[Example 5.13]{BG17a}, the dual canonical basis $\mathbf{C}$ of $\ff$ is given by
	\begin{align}
		\label{eq:Lus-dCBA2}
		\mathbf{C} =\{b_\ba:=v^{\frac{1}{2}(a_2-a_1)(a_{12}-a_{21})}\vartheta_1^{a_1}\vartheta_2^{a_2}\vartheta_{12}^{a_{12}}\vartheta_{21}^{a_{21}}\mid \ba=(a_1,a_2,a_{12},a_{21})\in\N^4, a_1a_2=0\}.
	\end{align}
	
	A direct computation shows that
	\begin{align*}    \vartheta_1\vartheta_{12}&=v\vartheta_{12}\vartheta_1,\qquad \vartheta_2\vartheta_{12}=v^{-1}\vartheta_{12}\vartheta_2,\quad \vartheta_{12}\vartheta_{21}=\vartheta_{21}\vartheta_{12},
		\\    &\vartheta_1\vartheta_{21}=v^{-1}\vartheta_{21}\vartheta_1,\qquad \vartheta_2\vartheta_{21}=v\vartheta_{21}\vartheta_2.
	\end{align*}
	Recall the iHopf algebra $\tB^\imath_\tau$ on the Borel of type $A_2$. Note that $\K_1=v^{-1/2}h_2$, $\K_2=v^{-1/2}h_1$. We have
	\begin{align*}
		h_i*h_j&=v^{c_{i,\tau j}}h_ih_j=\begin{cases}
			v^{-1}h_i^2 &\text{ if }i=j,
			\\
			v^2h_ih_j &\text{ if }i\neq j.
		\end{cases}
		\\
		h_i*\vartheta_j&=h_i\vartheta_j
		\qquad
		\vartheta_j*h_i=v^{c_{i,\tau j}}\vartheta_jh_i,\qquad
		h_i*\vartheta_{12}=h_i\vartheta_{12},
		\\
		h_i*\vartheta_{21}&=h_i\vartheta_{21},
		\qquad
		\vartheta_{12}*h_i=\vartheta_{12}h_i,\qquad \vartheta_{21}*h_i=\vartheta_{21}h_i.
	\end{align*}

	By Theorem \ref{thm:dCB}, we denote by $\K_\alpha\diamond C_{\vartheta_1^{a_1}\vartheta_2^{a_2}\vartheta_{12}^{a_{12}}\vartheta_{21}^{a_{21}}}$ the dual canonical basis of $\tB^\imath_\tau$ corresponding to $\K_\alpha\diamond b_{(a_1,a_2,a_{12},a_{21})}$; see \eqref{eq:Lus-dCBA2}.  Moreover, 
	we have
	\begin{align*}    C_{\vartheta_{12}^a\vartheta_{21}^b}*C_{\vartheta_{12}^c\vartheta_{21}^d}=C_{\vartheta_{12}^c\vartheta_{21}^d}*C_{\vartheta_{12}^a\vartheta_{21}^b}.
	\end{align*}
	
	Clearly, we have $C_{\vartheta_1}=\vartheta_1$, $C_{\vartheta_2}=\vartheta_2$. We denote 
	\[x_{12}=\frac{v^{\frac{1}{2}}\vartheta_1*\vartheta_2-v^{-\frac{1}{2}}\vartheta_2*\vartheta_1}{v-v^{-1}},\quad x_{21}=\frac{v^{\frac{1}{2}}\vartheta_2*\vartheta_1-v^{-\frac{1}{2}}\vartheta_1*\vartheta_2}{v-v^{-1}}.\]
	One checks that 
	$$\vartheta_{12}=x_{12}+\K_2-v\K_1, \qquad \vartheta_{21}=x_{21}+\K_1-v\K_2,
	$$
	and then 
	\begin{align}
		C_{\vartheta_{12}}=x_{12}+\K_2-[2]\K_1=\vartheta_{12}-v^{-1}\K_1,
		\\
		C_{\vartheta_{21}}=x_{21}+\K_1-[2]\K_2=\vartheta_{21}-v^{-1}\K_2.
	\end{align}

	\begin{proposition}
		\label{prop:recursive-form1}
		For any $a\geq0,b\geq0$, we have
		\begin{align}
			\label{eq:recursive-form1}
			C_{\vartheta_{12}^{a+1}\vartheta_{21}^b}=C_{\vartheta_{12}}*C_{\vartheta_{12}^a\vartheta_{21}^b}-\K_1*\K_2*C_{\vartheta_{12}^a\vartheta_{21}^{b-1}} -\K_1*C_{\vartheta_{12}^{a-1}\vartheta_{21}^{b+1}},
			\\
			\label{eq:recursive-form2}
			C_{\vartheta_{21}^{a+1}\vartheta_{12}^b}=C_{\vartheta_{21}}*C_{\vartheta_{21}^a\vartheta_{12}^b}-\K_1*\K_2*C_{\vartheta_{21}^a\vartheta_{12}^{b-1}} -\K_2*C_{\vartheta_{21}^{a-1}\vartheta_{12}^{b+1}}.
		\end{align}
	\end{proposition}
	
	\begin{proof}
		It suffices to prove \eqref{eq:recursive-form1}, as the other one follows by symmetry. We proceed by induction on $a,b$. 
		We have 
		\begin{align*}
			\Delta(\vartheta_{12})&=\vartheta_{12}\otimes 1+h_1h_2\otimes \vartheta_{12}+v^{\frac12}h_2\vartheta_1\otimes \vartheta_2,
			\\
			\Delta(\vartheta_{12}^a\vartheta_{21}^b)=&\, \big(\vartheta_{12}\otimes 1+h_1h_2\otimes \vartheta_{12}+v^{\frac12}h_2\vartheta_1\otimes \vartheta_2\big)^a
			\\
			&\cdot \big(\vartheta_{21}\otimes 1+h_1h_2\otimes \vartheta_{21}+v^{\frac12}h_1\vartheta_2\otimes \vartheta_1\big)^b.
		\end{align*}
		So by definition, we have
		\begin{align*}    &\vartheta_{12}*\vartheta_{12}^a\vartheta_{21}^b
			\\
			&=\, \vartheta_{12}^{a+1}\vartheta_{21}^b+\sum_{t=0}^{a-1}\varphi(\vartheta_{21},\vartheta_{12})\vartheta_{12}^{t}h_1h_2\vartheta_{12}^{a-1-t}\vartheta_{21}^b+\sum_{t=0}^{b-1}\varphi(\vartheta_{21},\vartheta_{21})\vartheta_{12}^a\vartheta_{21}^th_1h_2\vartheta_{21}^{b-1-t}
			\\
			&\quad+\sum_{s=0}^{a-1}\sum_{t=0}^{b-1}\varphi(\vartheta_{21},\vartheta_2\vartheta_1)\vartheta_{12}^s(v^{\frac12}h_2\vartheta_1)\vartheta_{12}^{a-1-s}\vartheta_{21}^t(v^{\frac12}h_1\vartheta_2)\vartheta_{21}^{b-1-t}
			\\
			&\quad+\sum_{s=0}^{a-1}\varphi(v^{\frac12}h_1\vartheta_{2},\vartheta_2)\vartheta_2\vartheta_{12}^s(v^{\frac12}h_2\vartheta_1)\vartheta_{12}^{a-1-s}\vartheta_{21}^b
			\\
			&=\, \vartheta_{12}^{a+1}\vartheta_{21}^b-v^{-1}(v-v^{-1})\sum_{t=0}^{a-1} v^{-1-2t}\K_1*\K_2*\vartheta_{12}^{a-1}\vartheta_{21}^b
			\\
			&\quad+(v-v^{-1})\sum_{t=0}^{b-1}v^{-2a-2t-1}\K_1*\K_2*\vartheta_{12}^a\vartheta_{21}^{b-1}
			\\
			&\quad+v^{-1/2}(v-v^{-1})^2\sum_{s=0}^{a-1}\sum_{t=0}^{b-1}v^{-2s-2t-2}\K_1*\K_2*(v^{\frac12}\vartheta_{12}+v^{-1/2}\vartheta_{21})\vartheta_{12}^{a-1}\vartheta_{21}^{b-1}
			\\
			& \quad+\sum_{s=0}^{a-1}(v-v^{-1})v^{-2s-1}\K_1*\vartheta_{12}^{a-1}\vartheta_{21}^{b+1} +\sum_{s=0}^{a-1}(v-v^{-1})v^{-2s-2}\K_1*\vartheta_{12}^{a}\vartheta_{21}^b
			\\&\in \, \vartheta_{12}^{a+1}\vartheta_{21}^b+\delta_{a,0}\K_1*\K_2*\vartheta_{12}^a\vartheta_{21}^{b-1}+\delta(a>0) \K_1*\K_2*\vartheta_{12}^a\vartheta_{21}^{b-1}+\K_1*\vartheta_{12}^{a-1}\vartheta_{21}^{b+1}
			\\   &\quad+\sum_{(\beta,b)\in\N^\I\times\mathbf{C}}v^{-1}\Z[v^{-1}]\cdot \K_\beta\diamond \iota(b)
			\\
			&=\, \vartheta_{12}^{a+1}\vartheta_{21}^b+\K_1*\K_2*\vartheta_{12}^a\vartheta_{21}^{b-1}+\K_1*\vartheta_{12}^{a-1}\vartheta_{21}^{b+1}
			+\sum_{(\beta,b)\in\N^\I\times \mathbf{C}}v^{-1}\Z[v^{-1}]\cdot \K_\beta\diamond \iota(b),
		\end{align*}
		since $\varphi(\vartheta_{21},\vartheta_{12})=-v^{-1}(v-v^{-1}) =-\varphi(h_1\vartheta_2,\vartheta_2)$, $\varphi(\vartheta_{21},\vartheta_2\vartheta_1)=v^{-1/2}(v-v^{-1})^2$, and $\vartheta_{2}\vartheta_1=v^{\frac12}\vartheta_{21}+v^{-1/2}\vartheta_{12}$. Here 
		$$\delta(a>0)=\begin{cases}1& \text{ if }a>0,
			\\
			0& \text{ otherwise.}
		\end{cases}$$
		By Theorem \ref{thm:dCB}, we know
		$C_{\vartheta_{12}^a\vartheta_{21}^b}\in \vartheta_{12}^a\vartheta_{21}^b+\sum_{(\beta,b)\in\N^\I\times \mathbf{C}}v^{-1}\Z[v^{-1}]\cdot \K_\beta\diamond \iota(b)$. Then 
		\begin{align*}    &C_{\vartheta_{12}}*C_{\vartheta_{12}^a\vartheta_{21}^b}-\K_1*\K_2*C_{\vartheta_{12}^a\vartheta_{21}^{b-1}} -\K_1*C_{\vartheta_{12}^{a-1}\vartheta_{21}^{b+1}}
			\\
			&\in\vartheta_{12}^{a+1}\vartheta_{21}^b+\sum_{(\beta,b)\in\N^\I\times\mathbf{C}}v^{-1}\Z[v^{-1}]\cdot \K_\beta\diamond \iota(b),
		\end{align*}
		which is bar invariant. So by Theorem \ref{thm:dCB} again, we have 
		\begin{align*}    C_{\vartheta_{12}^{a+1}\vartheta_{21}^b}=C_{\vartheta_{12}}*C_{\vartheta_{12}^a\vartheta_{21}^b}-\K_1*\K_2*C_{\vartheta_{12}^a\vartheta_{21}^{b-1}} -\K_1*C_{\vartheta_{12}^{a-1}\vartheta_{21}^{b+1}}.
		\end{align*}
		This completes the proof.
	\end{proof}
	
	\begin{proposition}
		\label{prop:recursive-form2}
		For any $a,b,c\geq0$,
		we have 
		\begin{align}
			\label{eq:dCBA21}     C_{\vartheta_1^{a+1}\vartheta_{12}^b\vartheta_{21}^c} &=v^{(c-b)/2}C_{\vartheta_1}*C_{\vartheta_1^a\vartheta_{12}^b\vartheta_{21}^c}-\delta_{a,0}\K_1\diamond C_{\vartheta_1^{a+1}\vartheta_{12}^{b-1}\vartheta_{21}^c},
			\\     \label{eq:dCBA22}C_{\vartheta_2^{a+1}\vartheta_{12}^b\vartheta_{21}^c} &=v^{(b-c)/2}C_{\vartheta_2}*C_{\vartheta_2^a\vartheta_{12}^b\vartheta_{21}^c}-\delta_{a,0}\K_1\diamond C_{\vartheta_2^{a+1}\vartheta_{12}^{b}\vartheta_{21}^{c-1}}.
		\end{align}
		Moreover, the following identity holds:
		\begin{align}
			\label{eq:comm-rel1}    C_{\vartheta_{1}}*C_{\vartheta_{1}^a\vartheta_{12}^b\vartheta_{21}^c} =v^{b-c}C_{\vartheta_{1}^a\vartheta_{12}^b\vartheta_{21}^c}*C_{\vartheta_{1}}.
		\end{align}
	\end{proposition}
	
	\begin{proof}
		We prove \eqref{eq:dCBA21} and 
		\eqref{eq:comm-rel1} by induction on $a$. The proof for \eqref{eq:dCBA22} is skipped. 
		
		Let us first prove \eqref{eq:comm-rel1} for $a=0$. Clearly, we have $C_{\vartheta_1}*C_{\vartheta_{12}}=vC_{\vartheta_{12}}*C_{\vartheta_1}$ and $C_{\vartheta_1}*C_{\vartheta_{21}}=v^{-1}C_{\vartheta_{21}}*C_{\vartheta_1}$.
		By the recursive formulas in Proposition \ref{prop:recursive-form1}, we obtain by induction on $b+c$ that
		\begin{align}
			\label{eq:comm-relthetabc}
			C_{\vartheta_{1}}*C_{\vartheta_{12}^b\vartheta_{21}^c}=v^{b-c}C_{\vartheta_{12}^b\vartheta_{21}^c}*C_{\vartheta_{1}}.
		\end{align}

		Assume that 
		\begin{align}
			\label{eq:assump}
			C_{\vartheta_1}*C_{\vartheta_{1}^k\vartheta_{12}^b\vartheta_{21}^c}=v^{b-c}C_{\vartheta_{1}^k\vartheta_{12}^b\vartheta_{21}^c}*C_{\vartheta_{1}}, \quad \text{ for } k\leq a.
		\end{align}
		Now we compute 
		\begin{align*}   \vartheta_1*\vartheta_1^a\vartheta_{12}^b\vartheta_{21}^c &=\vartheta_1^{a+1}\vartheta_{12}^b\vartheta_{21}^c+(v-v^{-1})\sum_{t=0}^{b-1}\vartheta_1^a\vartheta_{12}^t(v^{\frac12}h_2\vartheta_1)\vartheta_{12}^{b-1-t}\vartheta_{21}^c
			\\   &=\vartheta_1^{a+1}\vartheta_{12}^b\vartheta_{21}^c+v^{\frac12}(v-v^{-1}) \sum_{t=0}^{b-1}v^{a-2t}h_2\vartheta_1^{a+1}\vartheta_{12}^{b-1}\vartheta_{21}^c
			\\   &=\vartheta_1^{a+1}\vartheta_{12}^b\vartheta_{21}^c+(v-v^{-1}) \sum_{t=0}^{b-1}v^{a-2t+1}\K_1* \vartheta_1^{a+1}\vartheta_{12}^{b-1}\vartheta_{21}^c
			\\    &=\vartheta_1^{a+1}\vartheta_{12}^b\vartheta_{21}^c+(v-v^{-1}) \sum_{t=0}^{b-1}v^{-2t-\frac{a+1}{2}}\K_1\diamond \vartheta_1^{a+1}\vartheta_{12}^{b-1}\vartheta_{21}^c.
		\end{align*}
		So 
		\begin{align*}
			&v^{(c-b)/2}    \vartheta_1*(v^{a(c-b)/2}\vartheta_1^a\vartheta_{12}^b\vartheta_{21}^c)
			\\
			&=(v^{(a+1)(c-b)/2}\vartheta_1^{a+1}\vartheta_{12}^b\vartheta_{21}^c)+(1-v^{-2})\sum_{t=0}^{b-1}v^{-a-2t}\K_1\diamond(v^{(a+1)(c-b+1)/2}\vartheta_1^{a+1}\vartheta_{12}^{b-1}\vartheta_{21}^c).
		\end{align*}
		By assumption, $v^{(c-b)/2}C_{\vartheta_1}*C_{\vartheta_1^a\vartheta_{12}^b\vartheta_{21}^c}$ is bar invariant. So similar to the proof of Proposition~ \ref{prop:recursive-form1}, we have
		\begin{align}
			\label{eq:recursivea}
			C_{\vartheta_1^{a+1}\vartheta_{12}^b\vartheta_{21}^c}&=v^{(c-b)/2}C_{\vartheta_1}*C_{\vartheta_1^a\vartheta_{12}^b\vartheta_{21}^c}-\delta_{a,0}\K_1\diamond C_{\vartheta_1^{a+1}\vartheta_{12}^{b-1}\vartheta_{21}^c}.
		\end{align}
		
		Finally, we prove that $C_{\vartheta_1}*C_{\vartheta_{1}^{a+1}\vartheta_{12}^b\vartheta_{21}^c}=v^{b-c}C_{\vartheta_{1}^{a+1}\vartheta_{12}^b\vartheta_{21}^c}*C_{\vartheta_{1}}$. If $a\neq0$, then it follows by \eqref{eq:recursivea} and the inductive assumption \eqref{eq:assump}. For $a=0$, it follows from \eqref{eq:recursivea} that 
		\begin{align*}
			C_{\vartheta_1\vartheta_{12}^b\vartheta_{21}^c}&=v^{(c-b)/2}C_{\vartheta_1}*C_{\vartheta_{12}^b\vartheta_{21}^c}-\K_1\diamond C_{\vartheta_1\vartheta_{12}^{b-1}\vartheta_{21}^c}.
		\end{align*}
		By \eqref{eq:comm-relthetabc} and induction on $b+c$, the desired identity holds. The proof is completed.
	\end{proof}
	
	Propositions \ref{prop:recursive-form1} and \ref{prop:recursive-form2} provide us explicit recursive formulas for the dual canonical basis of $\tUi_v(\mathfrak{sl}_3)$.


\end{document}